\documentclass[11pt]{article}
\usepackage{graphicx, color}
\usepackage{amssymb, amsmath, amsthm, mathdots}
\usepackage{youngtab}
\usepackage[all]{xy}
\usepackage{mathtools, hyperref}
\usepackage{colonequals}
\usepackage{dynkin-diagrams}
\usepackage{float}
\begin{document}
\newtheorem{defn0}{Definition}[section]
\newtheorem{prop0}[defn0]{Proposition}
\newtheorem{thm0}[defn0]{Theorem}
\newtheorem{lemma0}[defn0]{Lemma}
\newtheorem{coro0}[defn0]{Corollary}
\newtheorem{exa}[defn0]{Example}
\newtheorem{exe}[defn0]{Exercise}
\newtheorem{rem0}[defn0]{Remark}

\numberwithin{equation}{section}
\def\rig#1{\smash{ \mathop{\longrightarrow}
    \limits^{#1}}}
    \def\drig#1{\smash{ \mathop{\dashrightarrow}
    \limits^{#1}}}
\def\nwar#1{\nwarrow
   \rlap{$\vcenter{\hbox{$\scriptstyle#1$}}$}}
\def\near#1{\nearrow
   \rlap{$\vcenter{\hbox{$\scriptstyle#1$}}$}}
\def\sear#1{\searrow
   \rlap{$\vcenter{\hbox{$\scriptstyle#1$}}$}}
\def\swar#1{\swarrow
   \rlap{$\vcenter{\hbox{$\scriptstyle#1$}}$}}
\def\dow#1{\Big\downarrow
   \rlap{$\vcenter{\hbox{$\scriptstyle#1$}}$}}
\def\up#1{\Big\uparrow
   \rlap{$\vcenter{\hbox{$\scriptstyle#1$}}$}}
\def\lef#1{\smash{ \mathop{\longleftarrow}
    \limits^{#1}}}
\def\O{{\mathcal O}}
\def\L{{\mathcal L}}
\def\M{{\cal M}}
\def\N{{\cal N}}
\def\P#1{{\mathbb P}^#1}
\def\PP{{\mathbb P}}
\def\CC{{\mathbb C}}
\def\QQ{{\mathbb Q}}
\def\R{{\mathbb R}}
\def\Z{{\mathbb Z}}
\def\N{{\mathbb N}}
\def\H{{\mathbb H}}
\def\sym{{\mathrm{Sym}}}
\newcommand{\defref}[1]{Def.~\ref{#1}}
\newcommand{\propref}[1]{Prop.~\ref{#1}}
\newcommand{\thmref}[1]{Thm.~\ref{#1}}
\newcommand{\lemref}[1]{Lemma~\ref{#1}}
\newcommand{\corref}[1]{Cor.~\ref{#1}}
\newcommand{\exref}[1]{Example~\ref{#1}}
\newcommand{\secref}[1]{Section~\ref{#1}}
\newcommand{\magma}{\textsf{Magma}}
\def\geq{\geqslant}
\def\leq{\leqslant}
\def\ge{\geqslant}
\def\le{\leqslant}

\def\Pf{\mathrm{Pf}}
\def\End{\textrm{End}}
\def\diag{\textrm{diag}}
\def\rk{\textrm{Rk}}
\def\ad{\mathfrak{ad}}
\def\Ker{\textrm{ker}}

\newcommand{\qedd}{\hfill\figurebox[2mm]{\ }\medskip}
\newcommand{\codim}{\textrm{codim}}
\setcounter{tocdepth}{1}



\title{Some Classical Invariants, from Harmonic Quadruples to Triangle Groups}


\author{Giorgio Ottaviani\smallskip\\
\MakeLowercase{with an appendix by} Vincenzo Galgano}

\date{}
\maketitle
\begin{center}
	{\it To the memory of Alan Huckleberry, a mathematics enthusiast.}
	\end{center}

\begin{abstract}
\noindent These notes are an expanded version of the lectures held in Troms\o{}, in May 2025
at the Lie-St\o{}rmer Summer School : 
{\it Invariant Theory from classics to modern developments}, in the framework of {\em TiME} events. We emphasize the analogy
between binary quartics and ternary cubics (and subsequently modular forms) based on their harmonic and equianharmonic invariants. Triangle groups are presented in both the elliptic and the hyperbolic setting with their associated tilings.
The topics include the discussion of a short Hilbert paper on polynomials which are powers, that was proposed to the participants.
The appendix contains some exercises, with sketches of solutions, and a section devoted to Pfaffians edited by Vincenzo Galgano.

\end{abstract}

\tableofcontents

\section{Introduction: the harmonic quadruples and the cross-ratio}
One of the first interesting configurations in invariant theory concerns the harmonic quadruples. Their use goes back to Pythagorean music scale, and the name {\it harmonic} was suggested by Archytas of Taranto \cite[\S 1, 15]{Sev}, three centuries BC. In the diatonic scale, three chords of length $1$, $\frac{4}{5}$. $\frac 23$ make the fundamental chord C-E-G (do-mi-sol).

\begin{figure}[H]
\centering
\includegraphics[width=3cm]{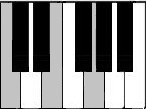}
\caption{\it chord C-E-G from an harmonic $4$-tuple.}\label{fig1}
\end{figure}

\noindent The second length is the harmonic mean of the other two. 

\begin{figure}[H]
\centering
\includegraphics[width=6cm]{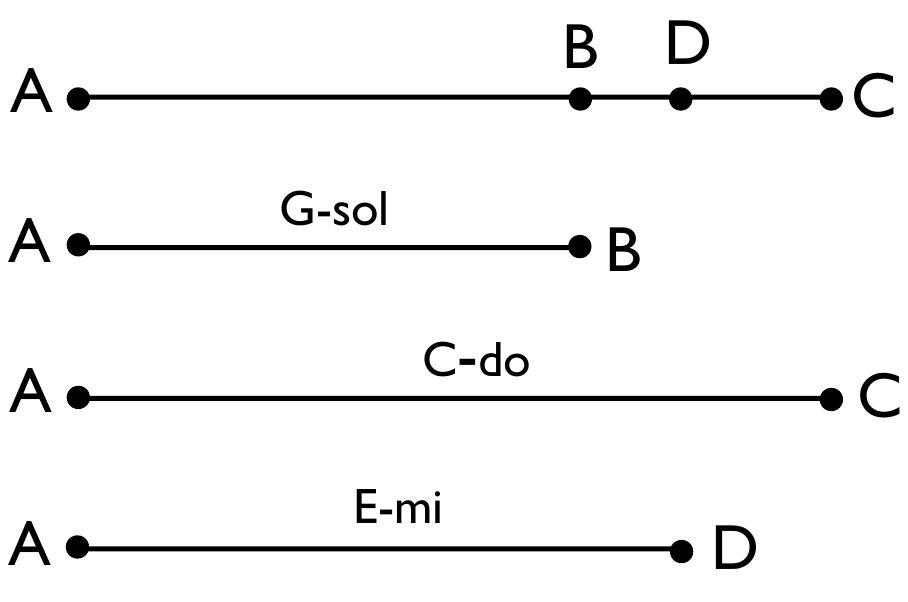}
\caption{\it AD is the harmonic mean between AB and AC. Recall that the sound frequencies are reciprocal to the lengths.}\label{fig2}
\end{figure}

\noindent The algebraic relation is

$$\frac{2}{AD}=\frac{1}{AB}+\frac{1}{AC}.$$ 

\noindent It gives $AB\cdot AC-AD\cdot AC=AB\cdot AD-AB\cdot AC$,
which is easily recognizable as equivalent to
\begin{equation}\label{eq:crossratio_harm}\frac{AC}{DC} = -\frac{AB}{DB}.\end{equation}

\noindent Equation (\ref{eq:crossratio_harm}) is the relation obtained by Apollonius in Prop. 37 of the third book of his famous treatise {\it Conics}.

\begin{figure}[H]
\centering{\includegraphics[width=6cm]{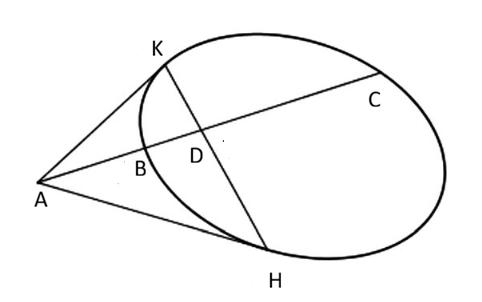}}
\caption{\it Apollonius construction: D is the fourth harmonic after A, B, C. The construction does not depend on the conic. The figure is from the introduction of \cite{DCG}.}\label{fig:apollonius}
\end{figure}

In figure \ref{fig:apollonius} consider $A$, $B$, $C$ as given and $D$ to be constructed.

\noindent The lines AH and AK are tangent to the (smooth) conic, so that HK is the polar line of A. The fact that the intersection point D does not depend on the conic is quite surprising at first glance and shows the depth reached by the mathematical culture of the Hellenistic period. The conic can be chosen as a circle, so that a ruler and compass construction can be carried out by a craftsman.
Ruler and compass were indeed the main computational tools of the antiquity. See the exercises in Appendix \ref{appendix:binary and ternary}.
I recommend the historical account of Projective Geometry {\it ``From Here to Infinity''} \space \cite{DCG} by Del Centina and Gimigliano for much more information.

\begin{figure}[H]
	\centering
\includegraphics[width=8cm]{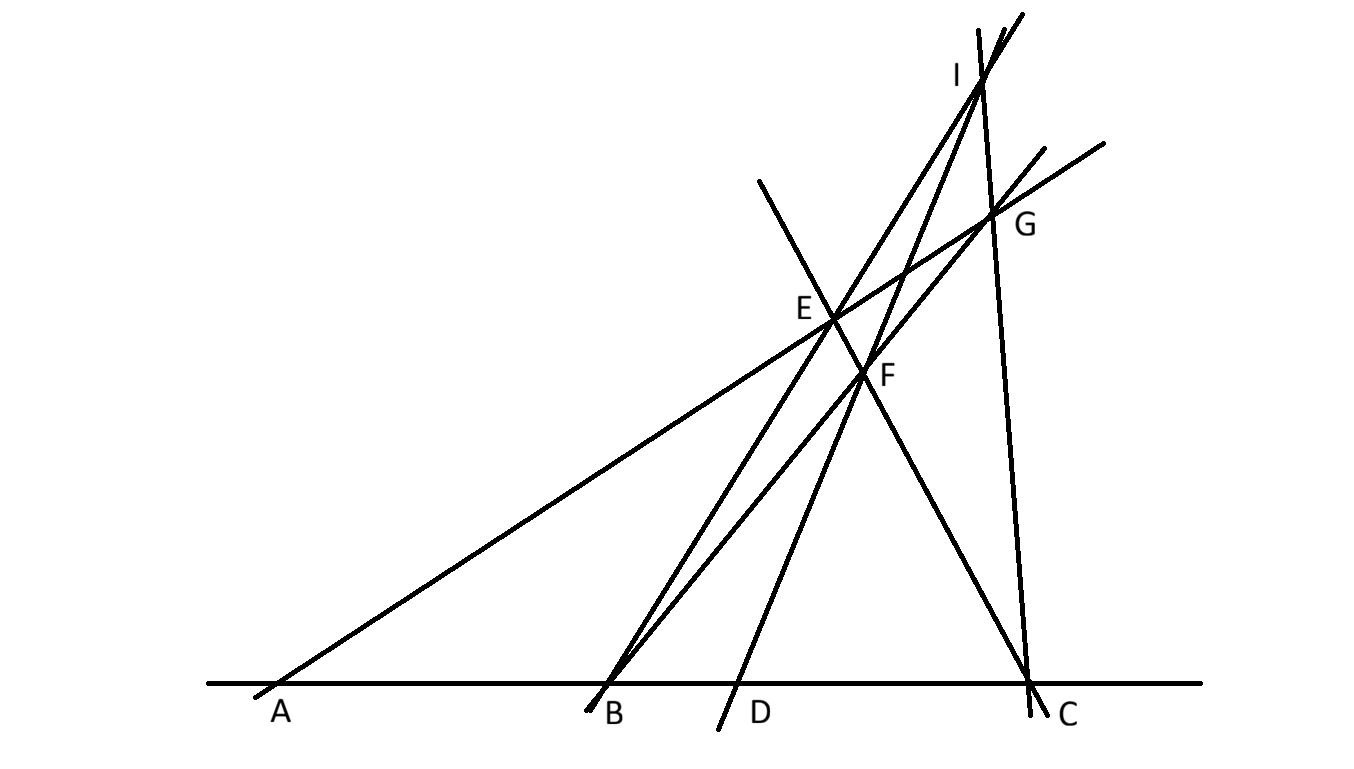}
\caption{\it The harmonic property of the complete quadrilateral. Given A, B, C collinear, choose E not collinear, choose F on CE. Then G, I and D are determined successively by the picture. The point D does not depend on the choices of E, F and it is the fourth harmonic after A, B, C.}\label{fig:complete_quadrilateral}
\end{figure}

Since the equation to find D is linear, it follows that the coordinates of $D$ lie in the same field of $A$, $B$ and $C$.
Square roots are not involved. Hence a construction without compass, with only the ruler, is possible.
This is the construction of the $4$-th harmonic with the complete quadrilateral, see figure \ref{fig:complete_quadrilateral}. 

This second construction is algebraically important but, in practice, it is beated by Apollonius' one, which requires fewer steps. One of the purposes of these TiME lectures is to keep the beauty of these constructions alive. I was encouraged by the participants of the Troms\o{} School to insert the two harmonic constructions in these notes.

It is worth remarking that the harmonic mean is a projective invariant, while the arithmetic mean is only an affine invariant. At the same time, it is wise to warn the reader that the fundamental chord C-E-G in music is a metric invariant. Applying an affine transformation, our ear still perceives a harmonic sound, corresponding to a musical key change, but a general projective transformation destroys the ratio between C and E. In other words, a musical piece is not projectively invariant. For example, the chord C-G-C (do-sol-do) still corresponds to a harmonic 4-tuple,
but it sounds completely different.

The harmonic quadruples of points in $\PP^1$ mark the birth of invariant theory.
Invariant theory provides a solid algebraic description of Geometry. The starting point is that geometric objects, like points, subspaces or subvarieties, are described in coordinates. Coordinate systems are not canonical but require some choices. Luckily, there is a group action on the set of these choices. The most relevant example is the group $GL(V)$, which describes exactly {\it the set of ordered basis of a vector space $V$}.  Hence this set comes equipped with a group structure, after a choice of an ordered basis has been made, to work as identity element. I learned this basic remark from Alan Huckleberry many years ago, at the beginning of my mathematical activity. We get the group action of $GL(V)$ on $V$. At projective level, the group $PGL(V)$ acts on $\PP(V)$, but it is often preferable to replace $PGL(V)$ with its universal cover $SL(V)$.
Given a group action, an invariant is a function on the geometric objects which does not change under the action.

Let's describe the harmonic quadruples in this setting.
Given $4$ ordered points on $\PP^1$, with coordinates respectively
$(a_0, a_1)$, $(b_0, b_1)$, $(c_0, c_1)$ $(d_0, d_1)$, denote 
$(ab):=\det\begin{pmatrix}a_0&a_1\\b_0&b_1\end{pmatrix}$, and so on.
Consider the action of $SL(2)$ on $\PP^1\times\PP^1\times\PP^1\times\PP^1$.
The three expressions
$$(ab)(cd)\qquad (ac)(bd)\qquad(ad)(bc)$$
are $SL(2)$-invariants, and they are linearly dependent, since their alternating sum is zero.
The invariant ring of four ordered points on $\PP^1$
is spanned by these expressions, see \cite[Theorem 35 in \S 5.3]{O12} for the precise statement.

Formula (\ref{eq:crossratio_harm}) is equivalent to

$$\frac{(ac)(bd)}{(ab)(cd)}= -1$$

or also to

$$\frac{(ac)(bd)}{(ad)(bc)}= \frac{1}{2}$$

where the left hand side is called the {\it cross-ratio} of the four points.
The cross-ratio of four distinct points may assume all the values $\lambda\in\CC$ with the exception of $\lambda=0,1,\infty$, values corresponding to quadruples with a double point.

The introduction of complex numbers into the history of quadruples led to the discovery of a new special class of quadruples, in some ways even more special than the harmonic quadruples.
They are less well-known, because they never have four real points. The quadruples belonging to this new class are called equianharmonic. In terms of  cross-ratio, they are determined by the following property, whose proof is a straightforward computation.

\begin{prop0}\label{prop:crossIJ} Let $A=(a_0, a_1)$, $B=(b_0, b_1)$, $C=(c_0, c_1)$ $D=(d_0, d_1)$ be four distinct points.
After permuting the $4$ points,  the cross-ratio $\frac{(ac)(bd)}{(ad)(bc)}=\lambda$
assumes  $6$ distinct values
 $$\lambda, \frac 1\lambda, 1-\lambda, \frac{1}{1-\lambda}, \frac{\lambda-1}{\lambda}, \frac{\lambda}{\lambda-1}$$
 
 unless
 \begin{enumerate}
 \item{} $\lambda\in\left\{-1, 2, \frac 12\right\}$, which can be taken as definition of harmonic $4$-tuple,
\item{} $\lambda\in\left\{e^{i\pi/3}, e^{-i\pi/3}\right\}$,  which can be taken as definition of equianharmonic $4$-tuple
\end{enumerate}
\end{prop0}

In closing this introduction, I feel it is important to remember that the tempered musical scale, based on $\sqrt[12]{2}$, used since the 19th century, has shifted the link between music and invariant theory from precision to approximation. The cross-ratio of the tempered C-E-G (do-mi-sol) is
$$\frac{(2^{-1/3}-2^{-7/12})}{2^{-1/3}(1-2^{-7/12})}=0,478\ldots$$
which approximates $\frac 12$ at $4.3\%$,  at the border of being perceptible to a trained ear.

\section{Binary quartics}\label{sec:binary_quartics}

Binary quartics are homogeneous polynomials of degree $4$ in two variables. They have $4$ roots which can be seen as four points in $\PP^1$, with multiplicity. Up to multiple scalars, they are determined by their roots,
if we work over $\CC$, which we will do from now on.
The space of binary quartics is easily identified with the symmetrized tensor product $\sym^4\CC^2$ . Up to multiple scalars, binary quartics can be thought as unordered quadruples of points on $\PP^1$,
so as quadruples up to the action of the symmetric group $\Sigma_4$. A polynomial $\sym^4\CC^2$ is called nonsingular if it has $4$ distinct roots.

\begin{prop0}\label{prop:IJ}
Let $Z=\left\{z_1,\ldots, z_4\right\}$ with $z_i\in\PP^1$ be the set of roots of a binary quartic $f=\sum_{i=0}^4a_i{4\choose i}x^{4-i}y^i$. 
\begin{itemize}
\item $Z$ is a harmonic $4$-tuple if and only if 
$J:= -a_{2}^{3}+2\,a_{1}a_{2}a_{3}-a_{0}a_{3}^{2}-a_{1}^{2}a_{4}+a_{0}a_{2}a_{4}=0$.

\item $Z$ is a equianharmonic $4$-tuple if and only if 
$I: =3\,a_{2}^{2}-4\,a_{1}a_{3}+a_{0}a_{4}=0$.

\item $Z$ is nonsingular if and only if $D:=I^3-27J^2\neq 0$. 
\end{itemize}
\end{prop0}

{\it First proof of Proposition \ref{prop:IJ}:}
As a computer algebra enthusiast, it is worth noting
that Proposition \ref{prop:IJ} has an immediate computational proof, based on Groebner basis eli\-mination, starting from the cross-ratio values achieved in Proposition \ref{prop:crossIJ}.
It is significant that the expressions for $I$ and $J$ appear automatically on the screen.
Here is the M2 code\cite{GS} with a few comments.

\begin{verbatim}
KK=toField(QQ[q]/ideal(q^2-q+1))---- q=exp(2*i*pi/3)
R=KK[a_0..a_4,x,y,t_0..t_3]
f=sum(5,i->a_i*binomial(4,i)*x^(4-i)*y^i)---binary quartic
x4=symmetricPower(4,matrix{{x,y}})
H=minors(2,contract(x4,product(4,i->(x-t_i*y)))||contract(x4,f))--coeffs versus roots
J1=H+ideal(2*(t_2-t_0)*(t_3-t_1)-(t_3-t_0)*(t_2-t_1))--cross-ratio=1/2
I1=H+ideal(q*(t_2-t_0)*(t_3-t_1)-(t_3-t_0)*(t_2-t_1))---cross-ratio=1/q
J=-(gens eliminate({t_0,t_1,t_2,t_3},J1))_(0,0)---deg 3--harmonic
I=3*(gens eliminate({t_0,t_1,t_2,t_3},I1))_(0,0)---deg 2---equianharmonic
\end{verbatim} 

\noindent The well known case of the discriminant can be treated in the same way. 
 \qed
\vskip 0.3cm

It is also worth providing a more theoretical understanding. Note indeed that (cf. \cite[Sec.\ 9.2]{Dolg}) 
\begin{equation}\label{eq:J}
J=\det\begin{pmatrix}a_0&a_1&a_2\\
a_1&a_2&a_3\\
a_2&a_3&a_4\end{pmatrix}.
\end{equation} 
Indeed, the $3\times 3$ matrix in (\ref{eq:J}) is ($12$ times) the matrix of the contraction operator
\begin{equation}\label{eq:midcat}\begin{array}{ccc}\sym^2V^\vee&\to&\sym^2V\\
\partial&\mapsto&\partial f\end{array}\end{equation}
in the standard bases $\{\partial_x^2,\partial_x\partial_y,\partial_y^2\}$ and $\{x^2, 2xy, y^2\}$.

\begin{lemma0}\label{lem:quartic}
$$\mathrm{rk}\begin{pmatrix}a_0&a_1&a_2\\
a_1&a_2&a_3\\
a_2&a_3&a_4\end{pmatrix}\le 1\Longleftrightarrow \exists(\lambda,\mu)\in\PP^1\textrm{\ such that\ }
a_i=\lambda^{4-i}\mu^{i}$$
\end{lemma0}

\begin{proof} The implication $\Longleftarrow $ is trivial since, by substituting the values of $a_i$, any two rows become proportional.
In order to prove the implication $\Longrightarrow $,
without loss of generality we can assume $a_0=1$, $a_1=t$. Then the dependency of the first two rows give
that the second row is $t$ times the first row; this fact implies
$a_2=t^2$ and $a_3=t^3$. In the same way $a_4=t^4$, looking at the third row.
\end{proof}
\vskip 0.3cm

{\it Second proof of Proposition \ref{prop:IJ}:}
The determinantal expression of $J$ allows us to interpret  the hypersurface $V(J)$ as the secant variety of the quartic rational normal curve $C_4$, the latter being the locus of quartics which are $4$-powers of a linear form. 
Indeed, denote $C_f:=\begin{pmatrix}a_0&a_1&a_2\\
a_1&a_2&a_3\\
a_2&a_3&a_4\end{pmatrix}$. Lemma \ref{lem:quartic} is equivalent to the statement that
$$\mathrm{rk} C_f\le 1\Longleftrightarrow f=(\lambda x_0+\mu x_1)^4.$$
Hence, for two linear forms $\ell_1$, $\ell_2$, we have
\begin{equation}\label{eq:rkJ}\mathrm{rk}\left( C_{\ell_1^4+\ell_2^4}\right)=\mathrm{rk}\left( C_{\ell_1^4}+C_{\ell_2^4}\right)\le 
\mathrm{rk} C_{\ell_1^4}+\mathrm{rk} C_{\ell_2^4}=1+1=2 .
\end{equation}

The secant variety $\sigma_2(C_4)$ is the Zariski closure of the locus of the secant lines
$\{\ell_1^4+\ell_2^4 \mid \ell_i\textrm{\ linear forms}\}$. This locus includes also the tangent variety, i.e. the union of the tangent lines,
which is the $SL(2)$-orbit of $x^3y$.

The inequality \eqref{eq:rkJ} implies that $\sigma_2(C_4)\subseteq V(J)$. Since $V(J)$ is a $3$-dimensional irreducible variety, it follows that the equality holds (see Exercise \ref{exe:sec2 and J invariant} for a different proof). Indeed $V(J)$ is the union of three $SL(2)$-orbits, which are respectively $C_4$,
the tangent variety to $C_4$ and the secant locus $\{\ell_1^4+\ell_2^4 \mid \ell_i \in \mathbb P^1\}$.

Hence a binary quartic is harmonic if and only if it is the sum of two $4$-powers, or a limit of these.
The typical representative is $f=x^4-y^4$, whose roots are the four $4$th roots of unity.
The binary quartic defined from our C-E-G (do-mi-sol) in the first page is
\begin{align*}f=x(x-y)(5x-4y)(3x-2y)=15\,x^{4}-37\,x^{3}y+30\,x^{2}y^{2}-8\,x\,y^{3}=\\
\frac{i}{4}\left(\frac{-5+i}{2}x+2y\right)^4-
\frac{i}{4}\left(\frac{-5-i}{2}x+2y\right)^4\end{align*}

The decomposition of $f$ as a sum of two $4$th powers has been found by Sylvester's technique (for a nice modern account see \cite{BGI}),
since it is straightforward to compute the kernel of (\ref{eq:midcat}) and check that
$$2\left(2\partial_x+\frac{5-i}{2}\partial_y\right)\left(2\partial_x+\frac{5+i}{2}\partial_y\right)f=\left(8\partial_x^2+20\partial_x\partial_y+13\partial_y^2\right)f=0$$

\begin{figure}[H]
\centering{\includegraphics[width=6cm]{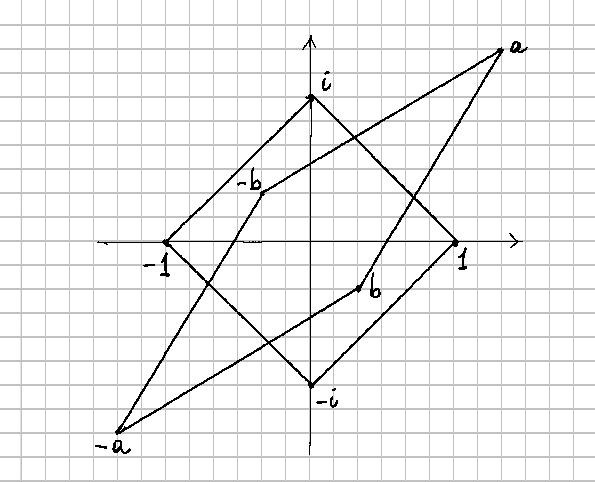}}
\caption{\it The four vertices of the square make a harmonic $4$-tuple, where $J=0$. The four vertices of the other rhombus make a equianharmonic $4$-tuple, where $I=0$. Here $a=\frac{1+i}{\sqrt{3}-1}$, $b=\frac{1-i}{\sqrt{3}+1}$.}\label{fig:IJ}
\end{figure}

The interpretation of $I$ is more delicate. A typical representative is $f=x(x^3-y^3)$, whose
 roots are the three $3$rd roots of unity with their midcenter. In the stereographic projection
 they correspond to the four vertices of a tetrahedron on the sphere, which is not regular. Klein choosed a different representative, namely
 
 $$\Psi=x^4-2\sqrt{-3}x^2y^2+y^4$$
 
 \noindent which has roots $\pm a, \pm b$ where $a=\frac{1+i}{\sqrt{3}-1}$, $b=\frac{1-i}{\sqrt{3}+1}$.
  In the stereographic projection
 they correspond to the four vertices of a regular tetrahedron on the sphere, see Figure \ref{fig:IJ}.
 
In \cite{O12} I gave three different proofs of the anharmonic case of Proposition \ref{prop:IJ}, one based on Reynolds operators,
a second one on Lie algebra action and a third one by using transvectants, see \cite[paragraph before Remark 1 in \S 1.2]{O12}. We will review the main property of transvectants
in Theorem \ref{thm:transvectant}
and the corresponding expressions for $I$, $J$ in (\ref{eq:IJtrans}).
\qed
\vskip 0.3cm 

\begin{prop0}
Let $f\in\sym^4\CC^2$ nonsingular. Then

\[\# Isotropy(f)=\begin{cases}
4 & IJ\neq 0 \\
8 & J=0 \\
12 & I=0 
\end{cases} \ .\]
\end{prop0}

The stabilizers are respectively the groups
$D_2=\Z_2\times\Z_2$, $D_4$ and $A_4$ (two dihedral and one alternate).

Correspondingly, the $SL(2)$-orbits in $\PP^4=\PP(\sym^4\CC^2)$ which are hypersurfaces
have degree $6$, and their equation has the form $I^3-kJ^2=0$ for some $k\in\CC$, unless the two special orbits
$V(I)$ of degree $2$ and $V(J)$ of degree $3$. In Remark \ref{rem:singdisc} we will see that
the discriminant is not isomorphic to the general orbit of degree $6$, since it has a larger singular locus.

We refer to K. DeVleming's lectures in this volume \cite{DeV} for the definition of semistable points.
The variety of not semistable points is called the nullcone and it coincides with the zero locus of all the invariants. 
For quadruples of points, the nullcone consists of quadruples where at least $3$ points coincide. For binary quartics the nullcone
is $V(I, J)$ (see Exercise \ref{exercise: I,J,Hess}). The nullcone coincides with the locus of quartics with a root of multiplicity $\ge 3$
and it consists of two orbits, namely $C_4$ and its tangent variety.
Note that the quartics which are squares (having two roots of multiplicity $2$) are semistable unless
they have a root of multiplicity $4$. We will meet the variety of squares again in Section \ref{sec:Hilbert1886}, when we discuss Hilbert's paper \cite{Hil86}.

\begin{prop0} \cite[Chap.\ II.2, Theorem 1, p.\ 23]{DO} For $p=(p_1,\ldots, p_4)\in \PP^1\times\PP^1\times\PP^1\times\PP^1 $,
 \begin{itemize}
 \item $p$ is semistable if and only if
at most $2$ among $p_1,\ldots, p_4$ coincide.
\item $p$ is stable if and only if the points $p_1,\ldots, p_4$ are distinct.
 \end{itemize}
 \end{prop0}

\noindent The map

$$\begin{array}{ccc}\bigg(\PP^1\times\PP^1\times\PP^1\times\PP^1 \bigg)^{ss}&\to&\PP^1\\
(a,b,c,d)&\mapsto&\left( (ac)(bd), (ad)(bc) \right)\end{array}$$

\noindent is the quotient map on the moduli space.

\noindent The symmetrized version for binary quartics is the following quotient map

$$\begin{array}{ccc}\left(\PP^4 \right)^{ss}=\PP^4\setminus V(I,J)&\to&\PP^1\\
f&\mapsto&\left( I^3, J^2 \right)\end{array}$$

\begin{rem0}\label{rem:singdisc} Among the semistable quartics there is also the form $x^2y(x+y)$, with a double point.
Its orbit is the discriminant $D=I^3-27J^2$, which contains in the singular locus the null cone
but also the variety of squares. The general orbit contains in the singular locus only the null cone.
\end{rem0}

A systematic construction of invariants is possible through transvectants.

\begin{thm0}\label{thm:transvectant}
Let $\dim U=2$. Let $\min(d, e)\ge n$. There is a unique bilinear $SL(2)$-invariant map

$$\sym^dU\times \sym^eU\to\sym^{d+e-2n}U$$
such that 
$$(u^d,v^e)\mapsto u^{d-n}v^{e-n}(u\wedge v)^n.$$
If $f\in S^dU$, $g\in S^eU$ , its expression is 

$$(f,g)_n=\frac{n!}{d!}\frac{n!}{e!}\sum_{i=0}^n(-1)^i{n\choose i}\frac{\partial f}{\partial x^{n-i}\partial y^i}
\frac{\partial g}{\partial x^{i}\partial y^{n-i}} \ , $$
which is called the $n$-th transvectant of $f$ and $g$.
\end{thm0}
\begin{proof}
The uniqueness is clear since any polynomial can be expressed as a sum of powers of linear forms (this statement
is equivalent to the fact that the linear span of the rational normal curve is the full ambient space).
To show the existence, we compute first 
\begin{align*}
(u^d, v^e)_n & =\frac{n!}{d!}\frac{n!}{e!}\sum_{i=0}^n(-1)^i{n\choose i}\frac{\partial u^d}{\partial x^{n-i}\partial y^i}
\frac{\partial v^e}{\partial x^{i}\partial y^{n-i}}= \\
& = \sum_{i=0}^n(-1)^i{n\choose i}u^{d-n}v^{e-n}(u_0v_1)^{n-i}(u_1v_0)^{i} =u^{d-n}v^{e-n}(u_0v_1-u_1v_0)^n \ .
\end{align*}
This computation makes easy to check the $SL(2)$-invariance, since by linearity it is sufficient to check it
when $f=u^d$, $g=v^e$. In this case, for any $h\in SL(2)$ we have

$$(h\cdot u^d,h\cdot v^e)\mapsto (h\cdot u)^{d-n}(h\cdot v)^{e-n}(u\wedge v)^n=
h\cdot(u^d, v^e)_n$$
The proof is complete.
\end{proof}
\vskip 0.3cm

 Note that $(f,g)_1$
is the Jacobian, while $(f,f)_2$ is the Hessian. A classical Theorem of Gordan states that all invariants of binary forms can be expressed
as an iteration of transvectants. 

For $f\in \sym^4\CC^2$ we can express the invariants $I$ and $J$ 
in terms of transvectants. Indeed, it is easy to check that
\begin{equation}\label{eq:IJtrans}
I=(f,f)_4\qquad J=(f,(f,f)_2)_4.\end{equation}

\begin{thm0}\label{thm:inv_binary_quartics}
The invariant ring $\CC[\sym^4\CC^2]^{SL(2)}$ is the polynomial ring $\CC[I, J]$.
\end{thm0}

For a proof we refer to \cite[Theor. 25]{O12}. Note that the invariant ring for $SL(2)$ coincides with the semi-invariant ring for $GL(2)$
(for semi-invariants see \eqref{eq:semi-inv} ).

\begin{rem0}
We have seen that, when projected on the sphere $S^2=\PP^1(\CC)$, the four roots of a equianharmonic quartic are the four vertices of a regular tetrahedron, while the four roots of a harmonic quartic are coplanar and are the vertices of a square.
This beautiful description has to be understood with a grain of salt.
Up to $SL(2)$-action the four roots of a harmonic quartic may loose the metric configuration of the square,
but they remain coplanar. The same holds for the four roots of a equianharmonic quartic, that may loose the metric configuration of the tetrahedron.
\end{rem0}

\begin{rem0}\label{rem:Ipfaff}
The invariant $I$ can be written as a Pfaffian: see Appendix \ref{appendix:Pfaffians} for more on Pfaffians. Indeed, up to scalar multiples,

$$I=\Pf\left(\!\begin{array}{cccc}
      0&a_{0}&2\,a_{1}&3\,a_{2}\\
      -a_{0}&0&a_{2}&2\,a_{3}\\
      -2\,a_{1}&-a_{2}&0&a_{4}\\
      -3\,a_{2}&-2\,a_{3}&-a_{4}&0
      \end{array}\!\right). $$
      This $4\times 4$ matrix corresponds to the composed linear map (for $\dim U=2$)
      
      $$\sym^3U^\vee \ \rig{a} \ \sym^3U \ \rig{b} \  \sym^3U$$
      where $a$ is induced as third symmetric power of the map
      
      $$\begin{array}{ccc}U^\vee&\to& U\\
      \partial_v&\mapsto&\partial_v(x\wedge y)\end{array}$$
      and $b(g)=(g,f)_2$.
\end{rem0}

The general binary quartic is $SL(2)$-equivalent (see Exercise \ref{exercise: generic binary quartic}) to a member of the pencil
\begin{equation}\label{eq:pencil_quartic}x^4+6tx^2y^2+y^4.\end{equation}
The Hessian of any quartic in the pencil  (\ref{eq:pencil_quartic}) is still a member of the pencil \cite[Remark 3.4]{CO}.
For this pencil of quartics, the two invariants have the simple form $I=1+3t^2$ and $J=t(1-t)(1+t)$, and the discriminant is a scalar multiple of $(3t-1)^2(3t+1)^2$.

The following proposition is an amusing consequence of Theorem \ref{thm:inv_binary_quartics}.
\begin{prop0}
Consider the Hessian map
$$H: \PP^4 \ \drig{} \ \PP^4 .$$
Then the Jacobian of the Hessian map is a scalar multiple of $I\cdot J$.
\end{prop0}
\begin{proof}
The Jacobian of the Hessian map is an invariant of degree $5$
(it is the determinant of a $5\times 5$ matrix with linear entries).
The statement follows from Theorem \ref{thm:inv_binary_quartics}.
The Jacobian is actually nonzero since the Hessian map is $2:1$ and it is dominant.
\end{proof}

\begin{prop0}[cf. Exercise \ref{exercise: I,J,Hess}]\label{prop:hessian polar of J}
The Hessian map is the polar map of $J$. It is generically $2:1$, and its branch locus coincides with the hypersurface $V(J)$. Moreover, given $f\in \sym^4\CC^2$, then $J(f)=0$ if and only if $H(f)$ is a square, while $I(f)=0$ if and only if $H(H(f))$ is proportional to $f$.
\end{prop0}

Enriques and Fano devoted the last part of their pioneering work \cite{EF} to special $SL(2)$-orbits. We reproduce here the following starting point, whose proof
is not difficult.

\begin{thm0}[Enriques-Fano \cite{EF}]\label{thm:EF}
The $SL(2)$-orbit of a binary form $f$ is a $3$-fold if and only if its stabilizer is finite.
This happens if $d\ge 3$ and $f$ has $d$ distinct roots. In this case, as a variety in $\PP(\sym^d\CC^2)$, the degree of the closure is
$$\deg \overline{\left(SL(2)\cdot [f]\right)}=\frac{d(d-1)(d-2)}{\# Stab(f)}$$
where $Stab(f)$ is a finite group in $PGL(2)$.
\end{thm0}

At the end of their paper Enriques and Fano list four smooth $3$-folds that have a $SL(2)$-action with a dense orbit, which come from polyhedral groups. This is likely the first geometric study of orbits coming from Klein's classification of polyhedral groups (that we will resume in Section \ref{sec:ADE}) .
Equations are obtained for all these cases from the covariant $(f, f)_4$, which is the anharmonic invariant in the case of binary quartics. The tetrahedral case is treated in
\cite[\S 25]{EF}, the octahedral case in \cite[\S 26]{EF} and finally the icosahedral case in \cite[\S 27]{EF}.
The icosahedral case is a tour de force,  it is five pages long, where the $17$ equations of the covariant $(f, f)_4$ are written explicitly. From these equations
it can be proved that the examples are smooth, since it is enough to compute the rank of the Jacobian at the representative $x^d$ of the closed orbit.

Only one century later, Aluffi and Faber prove the remarkable result that these examples (and their powers) are the only ones
which are smooth $3$-folds, namely:

\begin{thm0}[Aluffi-Faber \cite{AF}, Prop. 4.2]\label{thm:AF}
The smooth $3$-dimensional $SL(2)$-orbit closures in some $\PP(\sym^d\CC^2)$ are the following four classes:
\begin{enumerate}
\item{} the orbit closure of $(x^3 + y^3)^k$, $k\ge 1$, with stabilizer $D_3 = S_3\subset PSL(2)$;
\item{}  the orbit closure of $(x^4 + xy^3)^k$, $k\ge 1$,  (power of a anharmonic quartic), with stabilizer the tetrahedral group $A_4\subset PSL(2)$;
\item{}  the orbit closure of $(x^5y - xy^5)^k$, $k\ge 1$, with stabilizer the octahedral group $S_4\subset PSL(2)$;
\item{}  the orbit closure of $(x^{11}y + 11x^6y^6 - xy^{11})^k$, $k\ge 1$, with stabilizer the icosahedral group $A_5\subset PSL(2)$.
\end{enumerate}

\noindent The four $3$-folds obtained are respectively (for $k=1$)

\begin{enumerate}\label{Fanolist}
\item{} $\PP^3$;
\item{}  $Q_3$, the smooth quadric $3$-fold;
\item{} the quintic Del Pezzo $3$-fold $V_5$ in $\PP^6$;
\item{}  the prime Fano $3$-fold $U_{22}\subset\PP^{12}$ of genus $12$.
\end{enumerate}

The orbits for $k\ge 1$ are respectively isomorphic to these ones, embedded by the line bundle $\O(k)$.
\end{thm0}

In \cite{AF} one can find also a modern proof of Theorem \ref{thm:EF}. This story has an interesting sequel,  since all the four varieties in Theorem \ref{thm:AF} are examples of smooth Fano $3$-folds, meaning that their anticanonical bundle $-K$ is ample. This fact was recognized by Mukai and Umemura in \cite{MU},  who actually cited Enriques-Fano's paper \cite{EF}.
Since prime Fano $3$-folds of genus $12$ were omitted in \cite[chap. V, \S 7]{Roth}, sometimes it has even been claimed that Fano was unaware of the $U_{22}$ variety, which is paradoxical since $U_{22}$ was introduced by Fano himself, together with Enriques.

In \cite[Remark 6.5]{MU} a nice feature of the variety $X=U_{22}\subset\PP(\sym^{12}\CC^2)=\PP^{12}$ is understood. They realize that $H^0(-K_X)=\sym^{12}\CC^{2}\oplus\sym^0\CC^{2}$ as $SL(2)$-module, it has dimension $14$, so that the anticanonical embedding lives in $\PP^{13}$, where a hyperplane $H=\PP\sym^{12}\CC^{2}$ and a point $P=\PP\sym^{0}\CC^{2}$ are ``canonically'' defined.
The plane sections of the anticanonical embedding are canonical curves.
The anticanonical embedding is the closure of the  $SL(2)$-orbit of the pair $(f,P)\in \sym^{12}\CC^{2}\oplus\sym^0\CC^{2}$, where $f=x^{11}y + 11x^6y^6 - xy^{11}$
from item 4. of Theorem \ref{thm:AF}.
The intersection $H\cap X$ is a singular surface, isomorphic to the tangent variety of the rational normal curve of degree $12$. The projection from $P$ to $H$ of the anticanonical embedding of $X$ is isomorphic to $X$ itself and
 it coincides with the closure of the $SL(2)$-orbit of $f$. Note that the plane sections of $U_{22}$ are just projections of canonical curves. The saturated ideal is generated by $31$ quadrics. It is a pleasant exercise to get the $31$ quadrics by saturating with M2 \cite{GS} the ideal of the $17$ quadrics provided by Enriques and Fano.
 One subtle property of $U_{22}$ is that the  scheme of lined containes in $U_{22}$ is nonreduced and indeed 
corresponds to the nonreduced polynomial $(x^2+y^2+z^2)^2$, see \cite{Muk, Sch}, \cite[\S 5,\S 8.2]{Fae}.

\begin{rem0}
The Betti numbers of $U_{22}$ in its anticanonical embedding in $\PP^{13}$
are the same of its hyperplane section, which is the tangent surface.
They are
\[ \begin{matrix}
         & 0 & 1 & 2 & 3 & 4 & 5 & 6 & 7 & 8 & 9 & 10\\
      \text{total:}
         & 1 & 45 & 231 & 550 & 693 & 660 & 693 & 550 & 231 & 45 & 1\\
      0: & 1 & . & . & . & . & . & . & . & . & . & .\\
      1: & . & 45 & 231 & 550 & 693 & 330 & . & . & . & . & .\\
      2: & . & . & . & . & . & 330 & 693 & 550 & 231 & 45 & .\\
      3: & . & . & . & . & . & . & . & . & . & . & 1
      \end{matrix} \]     
These Betti numbers are computed in \cite{AFPRW}. This computation for the rational normal curve
of degree $d$ was the key tool for providing a proof of the Green conjecture for the general canonical curve,
simpler compared to Voisin's original proof.
For a nice account of this story see \cite{EL} and \cite{RS}.\\
It is worth to see the Betti numbers of $U_{22}\subset\PP^{12}$.
Using M2 \cite{GS}, they are
\[ \begin{matrix}
  & 0 & 1 & 2 & 3 & 4 & 5 & 6 & 7 & 8 & 9 & 10 & 11 & 12\\
      \text{total:} & 1 & 31 & 140 & 504 & 1\,298 & 2\,046 & 2\,079 & 1\,474 & 759 & 288 & 78 & 13 & 1\\
      0: & 1 & . & . & . & . & . & . & . & . & . & . & . & .\\
      1: & . & 31 & 109 & 108 & 11 & . & . & . & . & . & . & . & .\\
      2: & . & . & 31 & 396 & 1\,287 & 2\,046 & 2\,079 & 1\,474 & 759 & 287 & 78 & 13 & 1\\
      3: & . & . & . & . & . & . & . & . & . & 1 & . & . & .
      \end{matrix}\]
A look at the previous table, in particular at the column labelled with 9, shows that the icosahedron has still some subtlety to be revealed.
\end{rem0}

\begin{rem0}
The items 3. and 4. of Theorem \ref{thm:AF} can be found as varieties of sum of powers.
Precisely, the variety $V_5$  parametrizes the three summands $\sum_{i=0}^2 l_i^2$
which express a smooth conic $q$ \cite[Prop. 10]{Muk}.
The prime Fano $3$-fold $U_{22}$ parametrizes the summands in the squared expression $\sum_{i=0}^5l_i^4=q^2$ \cite[Theorem 14]{Muk}.
This last result has a real counterpart, when the six lines $l_i$ are real, they correspond exactly to the six pairs of vertices of a regular icosahedron, see \cite[Theorem 9.13]{Rez}.
This shows again the link of the icosahedron with $U_{22}$.

The variety $U_{22}$ has a $6$-dimensional family of deformations as Fano $3$-fold $X_{22}$
of genus $12$, which correspond to the summands  $\sum_{i=0}^5l_i^4=g$ where now $g$ is a general quartic plane curve,
parametrized indeed by the $6$-dimensional variety $M_3$ of curves of genus $3$ \cite[Theorem 11]{Muk}.

\end{rem0}

\begin{rem0}\label{rem:lastfano}
There is a different family of Fano $3$-folds of degree $22$ in $\PP^{13}$,
that was described by Fano in his last paper. This has index $2$  and has been revised by Andreatta and Pignatelli
in \cite{AP}, who called this example the ``Fano's last Fano''. 

\end{rem0}

\section{Ternary cubics}\label{sec:ternary_cubics}
Consider ternary cubics $f\in\sym^3\CC^3$. In this section we give for granted the well known correspondence between elliptic curves and quotients $\CC/\Lambda$
where $\Lambda$ is a lattice, see for example \cite[IV \S 4]{Hart}. This is essentially Abel's Theorem for genus one curves. Under the action of $SL(3)$, the behaviour of ternary cubics and binary quartics is quite similar.

The reason comes from the following construction. 

\begin{figure}[H]
	\centering
\includegraphics[width=6cm]{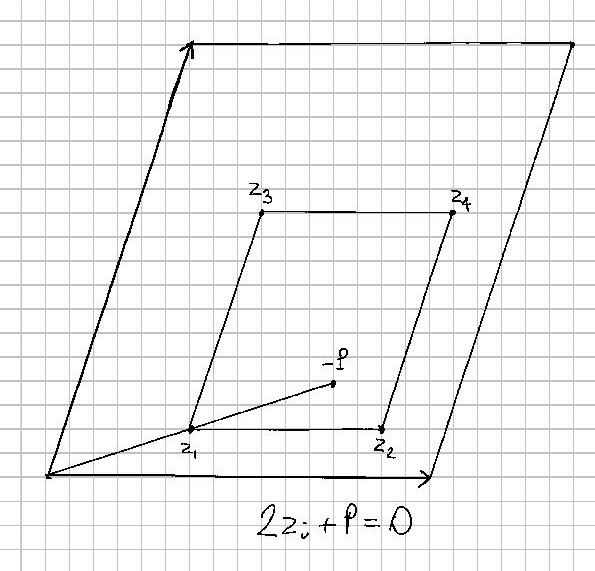}
\caption{\it The points $z_1,\ldots z_4$ are the solutions of the equation $2z+p=0$ in the lattice
$\CC/\Lambda$, where the generators of $\Lambda$ are the two independent sides of the parallelogram.}\label{fig:2zp}
\end{figure}

Let $f\in\sym^3\CC^3$ and pick $p\in V(f)\subset\PP^2$.
The projection $\pi_p\colon V (f)\to\PP^1$ is $2:1$ and by Riemann-Hurwitz has exactly four ramification points
$z_i\in V(f)$ and four branch points $\pi_p(z_i)\in\PP^1$.
We give for granted here the structure of complex torus on any elliptic curve,
which is Abel's Theorem for elliptic curves. Hence, after a flex has been fixed as origin,
three points on an elliptic curve are collinear if and only if their sum is zero in the group structure.
In this setting the four ramification points of $\pi_p$ are the four solutions $z$ of the equation
$2z+p=0$, clear from the lattice structure, see Figure \ref{fig:2zp}.

\begin{thm0}[Salmon]
The $SL(2)$-class of the four (unordered) branch points does not depend on $p\in V(f)$. 
\end{thm0}
\begin{proof}
This is essentially \cite[Lemma IV 4.4]{Hart}. We provide an alternative proof.

The four branch points give a map on the invariant ring of binary quartics.
After Theorem \ref{thm:inv_binary_quartics} we get a map 
\[ V(f) \ \rig{I^3,J^2} \ \PP^1. \] 
Since the four points are always distinct (see Figure \ref{fig:2zp}) , the image of the map sits in $\PP^1\setminus (27,1)$, in other words it avoids the discriminant.
Since the image is contained in a affine variety, or in a simpler way it is a global regular function, it must be constant.
\end{proof}
\vskip 0.4cm

The importance of Salmon's Theorem in our story is that it gives a link from ternary cubics to
binary quartics.

Precisely, Salmon's Theorem gives a map
from $SL(3)$-classes of ternary cubics to $SL(2)$-classes of binary cubics.

{\it The two invariant rings are essentially the same.}
Indeed the pullback of $I$ is the invariant $\mathcal{A}$ of degree $4$, called the Aronhold invariant.
The pullback of $J$ is the invariant $\mathcal{T}$ of degree $6$.

Let's state the result precisely.

\begin{thm0}\label{thm:inv_ternary_cubics} The invariant ring of ternary cubics
$\CC[\sym^3\CC^3]^{SL(3)}$ is the polynomial ring $\CC[\mathcal{A},\mathcal{T}]$.
\end{thm0}

 A ternary cubic $g$ with $\mathcal{A}(g)=0$ is called equianharmonic. In other words a cubic is equianharmonic
 when the four branch points of the projection, centered at any of its points, make an equianharmonic binary quartic.
Salmon's Theorem says that if projecting from one point we get an equianharmonic binary quartic,
then projecting from any point we get an equianharmonic binary quartic.
 A basic example is $g=x^4+y^4+z^4$, the Fermat cubic. The Aronhold invariant can be written as a Pfaffian, see \cite[\S 4.7]{O12} and Appendix \ref{appendix:Pfaffians}.

  A ternary cubic $g$ with $\mathcal{T}(g)=0$ is called harmonic. In other words a cubic is harmonic
 when the four branch points of the projection centered at any of its points make a harmonic binary quartic.

A nice feature of these two special cases is that they correspond to the two special cases of lattices.

\begin{prop0}\label{prop:G2G3}
Given a lattice $\Lambda\subset\CC$, $\Lambda=\langle 1, \tau\rangle$ with $\mathrm{Im}\tau>0$, there is a correspondence between
automorphisms of the elliptic curve $\CC/\Lambda$ which fix a point $P_0\in\CC/\Lambda$ and complex numbers $\alpha$ such that $\alpha\Lambda=\Lambda$.

We have moreover

$$\left\{\alpha\in\CC | \alpha\Lambda=\Lambda\right\}\simeq\left\{\begin{array}{ll}\Z_6=\{e^{\pi ik/3}|k=0,\ldots 5\}&\textrm{\ if\ }\Lambda=\langle 1, e^{\pi i/3}\rangle\\
\Z_4=\{i^{k}|k=0,\ldots 3\}&\textrm{\ if\ }\Lambda=\langle 1, i\rangle\\
\Z_2=\{\pm 1\}&\textrm{\ otherwise.}
\end{array}\right.$$
The case $\Lambda=\langle 1, e^{\pi i/3}\rangle$ corresponds to the equianharmonic curves. 
The case $\Lambda=\langle 1, i\rangle$ corresponds to the harmonic curves. 
\end{prop0}
\begin{proof}
The equianharmonic case is \cite[Example IV 4.20.2]{Hart}. In terms of the lattice $\Lambda$ it corresponds to the vanishing of $\sum_{\omega\in\Lambda\setminus\{0\}}\frac{1}{\omega^4}$
(see \cite[Theorem IV 4.12B]{Hart}). 

The vanishing follows from the $\Z_6$-invariance of the lattice in the following way.
Due to absolute convergence, we may reorder the summation and sum over each $\Z_6$-orbit.
Let $\tau=e^{\pi i/3}$, which satisfies $\tau^6=1$. The six summands of the orbit of any fixed $\omega\in\Lambda\setminus\{0\}$ give
$\sum_{j=0}^5(\omega\tau^j)^{-4}=\omega^{-4}(\tau^0+\tau^{-4}+\tau^{-2}+\tau^{0}+\tau^{-4}+\tau^{-2})=0$.
Hence the sum over any orbit is zero and it follows $\sum_{\omega\in\Lambda\setminus\{0\}}\frac{1}{\omega^4}=0$.

The harmonic case is \cite[Example IV 4.20.1]{Hart}.  In terms of the lattice $\Lambda$ it corresponds to the vanishing of $\sum_{\omega\in\Lambda\setminus\{0\}}\frac{1}{\omega^6}$.

The vanishing follows from the $\Z_4$-invariance of the lattice in the following way.
Due to absolute convergence, we may reorder the summation and sum over each $\Z_4$-orbit.
The four summands of the orbit of any fixed $\omega\in\Lambda\setminus\{0\}$ give
$\sum_{j=0}^3(\omega i^j)^{-6}=\omega^{-6}(i^0+i^{-2}+i^{0}+i^{-2})=0$.
Hence the sum over any orbit is zero and it follows $\sum_{\omega\in\Lambda\setminus\{0\}}\frac{1}{\omega^6}=0$.

These sums are called Eisenstein series, see Theorem \ref{thm:G2G3}.

\end{proof}
 
 The nullcone for ternary cubics consists of the closure of the orbit of $y^2z-x^3$, a curve with a cusp.
 The closure contains the following five orbits, we list the representatives and the dimension
 
 \begin{figure}[H]
 $$\begin{array}{lll}&representative&dimension\\
 \hline\\
 \includegraphics[height=1cm, angle=90]{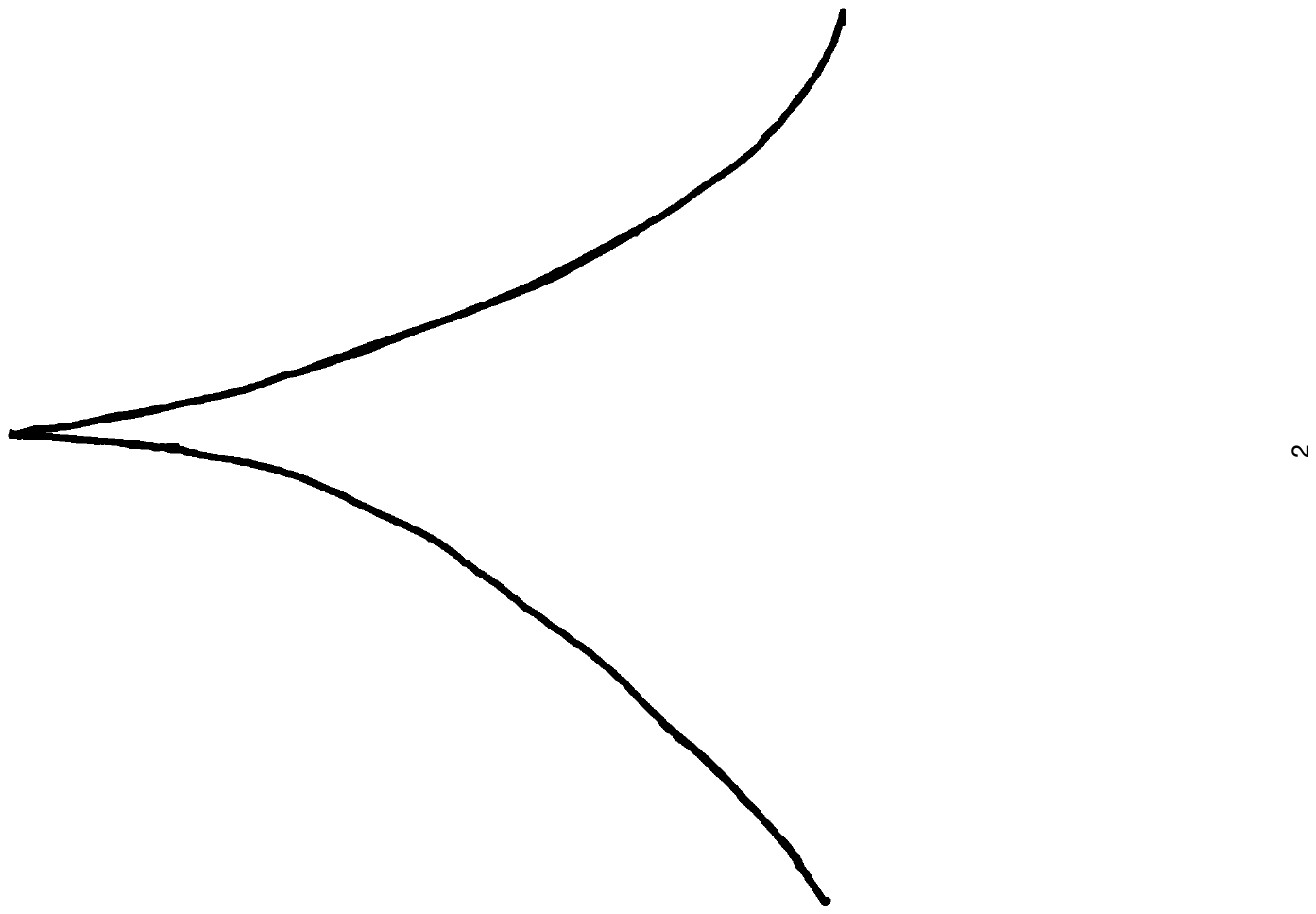}&y^2z-x^3&7\\
 \includegraphics[height=1cm, angle=90]{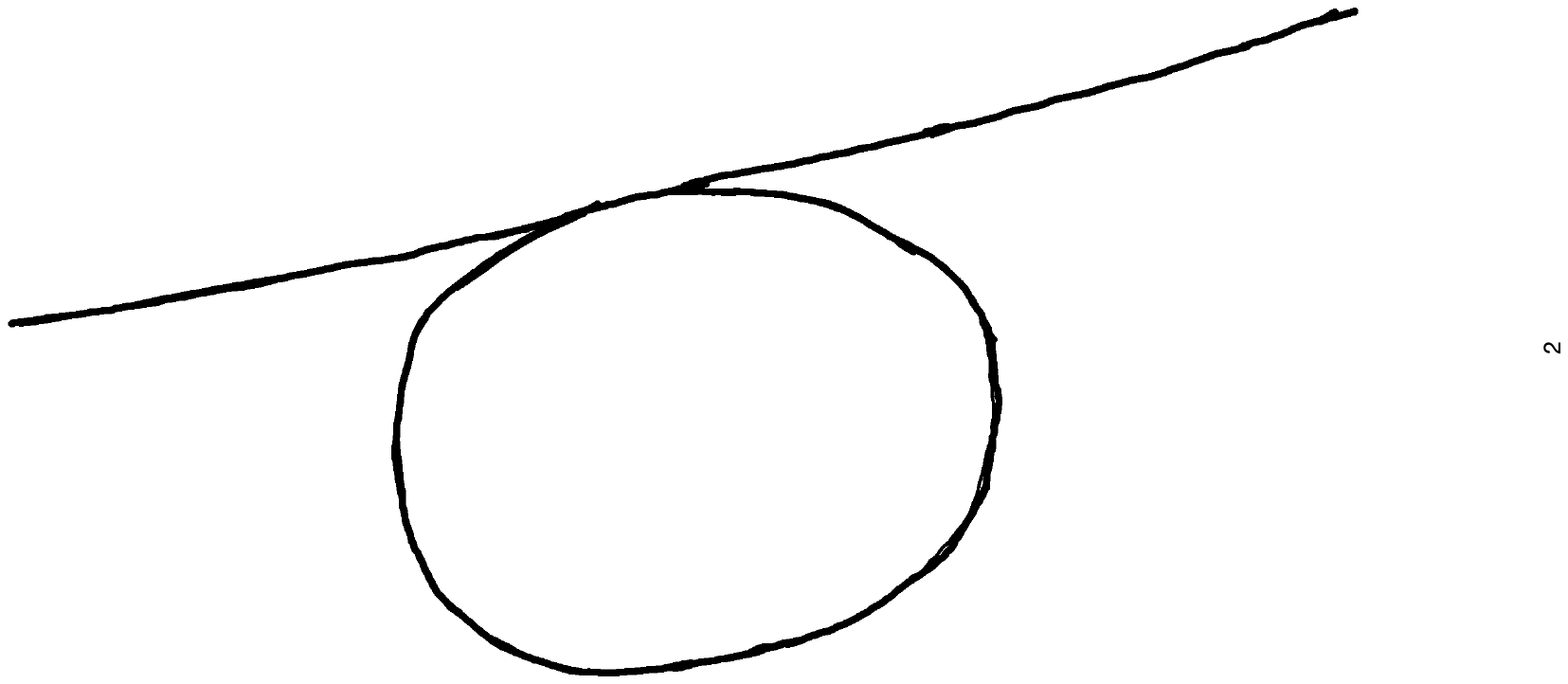}&x(y^2-xz)&6\\
  \includegraphics[height=1cm, angle=90]{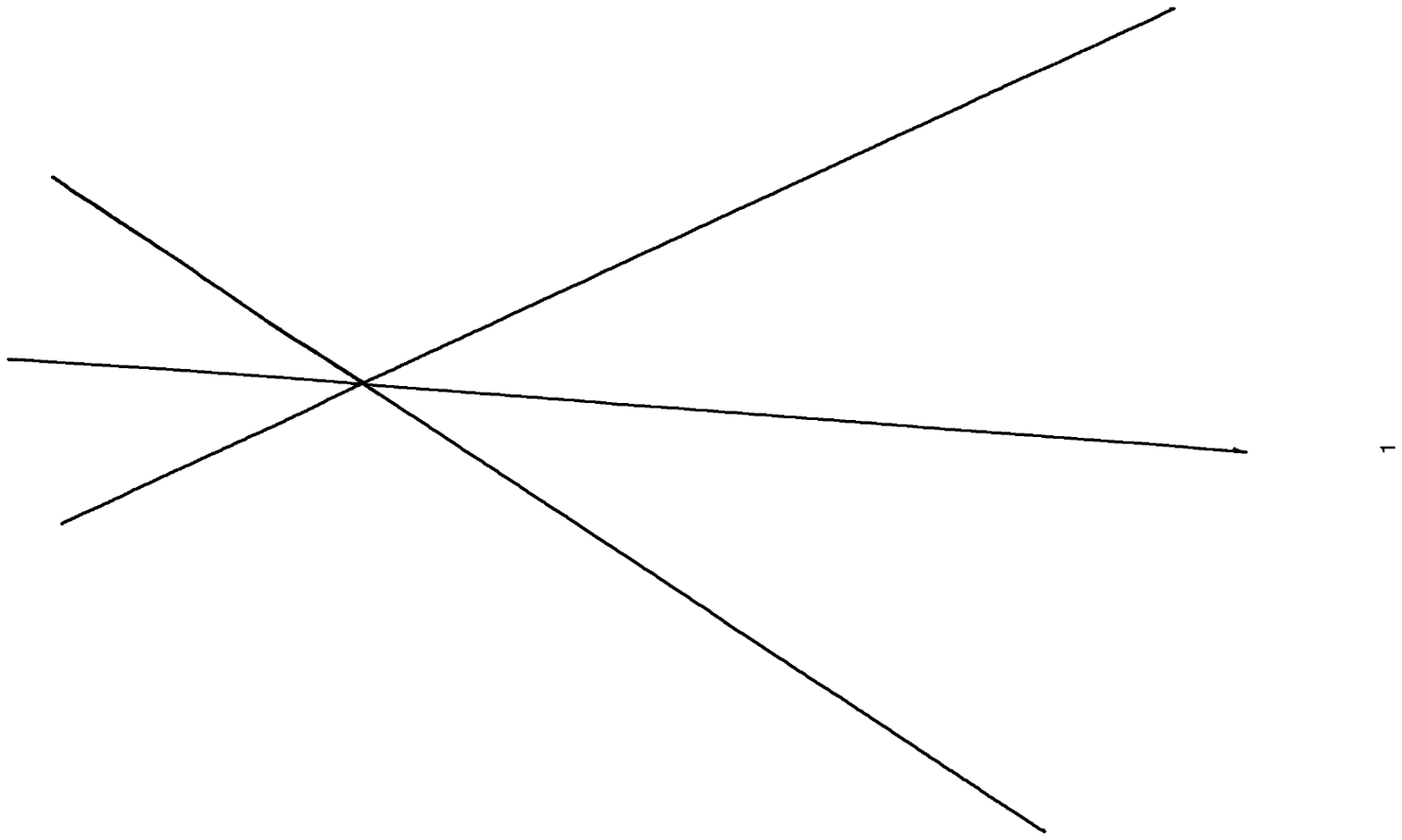}&xy(x+y)&5\\
  \includegraphics[height=1cm, angle=90]{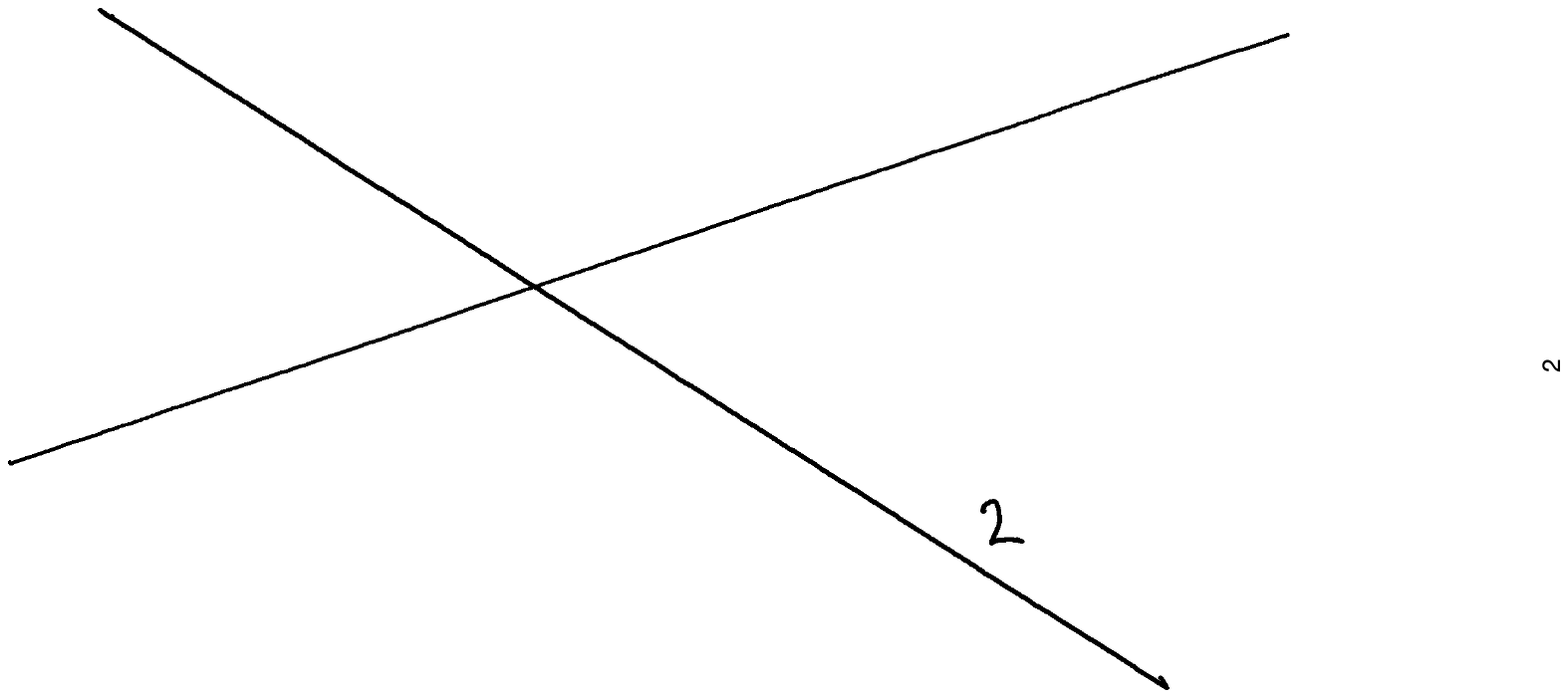}&x^2y&4\\
   \includegraphics[height=1cm, angle=90]{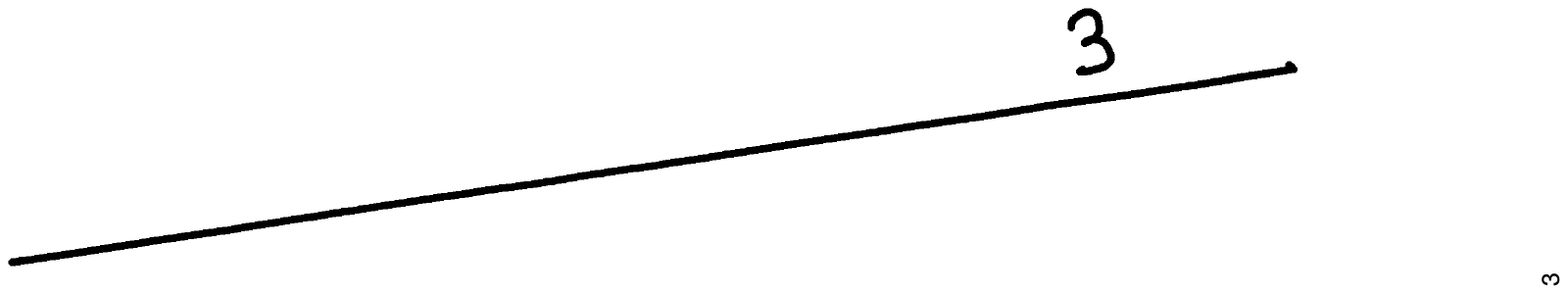}&x^3&2\\
 \end{array}$$
 \caption{The $SL(3)$-orbits in the nullcone for ternary cubics}\label{fig7}
 \end{figure}

 Note that the ``triangle'' $xyz$ is not in the nullcone and it is indeed semistable.
 
The coefficients of the Hessian (the determinant of the Hessian matrix) span a subspace of cubics, of dimension $10$, and they generate an ideal that we call $\mathcal{H}$.

Partial derivates of Aronhold invariant span another subspace of cubics, of dimension $10$, and they generate a second ideal, that we call $\mathcal{A}$.

The following computation, straightfoward with M2 \cite{GS}, shows that the zero locus of $\mathcal{H}$
has a scheme-theoretical non reduced structure (saturated and radical ideal are different!), see \cite[Proposition 4.3]{CO}. 

\begin{prop0}  Let $\mathcal{H}$ be the ideal generated by the $10$ cubics which are the coefficients of the Hessian,
let $\mathcal{A}$ be the  ideal generated by the $10$ partial derivatives of the Aronhold invariant, let $\mathcal{M}=(a_0,\ldots, a_9)$ be the maximal homogeneous ideal. Then the following holds:

$\mathcal{H}\colon\mathcal{M}^\infty = \mathcal{H}$

$\sqrt{\mathcal{H}} = \mathcal{A}\oplus\mathcal{H}$

$\mathcal{A}\colon\mathcal{M}^\infty = \mathcal{A}\oplus\mathcal{H}$

$\sqrt{\mathcal{A}} = \mathcal{A}\oplus\mathcal{H}$
\end{prop0}

The Hessian map
$\PP^9\drig{\mathcal{H}}\PP^9$ has degree $3$\cite[Theor. 4.7]{CO}. 

The map
$\PP^9\drig{\mathcal{A}}\PP^9$ has again degree $3$. 

The analog of the pencil (\ref{eq:pencil_quartic}) of binary quartics is the Hesse pencil
\begin{equation}\label{eq:hessepencil}x^3+y^3+z^3+6sxyz\end{equation}
This family shares the same $9$ flexes and it has the property that the Hessian of a member 
of the Hesse pencil is still a member of the Hesse pencil \cite[Proof of Theorem 4.7]{CO}.

The anharmonic invariant of the Hesse pencil is $S=I'=s(s^3-1)$. The harmonic invariant is $T=J'=8s^6+20s^3-1$,
the discriminant is a multiple of $J'^2-64I'^3$ and it is proportional to $8s^3+1$.

The binary quartic computed by Salmon's Theorem after projecting the Hesse pencil from $(1,-1,0)$
is easily seen to be
$(x+2sy)(x^3-6x^2ys-4y^3)$.

The correspondence given by Salmon's Theorem between (\ref{eq:pencil_quartic}) and (\ref{eq:hessepencil})
is not trivial. By computing the invariant $j=I^3/D$ in both cases we get the correspondence

$$\frac{(1+3t^2)^3}{(3t-1)^2(3t+1)^2}=\frac{s^3(1-s^3)^3}{8(s^3+1/8)^3}.$$

A real picture of this algebraic correspondence is the following.

\begin{figure}[H]
	\centering
\includegraphics[height=8cm]{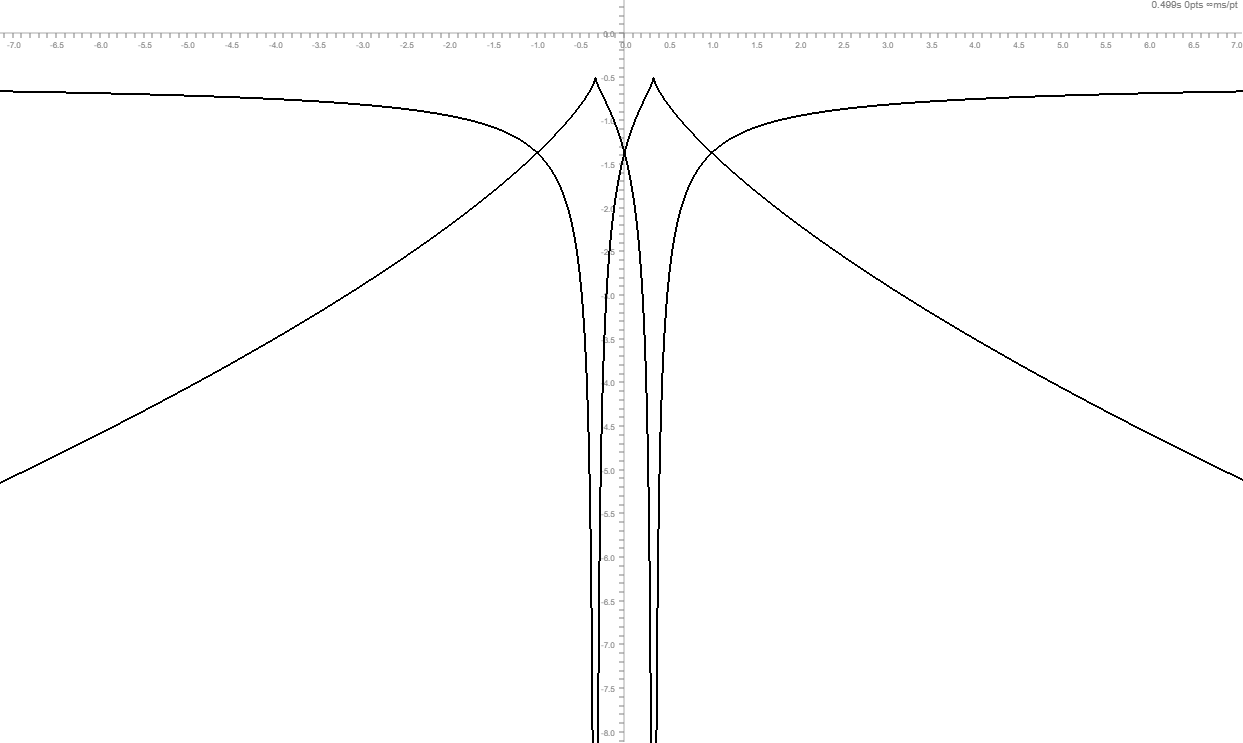}
\caption{\it Correspondence, according to Salmon's Theorem, between real canonical pencils of binary quartics (abscissa) and ternary cubics (ordinate). In the two cusps the discriminant vanishes. The three nodes
correspond to the harmonic case. The equianharmonic case does not appear on $\R$.}\label{fig8}
\end{figure}

On the abscissa the parameter $t$ of the binary quartic \eqref{eq:pencil_quartic}, on the ordinate the parameter $s$ of the ternary cubic \eqref{eq:hessepencil}.
The two cusps correspond to the vanishing of discriminant, $(t,s)=(\pm 1/3, -1/2)$. The three nodes
correspond to the harmonic case and have the same ordinate $s=-1.366\ldots$, the only negative real root of $8s^6+20s^3-1$
and the abscissa $t\in\{-1, 0, 1\}$. Note the Hesse pencil consists of two connected components for $s\in\R$, $s<-1/2$.

\begin{rem0}
Elisabetta Rocchi has studied in \cite{Roc} a further analogy between binary quartics and ternary cubics.
The group $SL(2)$ acts on the pencils of binary quartics inside the Grassmannian $Gr(\PP^1,\PP^4)$. The orbit of the pencil  (\ref{eq:pencil_quartic})  is exactly the Fano 3-fold $V_5$
of item 3. of Theorem \ref{thm:AF} . In the same way, the group $SL(3)$ acts on the pencils of ternary cubics.
The orbit of the Hesse pencil, inside the Grassmannian  $Gr(\PP^1,\PP^9)$ is a $8$-fold $X$ described in \cite{Roc}. It consists of nine orbits and has a $4$-dimensional singular locus made up of two orbits. It should be interesting to know if $X$ is $\QQ$-Cartier and, in the affirmative case, if $X$ is Fano.
\end{rem0}

\section{Finite polyhedral groups and the ADE classification}\label{sec:ADE}

We recall that over the real numbers we have the $2:1$ covering of real Lie groups $SU(2)\to SO(3,\R)$.

The polyhedral groups can be found as finite subgroups of $SO(3,\R)$.
Their pullback in $SU(2)$ are called {\it binary polyhedral groups}.

The complex version is the $2\colon 1$ covering $SL(2,\CC)\to SO(3,\CC)$.

A finite subgroup of $SO(3,\R)$ acts on the sphere.
The cyclic subgroup $\Z_n$ is conjugated to the group generated by
\begin{equation}\label{eq:dn}\begin{pmatrix}\cos(\pi/n)&-\sin(\pi/n)&0\\
\sin(\pi/n)&\cos(\pi/n)&0\\
0&0&1\end{pmatrix},\end{equation} which acts as rotation of $2\pi/n$ around the axis $South-North$.
 The Earth covered by time zones gives a good representation.
Each time zone is a fundamental domain $T$ for the action of the cyclic group, meaning that
$\cup_{g} gT=S^2$, $\mathring{T}\cap g\mathring{T} = \emptyset\quad\forall g\in\Z_n, g\neq 1$.
A time zone can be seen as a degenerate spherical triangle with angles $2\pi/n$, $2\pi/n$, $\pi$.

Klein classified  the finite groups $G\subset SO(3,\R)$ and we review his classification.

There are two basic classes, namely the cyclic group just described 
  and the dihedral group $D_n$ with $2n$ elements.

In the cyclic case, the pullback to $SU(2)$ is again a cyclic group of double order $\Z_{2n}$,
generated by one of the two pullbacks of a generator of $\Z_n$.
The fundamental domain for the action of the group $\Z_{2n}$ is just half of the previous domain,
drawing the bisector of the angle $2\pi/n$.

The dihedral group $D_n$ of size $2n$ is the symmetry group of a regular $n$-gon
and it
can be embedded in $SO(3)$
  with the trick of swapping North and South Pole, like in
\begin{equation}\label{eq:dn1}\begin{pmatrix}\cos(\pi/n)&\sin(\pi/n)&0\\
\sin(\pi/n)&-\cos(\pi/n)&0\\
0&0&-1\end{pmatrix}\end{equation}
 
 The $k$-powers of (\ref{eq:dn1}) and (\ref{eq:dn1}) for $k=0,\ldots, n-1$ are all the elements of $D_n$.
The sphere can be tasselated by $2n$ geodesic triangles, $n$ above the equator and $n$ 
below the equator, each one with angles
$\frac{\pi}{2},\frac{\pi}{2},\frac{2\pi}{n}$.
  The pullback to $SU(n)$ is the binary dihedral group, of size $4n$.
  
    Fundamental domain is again a triangle which is half of the previous ones,
  namely with angles $\frac{\pi}{2},\frac{\pi}{2},\frac{\pi}{n}$.

The three most interesting subgroups correspond to the classes of Platonic solids under duality, 
namely the tetrahedral $A_4$ , the esahedral/octahedral group $S_4$ and the dodecahedral/icosahedral group $A_5$. Each face is a regular polygon with $p$ sides. Join the center of the face with all the vertices of the same face. We get $p$ triangles, which are the fundamental domains for the groups in $SO(3)$.

The following table shows the group size and the angles of a geodesic triangle (here $p$ means the angle $\pi/p$), which is a fundamental domain for the binary group. The size of any binary group is twice the size of the corresponding group. Note the first table covers exactly the binary polyhedral groups.

\begin{equation}\label{eq:Dynkin}\begin{array}{c|c|c||c|c|c}
&&&p&q&r\\
\hline
\textrm{cyclic}&\Z_{n}&n&1&n&n\\
\hline
\textrm{dihedral}&D_{r}&2r&2&2&r\\
\hline
\textrm{tetrahedral}&A_4&12&2&3&3\\
\hline
\textrm{octahedral}&S_4&24&2&3&4\\
\hline
\textrm{icosahedral}&A_5&60&2&3&5
\end{array}\end{equation}

With this interpretation the polyhedral group are seen as {\it elliptic triangle groups}.

Each oriented graph $Q$ is defined by a set  $I$ of vertices, a set $\Omega$ of  edges  and source and target maps
$s, t\colon\Omega\to I$.  This is called a quiver. Every finite quiver $Q$ carries a quadratic form on $\Z^I$
$$\langle v, w\rangle = \sum_{i\in I}v_iw_i-\sum_{h\in\Omega}v_{s(h)}w_{t(h)}$$
and its symmetrized version
$$(v, w)=\langle v, w\rangle +\langle w, v\rangle$$

The associated quadratic form $q_Q(v)=\frac 12(v, v)$ is called the Tits form\cite[\S 1.5]{Kir}.

Let now $\Gamma(p, q, r)$ be the ``star'' graph  consisting of three branches of lengths $(p, q, r)$ meeting at a central vertex. If one of $p$, $q$, $r$ is $1$, then there are only two branches. It is an elementary fact\cite[Theor. 1.30]{Kir} that $q_{\Gamma(p, q, r)}$ is positive definite if and only if
\begin{equation}\label{eq:pqr}\frac 1p +\frac 1q +\frac 1r >1\end{equation}

It is time to list all the integer solution to the following inequalities

$$\begin{array}{c}(1/p+1/q+1/r)>1\\
{\begin{array}{ccc}p&q&r\\
\hline\\
2&2&r\ge 2\\
2&3&3\\
2&3&4\\
2&3&5\end{array}}\end{array}\qquad
\begin{array}{c}(1/p+1/q+1/r)=0\\
{\begin{array}{ccc}p&q&r\\
\hline\\
2&3&6\\
2&4&4\\
3&3&3\\
&&\phantom{}
\end{array}}\end{array}$$

Each of the groups in \ref{eq:Dynkin} corresponds to a ``star'' graph with the corresponding $(p, q, r)$
The correspondence is the following
\begin{figure}[H]
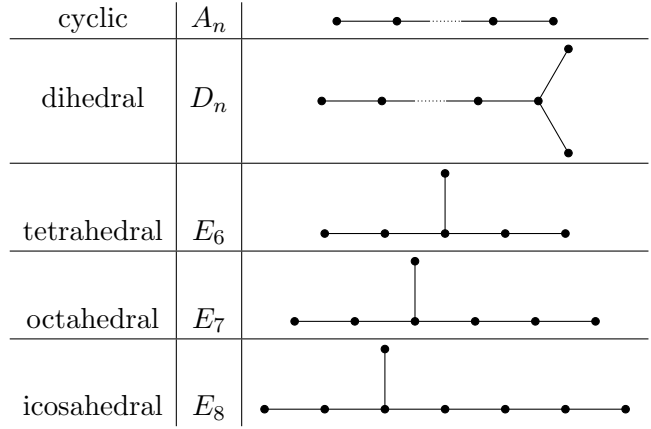

$$\begin{array}{c|c|c}
\textrm{cyclic} & A_n & \dynkin[edge length=0.8cm] A{}\\
\hline
\textrm{dihedral} & D_n & \dynkin[edge length=0.8cm] D{}\\
\hline
\textrm{tetrahedral} & E_6  &\dynkin[edge length=0.8cm]{E}{6}\\
\hline
\textrm{octahedral} & E_7  &\dynkin[edge length=0.8cm]{E}{7}\\
\hline
\textrm{icosahedral} & E_8  &\dynkin[edge length=0.8cm]{E}{8}
\end{array}$$
\caption{Dynkin diagrams for the polyhedral groups. The ADE classification.}\label{fig:DynkinADE}
\end{figure}

The graphs listed are exactly the Dynkin diagrams and they coincide with the 
positive integer solution to (\ref{eq:pqr}).

The fundamental theorem is that $q_Q$ is positive definite if and only if $Q$ is one of the graphs in 
Figure \ref{fig:DynkinADE}.
\cite[Theor. 1.28]{Kir}.
This does not depend on the orientation of the graph.

We will sketch in Remark \ref{rem:McKay} the McKay correspondence, which is an alternative and more ``modern'' way to see the same ADE correspondance between graphs and groups.

We summarize the tesselations of the sphere in geodesic triangles discussed before, as fundamental domains
for the group action.
We recall the elementary

\begin{prop0}\label{prop:areatriang}
The area of a geodesic triangle on the unitary sphere, with angles
$\frac{\pi}{p}$, $\frac{\pi}{q}$, $\frac{\pi}{r}$ is equal to
$ \frac{\pi}{p} +\frac{\pi}{q}+\frac{\pi}{r}-\pi$ .
In particular we have the inequality
$$\frac 1p +\frac 1q +\frac 1r >1$$
\end{prop0}

Conversely, given a triangle with angles $\frac{\pi}{p}$, $\frac{\pi}{q}$, $\frac{\pi}{r}$,
such that
$$\frac 1p +\frac 1q +\frac 1r >1$$
we have a tessellation of the sphere with triangles congruent to the given one and a finite group $G\subset SO(3)$ acting on the sphere such that the triangle is a fundamental domain.
The number of triangles coincides with $|G|$ and it is equal to the ratio of areas

$$\frac{4}{\frac 1p +\frac 1q +\frac 1r -1}$$ 

Klein had the wonderful idea to recover the semi-invariants of the polyhedral groups from the forms giving the vertices. The computations are similar for all the polyhedra. The yoga of this construction is contained in these two amazing facts

$$\left\{\begin{array}{rcc}\textrm{Hessian(Vertices)} &=& \textrm{Faces}\\
\textrm{T(Vertices)} &=& \textrm{Edges}\end{array}\right.$$

where we recall that the covariant $T$ of a binary form $f$  is
$T(f)=\det\begin{pmatrix}f_x&f_y\\
H_x&H_y\end{pmatrix}$, where $H=Hess(f)$. The invariant $T$ will be crucial in \S  \ref{sec:Hilbert1886},
discussing Hilbert paper, see Theorem \ref{thm:invT}.

We start from the tetrahedron, which, as can be expected by Section \ref{sec:binary_quartics}, it involves the equianharmonic binary quartics seen in Prop. \ref{prop:IJ}.
We use the same notations as in \cite{Klein}.

\begin{prop0}
Let $\Phi= x^4+2\sqrt{-3}x^2y^2+y^4$ be a equianharmonic quartic. Its roots are vertices of a tetrahedron.
Then

$\Psi=Hess(\Phi)=x^4-2\sqrt{-3}x^2y^2+y^4$  is again equianharmonic. Its roots correspond to the
centroids of every face (this is sometimes called the countertetrahedron)

$t=T(\Phi)=xy(x^4-y^4)$ is a sextic whose roots correspond to the midpoint of each edge. They are the vertices of a octahedron.
\end{prop0}



Given a ring $R$ with a $G$-action, the invariant ring is
$$R^G:=\{f\in R| g\cdot f=f\quad\forall g\in G\}.$$
The semi-invariant ring is defined through the one dimensional representation, namely group morphism $\chi\colon G\to \CC^*$. They form a group $\Gamma_G=\mathrm{Hom}(G, \CC^*)$, sometimes called the character group of $G$ (don't confuse
these one-dimensional characters with higher dimensional characters, which are just class functions
and may vanish).
The semi-invariant ring is 
\begin{equation}\label{eq:semi-inv}R^G_s:=\{f\in R| \exists\chi\in\Gamma_G\textrm{\ such that\ } g\cdot f=\chi(g)f\quad\forall g\in G\}.\end{equation}
Note the inclusion $R^G\subset R^G_s$.

\textbf{Teorema 1 (Molien).} \cite[§2.2]{Stu}\label{thm:molien}

\[
\sum_{d=0}^{\infty}
\dim K[u_{1},\dots,u_{n}]^{G}_{d}\, t^{d}
=
\frac{1}{|G|}
\sum_{A\in G}
\frac{1}{\det(I_{n}-tA)}.
\]

Molien Theorem has the following extension to the semi-invariants.

\begin{prop0}\label{prop:semiMolien}
Let $G$ be a finite group acting on $U$.

\begin{enumerate}
\item{}The invariant subring induced on the symmetric algebra $S^*U$ has Hilbert series
\begin{equation}\label{molienformula}\sum_{i=0}^{+\infty}\dim \left(S^iU\right)^Gt^i=\frac{1}{|G|}\sum_{g\in G}\frac{1}{\det(1-tg)}\end{equation}
where $g$ acts on $U$.

\item{}If $\chi\colon G\to\CC^*$ is a one dimensional character,
the subring of semi-invariants with respect to $\chi$
is defined by $(S^iU)_{\chi}=\left\{f\in S^i U | g\cdot f=\chi(g) f\right\}$
Its Hilbert series is

\begin{equation}\label{molienchi}\sum_{i=0}^{+\infty}\dim \left(S^iU\right)_\chi t^i=\frac{1}{|G|}\sum_{g\in G}\frac{\overline{\chi(g)}}{\det(1-tg)}\end{equation}
 \end{enumerate}
\end{prop0}
\begin{proof} The first case is \cite{Stu} Theor. 2.2.1
The second case is an easy modification, note that when $\chi=id$ 
$(S^iU)_{id}=(S^iU)^G$.
 \end{proof}

The character table of $A_4\subset SO(3,\R)$ is (see \cite[Exerc. 2.26]{FH} )

\begin{equation}\label{eq:A4char}\begin{array}{r|c|c|c|c}
\textrm{size}&1&3&4&4\\
\hline
\rho_1&1&1&1&1\\
\rho_2&1&1&\omega&\omega^2\\
\rho_3&1&1&\omega^2&\omega\\
\rho_4&3&-1&0&0
\end{array}
\end{equation}

This knowledge is enough to understand the semi-invariant rings of both the tetrahedral and the
binary tetrahedral groups.

The interesting feature of (\ref{eq:A4char}) is the appearance of three one dimensional representations.
The first row gives the identity, while the second and the third row are dual to each other.

The tetrahedral group, as a subset of $SO(3)$, is Alt(4), it has $12$ elements.

Let $\varphi\colon SU(2)\to SO(3,\R)$.
Define
\begin{equation}\label{eq:R}
\rho' : \mathbb{C}^{2} \to \mathbb{C}^{2}, \qquad
\begin{pmatrix} x \\ y \end{pmatrix}
\mapsto
\begin{pmatrix}
i & 0 \\
0 & -i
\end{pmatrix}
\begin{pmatrix} x \\ y \end{pmatrix}
= R \begin{pmatrix} x \\ y \end{pmatrix},
\end{equation}

\begin{equation}\label{eq:S}
\sigma' : \mathbb{C}^{2} \to \mathbb{C}^{2}, \qquad
\begin{pmatrix} x \\ y \end{pmatrix}
\mapsto
\frac{1}{2}
\begin{pmatrix}
1+i & -(1-i) \\
1+i & 1-i
\end{pmatrix}
\begin{pmatrix} x \\ y \end{pmatrix}
= S \begin{pmatrix} x \\ y \end{pmatrix}.
\end{equation}

The binary tetrahedral group is
$H_{24} = \varphi^{-1}(\mathrm{Alt}(4))$, it has $24$ elements and it is generated by
$R$ and $S$. 
We have

\[
\rho = \varphi(R) =
\begin{pmatrix}
1 & 0 & 0 \\
0 & -1 & 0 \\
0 & 0 & -1
\end{pmatrix},
\qquad
\sigma = \varphi(S) =
\begin{pmatrix}
0 & -1 & 0 \\
0 & 0 & -1 \\
-1 & 0 & 0
\end{pmatrix}.
\]

where $\rho$ and $\sigma$ generate $Alt(4)$,
Since for any $A \in K^{n \times n}$ the characteristic polynomial $\det(I_{n} - tA)$ is invariant by similarity
the Hilbert series of $\mathbb{C}[u,v]^{H_{24}}$ is

\[
\begin{aligned}
\sum_{d=0}^{+\infty}\dim\bigl(\mathbb{C}[u,v]_{d}^{H_{24}}\bigr)t^{d}
&{=}
\frac{1}{|H_{24}|}\sum_{A\in H_{24}}\frac{1}{\det(I_{2}-tA)}\\[1mm]
&=\frac{1}{24}\Biggl(
\frac{1}{\det(I_{2}-tI_{2})}
+\frac{1}{\det(I_{2}+tI_{2})}
+\frac{3}{\det(I_{2}-tR)}
+\frac{3}{\det(I_{2}+tR)}\\
&\hspace{1.4cm}
+\frac{4}{\det(I_{2}-tS)}
+\frac{4}{\det(I_{2}+tS)}
+\frac{4}{\det(I_{2}-tS^{2})}
+\frac{4}{\det(I_{2}+tS^{2})}\Biggr)\\[1mm]
&=\frac{1}{24}\Biggl(
\frac{1}{(1-t)^{2}}
+\frac{1}{(1+t)^{2}}
+\frac{3}{1+t^{2}}
+\frac{3}{1+t^{2}}\\
&\hspace{1.4cm}
+\frac{4}{1-t+t^{2}}
+\frac{4}{1+t+t^{2}}
+\frac{4}{1+t+t^{2}}
+\frac{4}{1-t+t^{2}}
\Biggr).
\end{aligned}
\]

Now note that $(1-t^6)=(1-t^2)(1+t+t^2)(1-t+t^2)$.
Hence the Hilbert series simplifies to

$$\frac{1-t^4+t^8}{(1-t^6)(1-t^4)}$$

Multiply the numerator and the denominator by $(1+t^{4})(1-t^{12})$, we get
\[
\sum_{d=0}^{+\infty}\dim\bigl(\mathbb{C}[u,v]_{d}^{H_{24}}\bigr)t^{d}=\frac{1-t^{24}}{(1-t^{6})(1-t^{8})(1-t^{12})}.
\]

This gives the presumption that  $\mathbb{C}[u,v]^{H_{24}}$ has $3$ generators of degree respectively $6,8,12$, linked by a unique algebraic relation of degree $24$.
This claim is indeed correct and one can show that three possible invariants of these degrees are respectively

\begin{enumerate}
\item $x^5y-xy^5$, the midpoints of the edges of the tetrahedron,

\item $(x^4-2\sqrt{3}ix^2y^2+y^4)(x^4+2\sqrt{3}ix^2y^2+y^4)$,
where each of the two factors is equianharmonic and vanishes respectively on the vertices $\{\pm a, \pm b\}$ of the rhombus in Figure \ref{fig:IJ} and on their rotation by $\pi/2$,
they correspnd to the vertices of the tetrahedon and its dual tetrahedron,

\item  $(x^4-2\sqrt{3}ix^2y^2+y^4)^3$.

\end{enumerate}

The interest of this invariant ring is milder than expected. Note that $g=(x^4-2\sqrt{3}ix^2y^2+y^4)$ is not an invariant but it is only a semi-invariant, indeed
the action of $R$ in (\ref{eq:R}) leaves $g$ invariant while the action of $S$ in (\ref{eq:S}) multiplies $g$ by $e^{2\pi i/3}$. In the next subsections we perform the computation of the semi-invariant ring of the three famous polyhedral groups, which turn out to contain the basic geometrical informations of the three polyhedra. 

\subsection{Semi-invariants of binary tetrahedral group}

\begin{figure}[H]
	\centering
 \includegraphics[width=5cm]{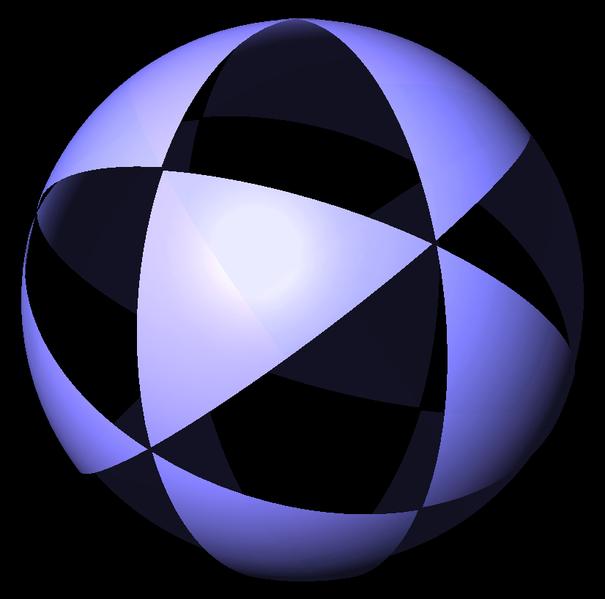} 
\caption{Binary tetrahedral group, $(2,3,3)$, $24$ triangles.}\label{fig:bintetra}
\end{figure}
\vskip 0.8cm

The tesselation with $24$ triangles shown in Figure\ref{fig:bintetra} arises
if we inscribe a tetrahedron in the sphere and we proceed with the barycentric subdivision of each face, joining
the center of the face with each vertex and with each midpoint of an edge. The division is projected on the sphere from its center so that
each triangle becomes a geodesic triangle on the sphere.
Note that a small neighborhood of a vertex may meet more fundamental regions than a small 
neighborhood of a point on an edge. This marks a different behaviour on the quotient space,
that can be understood in the setting of orbifolds, which we do not consider in these notes.

The conjugacy classes of $\mathrm{Alt}(4)$ are $4$, of respective cardinality 
$1,3,4,4$ with representativesi $\varphi(I_{2})$, $\varphi(R)$, $\varphi(S)$, $\varphi(S^{2})$; 
hence the conjugacy classes of the binary tetrahedral group $H_{24}$ are $8$, with cardinality 
$1,1,3,3,4,4,4,4$ and representatives 
$I_{2}, -I_{2}, R, -R, S, -S, S^{2}, -S^{2}$.

\begin{prop0}
Consider the binary tetrahedral group $G$ of order $24$ acting on $S^2$.
It has three one dimensional characters with the following characters (the first row represents the size of each conjugacy class)
$$\begin{array}{c|c|c|c|c|c|c|c}
1&1&6&4&4&4&4\\
\hline
id&1&1&1&1&1&1&1\\
\chi&1&1&1&\omega^2&\omega&\omega&\omega^2\\
\chi'&1&1&1&\omega&\omega^2&\omega^2&\omega\\
\end{array}$$
We use Molien formula in Prop. \ref{prop:semiMolien}.

The semi-invariant ring collects the contribution of the three one dimensional characters.
We get that the Hilbert series of the semi-invariant ring is

$$\frac{1}{24}\left(\frac{3}{(1-t)^2}+\frac{3}{(1+t)^2}+\frac{6\cdot 3}{1+t^2}\right) =
\frac{1-t^2+t^4}{(1-t^2)^2(1+t^2)} =\frac{1+t^6}{(1-t^4)^2}=\frac{1-t^{12}}{(1-t^4)^2(1-t^6)}$$
\end{prop0}
The last expression suggests three generators of degree $4$, $4$, $6$, with a relation in degree $12$.
This is correct and it is shown in the following Proposition.

 \begin{prop0}
 $\forall g\in G$, $g\cdot \Phi = \chi(g)\Phi$,
$g\cdot \Psi = \chi'(g)\Psi$,
$g\cdot t = t$.
 Hence $\Phi$, $\Psi$, $t$ are semi-invariants. The semi-invariant is ring generated by
$\Phi$, $\Psi$, $t$ with the relation
$$12\sqrt{-3}t^2-\Phi^3+\Psi^3=0$$
\end{prop0}

\begin{proof} The ring generated by $\Phi$, $\Psi$, $t$ is a subring of the semi-invariant ring and has the same Hilbert series.
Hence the two rings coincide. This argument works verbatim also in the octahedral and icosahedral cases and we will not repeat it.
\end{proof}

\begin{exe} The semi-invariant ring of the tetrahedral group acting on $\CC^3$
has Hilbert series
$$\frac{1-t^6}{(1-t^2)^3(1-t^3)}$$
and it is generated by $3$ quadratic polynomials and one cubic polynomial with a relation in degree $6$.
\end{exe}

\subsection{Semi-invariants of binary octahedral group}

\begin{figure}[H]
 \includegraphics[width=7cm]{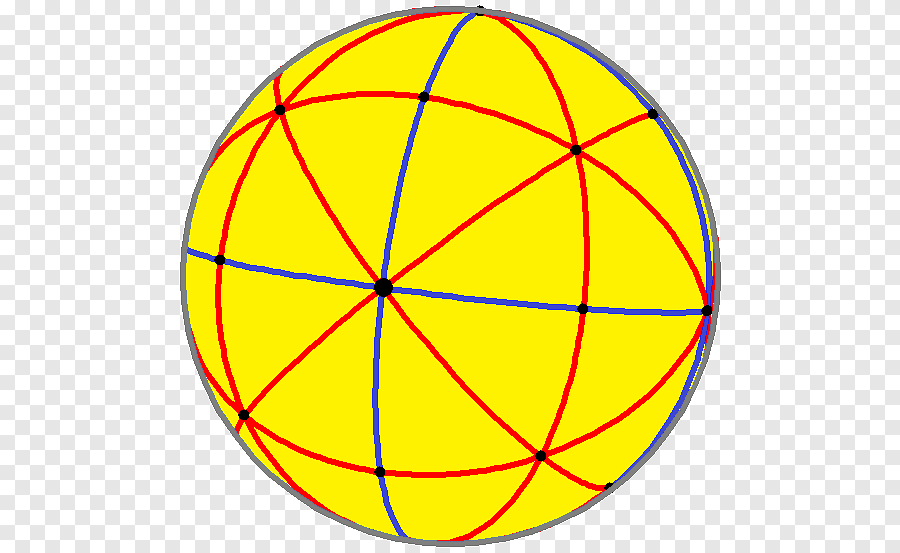} 
\vskip 0.8cm

\caption{Binary octahedral group, $(2, 3, 4)$, $48$ triangles}\label{fig11}

\end{figure}

In the octahedral case there are only two $1$-dimensional characters. We have the table
$$\begin{array}{c|c|c|c|c|c|c|c|c}
1&1&6&6&8&8&18&18\\
\hline
id&1&1&1&1&1&1&1&1\\
\chi&1&1&-1&-1&1&1&-1&1\\
\end{array}$$
By using Molien formula, the Hilbert series of the semi-invariant ring is

$$\frac{1}{48}\left(\frac{2}{(1-t)^2}+\frac{2}{(1+t)^2}+\frac{8\cdot 2}{1-t+t^2}+\frac{8\cdot 2}{1+t+t^2}+
\frac{6\cdot 2}{1+t^2}\right) = \frac{1-t^{24}}{(1-t^6)(1-t^8)(1-t^{12})}$$
The last expression suggests three generators of degree $6$, $8$, $12$, with a relation in degree $24$.

\begin{prop0}
$t=xy(x^4-y^4)$ is a sextic whose roots correspond to the vertices of a octahedron.

Then

$W=Hess(t)=x^8+14x^4y^4+y^8$  has roots correspond to the
centroids of every face.

$\Xi=T(t)=x^{12}-33x^8y^4-33x^4y^8+y^{12}$ has roots correspond to the midpoint of each edge. They are the vertices of a icosahedron.
\end{prop0}
 
 \begin{prop0}
 
$\forall g\in G$, $g\cdot t = \chi(g)t$,
$g\cdot W = W$,
$g\cdot \Xi = \chi(g)\Xi$.

 Hence $t$, $W$, $\Xi$ are semi-invariants. The semi-invariant is ring generated by
$t$, $W$, $\Xi$ with the relation
$$108t^4-W^3+\Xi^2=0$$
\end{prop0}

\subsection{Semi-invariants of binary icosahedral group}
\begin{figure}[H]
	\centering
 \includegraphics[width=6cm]{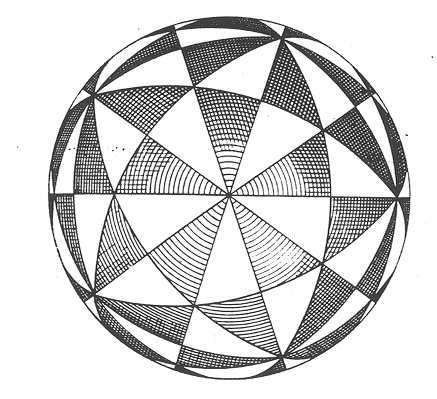} 
\vskip 0.8cm
\caption{Binary icosahedral group, $(2, 3, 5)$, $120$ triangles.}\label{fig12}
\end{figure}

The binary icosahedral group has $10$ conjugacy classes, of cardinality
$1,1,12,12,12,12,15,15,20,20$ 
The icosahedral group $a_5$ is the first simple group, so a new phenomenon arrives, the group character is trivial.
Hence the only semi-invariants are invariant, simplifying the previous computations so that the semi-invariant ring coincides with the invariant ring.
In other words, just the standard form of Molien Theorem \ref{thm:molien} is enough, we do not need Prop. \ref{prop:semiMolien}.

The Hilbert series is

{\footnotesize
$$\begin{aligned} 
\frac{1}{120}\left(
\frac{1}{1-2t+t^2}+\frac{1}{1+2t+t^2}+
12\frac{1}{1-\rho_0 t + t^2}+12\frac{1}{1+\rho_0 t + t^2}
+12\frac{1}{1-\rho_1 t + t^2}+12\frac{1}{1+\rho_1 t + t^2}\right.\\
\left.+15\frac{1}{1+t^2}+15\frac{1}{1+t^2}
+20\frac{1}{1+t+t^2}+20\frac{1}{1-t+t^2}
\right),
\end{aligned}$$}

where $\rho_0 = 2\cos(\frac{2\pi}{5})$ and $\rho_1 = 1 + \rho_0$.

After simplification and multiplication by
\[
\prod_{j=0}^{2} (t^{2j} + t^{2j\pm 1} + t^{2j\pm 20})
\prod_{j=5}^{5} (t^{2j} + t^{2j\pm 6} + t^{2j\pm 15}),
\]
we get

\[
\frac{(1 - t^{12})(1 - t^{20})(1 - t^{30})}{1 - t^{60}}.
\]

This suggests we have generators of degree $12$, $20$, $30$, with a relation in degree $60$..
Again this is correct.
We have indeed

\begin{prop0}
$f=(x^{11}y + 11x^6y^6 - xy^{11})$ is an invariant of degree $12$, representing the vertices of the icosahedron

$H(f)$ has degree 20, representing the faces

$T(f)$ has degree 30, representing the edges.
The invariant ring is generated by $f$, $H(f)$, $T(f)$ with the relation
$$T^2+H^3-1728f^5=0$$
\end{prop0}

The $4$th transvectant (see Theorem \ref{thm:transvectant}) of a binary from $f$ has the expression 
$$(f, f)_4=f_{xxxx}f_{yyyy}-4f_{xxxy}f_{xyyy}+6f_{xxyy}^2$$

\begin{thm0}[Klein]\label{thm:klein}\cite[chapt. I \S 13]{Klein}
A nonsingular real binary form $f$ of degree $4$ vanishes on the vertices of a regular tetrahedron if and only if $(f, f)_4=0$

A  nonsingular real binary form $f$ of degree $6$ vanishes on the vertices of a regular octahedron if and only if $(f, f)_4=0$

A  nonsingular real binary form $f$ of degree $12$ vanishes on the vertices of a regular icosahedron if and only if $(f, f)_4=0$
\end{thm0}

The key to understand the reality condition in Theorem \ref{thm:klein} is that the compact part of $SL(2,\CC)$ is $SU(2)$
which is the universal covering of $SO(3,\R)$. Hence the real elements act as isometries on the sphere
and preserve any regular polyhedron.

As a interlude between the spherical tilings and the hyperbolic tilings of \S \ref{sec:hyptriangle},
it is the right moment to present the next plane tiling.

\begin{figure}[H]
	\centering
 \includegraphics[width=6cm]{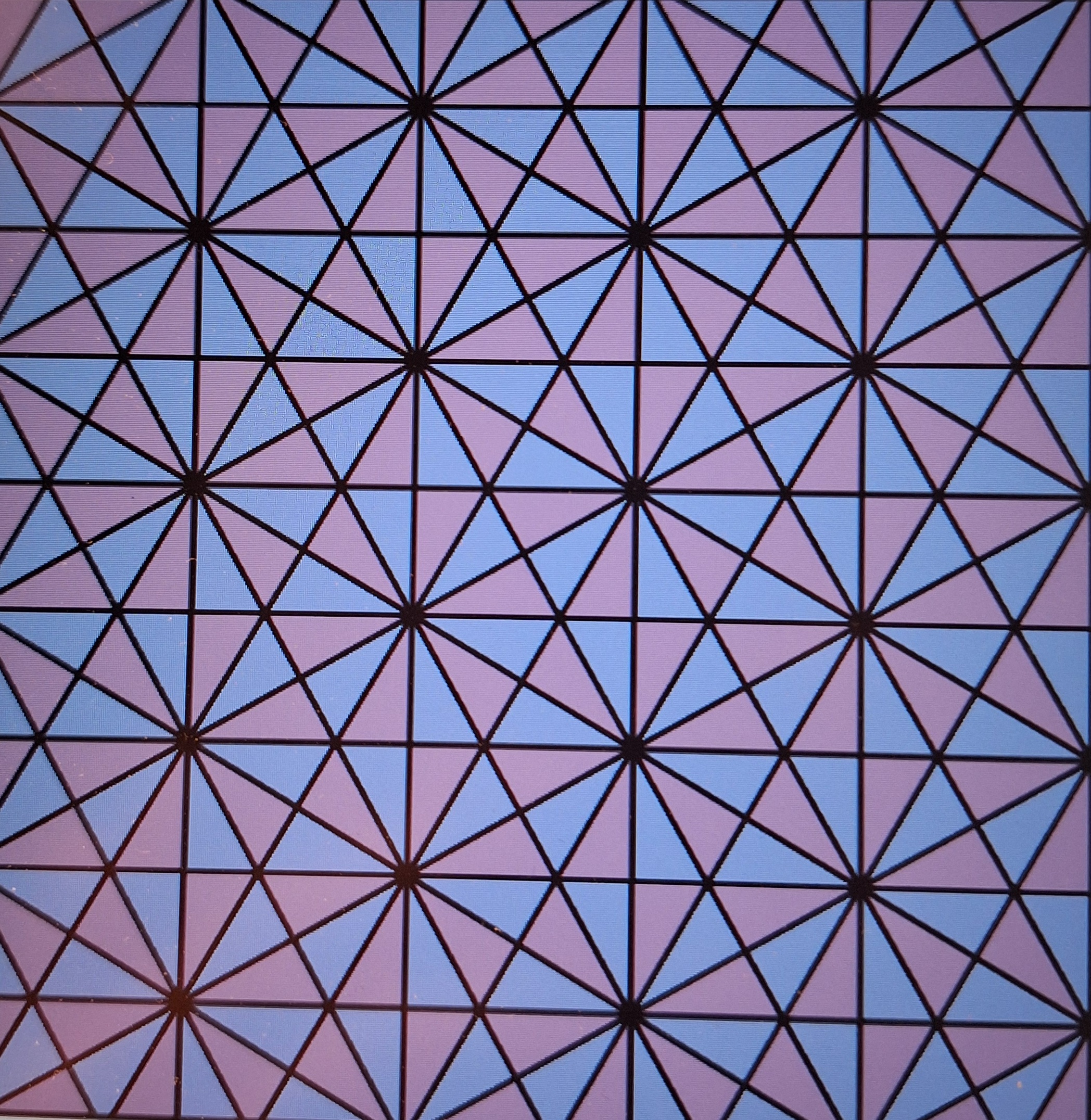} 
\caption{Bisected Hexagonal group, $(2, 3, 6)$}\label{fig13}
 \end{figure}

\begin{rem0}\label{rem:McKay}
In this remark we sketch the McKay correspondence which allows us to construct the extended Dynkin diagrams from the polyhedral groups, namely from their character table.

The irreducible representations of a group $G\subset SL(2)$ are seen as vertices $V_i$ of a graph. Let $V$ be the standard $2$-dimensional representation. The number of arrows
$V_i\to V_j$ is by definition $\dim_G\mathrm{Hom}(V_i, V\otimes V_j)$.

Get the quivers of the extended Dynkin diagrams exactly as in Figure \ref{fig:extendedDynkin}. For details see \cite[\S 8.3]{Kir}.

\begin{figure}[H]
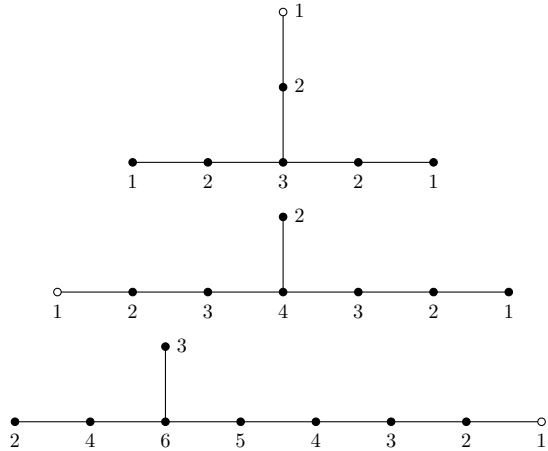

	\centering
\dynkin [extended, labels={1,1,2,2,3,2,1}, edge length=1cm]E6

\dynkin [extended, labels={1,2,2,3,4,3,2,1}, edge length=1cm]E7

\dynkin [extended, labels={1,2,3,4,6,5,4,3,2}, edge length=1cm]E8
\caption{Extended Dynkin diagrams correponding respectively to tetrahedron, according to McKay correspondence.}
\label{fig:extendedDynkin}
\end{figure}
\end{rem0}

\section{Covariants for the variety of polynomials which are powers. A short paper by Hilbert. }
\label{sec:Hilbert1886}
During the school, participants were divided into groups, discussing and trying to understand the main points
of the paper \cite{Hilb}, published by Hilbert in Math. Ann. in 1886. This section is a interlude in the lectures and reflects the chronological order followed at the School.

 Hilbert was born in  1862 and received his doctorate in Koenigsberg, now Kaliningrad, in 1885. So it is a youthful paper, only four pages long.
The question raised in the paper is the following

{\bf Question} {\it Given a homogeneous polynomial $f$ of degree $n=\mu\nu$, what are the conditions for the existence of a homogeneous polynomial $g$ of degree $\nu$ such that }
$$f=g^\mu.$$
\vskip 0.3cm

Hilbert gives a complete answer in the case of binary forms, by constructing a covariant vanishing exactly
on the above variety of $\mu$-powers. 
It is rather strange that this paper by Hilbert was forgotten for more than one century, likely because it was not included in the valuable book \cite{Hilinv}, collecting the English translation
of most of Hilbert's Invariant Theory papers. Rota in 1999 pointed out
the importance of Hilbert result in the colloquium lecture \cite{Ro}.
Let's report Rota's own words

{\it For a long time, I thought the answer to this question to 
be beyond reach, until one day, while leafing through the 
second volume of Hilbert's collected papers, I accidentally 
discovered that Hilbert had completely solved it. The solution can be elegantly
 expressed in umbral notation. This 
is only one of several striking results of Hilbert's on invariant theory that have been forgotten. }

In 2014 Abdesselam and Chipalkatti published the beautiful paper \cite{AC}, where they constructed the same
Hilbert covariant, in a completely different way. They even extend their construction in arbitrary dimension.

Let $d=\mu\nu$. Define
$$X_{d,\mu,n} := \{g^\mu|g\in\sym^{\nu}\CC^{n+1}\},\qquad\textrm{the affine variety of\ }\mu
\textrm{-powers}.$$

\begin{prop0}
$X_{d,\mu,n} $ is smooth (as a projective variety, meaning the affine variety is singular only at the origin), obtained as the cone of a projection of $v_\mu(\PP\sym^{\nu}\CC^{n+1})$.
\vskip 0.4cm
 $T_{g^\mu}X_{d,\mu,n}  = \{ g^{\mu-1}h|h\in \sym^{\nu}\CC^{n+1}\}$.
\end{prop0}
\begin{proof} The tangent space is computed by expanding $(g+\epsilon h)^\mu=g^\mu+\epsilon\mu g^{\mu-1}h+\ldots$. Since all tangent spaces have the same dimension,
the variety is smooth.
\end{proof}
\vskip 0.8cm

We want to sketch the first interesting cases of Hilbert results in \cite{Hilb}. We construct equations
for the varieties of quartics which are square and sextics which are cubes.
Let's recall the following characterization of powers of linear forms, giving equations for the rational normal curve.
\begin{prop0}
Given $f\in\sym^d\CC^2$, 
$$\exists\textrm{\ a linear form\ }\ell\textrm{\ such that\ }f=\ell^d
\Longleftrightarrow Hess(f)=0$$
\end{prop0}
\begin{proof}
$\Longrightarrow$ Trivial.

$\Longleftarrow$ The hypothesis amounts to the vanishing of the Jacobian of the map $\PP^1\rig{(f_x, f_y)}\PP^1$.
Hence the map drops dimension and its image is reduced to a point.
Then there exists scalars $(\lambda,\mu)\neq (0,0)$ such that $\lambda f_x+\mu f_y=0$
which implies $f_{\lambda x+\mu y} =0$. Hence $f$ depends on only one variable, after a linear coordinate change.
The result is proved.
\end{proof}

The starting example is given by binary quartic which are squares. This is the Veronese surface in $\PP^4$,
to distinguish it from her famous parent, which is the Veronese surface in $\PP^5$, embedded with the complete linear system $|{\mathcal O}_{\PP^2}(2)|$. The embedding in $\PP^4$
is given by a codimension one subspace of $H^0({\mathcal O}_{\PP^2}(2))$. Geometrically,
it is given by the projection from a general point $p\in\PP^5$.

\begin{lemma0}\label{lem:rootfH}
Assume that a root of $f$ is also a root of the Hessian $H$. Then it is at least a double root of $f$.
\end{lemma0}
\begin{proof}
We may assume $f=xg$, then by assumption $x$ divides $H$.
Compute

 $f_x=g+xg_x$, $f_y=xf_y$

$f_{xx}=2g_x+xg_{xx}$, $f_{xy}=g_y+xg_{xy}$, $f_{yy}=xg_{yy}$

Hence if $x$ divides $H$ it follows that $x$ divides $g_y$.
By Euler relation, this implies that $x$ divides $g$, hence $x^2$ divides $f$, as we wanted.
\end{proof} 

The covariant $T$, introduced in next Proposition \ref{prop:verp4}, is crucial
for this section, see Theorem \ref{thm:invT}.

\begin{prop0}\label{prop:verp4}
The Veronese surface of squares in $\PP^4$, which is the projective variety associated to $X_{4, 2, 1}$ has equations given by the covariant $$T=\det\begin{pmatrix}f_x&f_y\\
H_x&H_y\end{pmatrix},$$ where $f$ is a binary quartic and $H=Hess(f)$.
\end{prop0}
\begin{proof}
The Veronese surface on $\PP^4$ is the closure of the orbit of $f=(xy)^2$. We compute
$H(f)=-12x^2y^2$, then $T(f)=0$.
Conversely, assume $T(f)=0$. $T(f)$ is bihomogeneous of bidegree $(6, 3)$, meaning it has degree $6$ in $(x, y)$ and degree $3$ in $a_i$.
The assumption $T(f)=0$ means that expanding $T=\sum_{i=0}^6x^{6-i}y^iT^{(i)}(a_0,\ldots, a_4)$, all seven polynomials $T^{(i)}$ vanish.
The determinant defining $T$ is the jacobian of the map
$$\PP^1\rig{(f,H)}\PP^1$$
which vanishes if and only the map drops dimension if and only if there exists a homogeneous polynomial $G$
such that $G(f,H)=0$. We have the relation
$$\sum_{i=0}^ka_if^iH^{k-i}=0$$
and we may assume $a_0a_k\neq 0$ (otherwise we can divide by convenient powers of $f$ and $H$,
reducing the degree of $G$).
The relation implies that every root of $f$ is a root of $H$ and conversely. What we argued up to now works for $f$ of arbitrary degree (the map $(f, H)$ from $\PP^1$ to $\PP^1$ has to be modified as $(f^{2d-4}, H^d)$ 
but the conclusion that $f$ and $H$ share the same roots continues to hold).
Only now we use that $\deg f=4$.
By Lemma \ref{lem:rootfH} it follows that $f$ has no simple roots, so it is a square. This concludes the proof.
\end{proof}

\begin{prop0}\label{prop:sexticube}
The surface of cubes in $\PP^6$, which is the projective variety associated to $X_{6, 3, 1}$ has equations given by the covariant $$T=\det\begin{pmatrix}f_x&f_y\\
H_x&H_y\end{pmatrix},$$ where $f$ is a binary sextic and $H=Hess(f)$.
\end{prop0}

\begin{proof} The surface of cubes is the closure of the orbit of $f=(xy)^3$.  we compute $H(f)=-45(xy)^4$ and $T(f)=0$. Conversely, assume $T(f)=0$. We want to prove that $f$ is a cube. Exactly as in the proof of Prop. \ref{prop:verp4} we see that $f$ has no simple roots, so it may be factored as
$f=(g_1g_2g_3)^2$ or $f=(h_1h_2)^3$ where $\deg g_i=\deg h_j=1$. 
We have to exclude that $f=(g_1g_2g_3)^2$. Assume by contradiction  that $f=(g_1g_2g_3)^2$. 

Up to $SL(2)$-action we have three cases, $f=\left\{\begin{array}{c}(x^3+y^3)^2\\(x^2y)^2\\(x^3)^2
\end{array}\right.$
In the first two cases we compute $T(f)\neq 0$, while in the third case $f$ is a cube and it is admissible.
This concludes the proof.

\end{proof}

These examples should illustrate that even the simplest case of $d$ even and binary forms which are a power of a quadratic form is not trivial.
Indeed Hilbert construction generalizes both Propositions \ref{prop:verp4} and \ref{prop:sexticube} and shows

\begin{thm0}\label{thm:invT} Let $d$ be even.
The surface of $d/2$-powers in $\PP^d$, which is the projective variety associated to $X_{d, d/2, 1}$ has equations given by the covariant $$T=\det\begin{pmatrix}f_x&f_y\\
H_x&H_y\end{pmatrix},$$ where $f$ is a binary form of degree $d$ and $H=Hess(f)$.
\end{thm0}

Theorem \ref{thm:invT} is just a particular case (for $\nu=2$) of next Theorem \ref{thm:hilbertmain}.
In order to explain Hilbert technique, we have to recall two differential operators acting on covariants.

We recall from \cite{O12} that a covariant for a binary form $f=\sum_{i=0}^da_i{d\choose i}x^{d-i}y^i$  is a bihomogeneous polynomial $F(a,x)
=\sum_{k=0}^q F_k(a){q\choose k}x^{q-k}y^k$,
which has bidegree $(m,q)$. Then the weight of $\phi_k$ is $\frac{1}{2}(dm-q)+k$, so that the weight
of the minimal weight $\phi_0$ is $\frac{1}{2}(dm-q)$ and the weight of the maximal weight $\phi_q$ is $\frac{1}{2}(dm+q)$.

We have the two differential operators, that may be applied to any $F(a, x)$ as above.

$D=\sum_{i=0}^{d-1}(i+1)a_i\frac{\partial}{\partial a_{i+1}}$

$\Delta=\sum_{i=0}^{d-1}(d-i)a_{i+1}\frac{\partial}{\partial a_{i}}$

$D$ and $\Delta$ leave the degree unchanged. $D$ decreases the weight by $1$, while $\Delta$ increases the weight by $1$. The operator $D$ corresponds to $\begin{pmatrix}0&1\\0&0\end{pmatrix}\in{\mathfrak sl}(2)$,
while $\Delta$ corresponds to  $\begin{pmatrix}0&0\\1&0\end{pmatrix}\in{\mathfrak sl}(2)$,
see \cite[\S 3.3]{O12}.

We have (see \cite[Prop. 11]{O12}

\begin{prop0}

(i) $D\Delta-\Delta D=\sum_{i=0}^{d}(d-2i)a_{i}\frac{\partial}{\partial a_i}.$

(ii) Let $a_0^{\nu_0}\ldots a_d^{\nu_d}$ be a monomial with degree $\sum_{i=0}^d\nu_i=m$ and weight 
$\sum_{i=0}^di\nu_i=p$. Then

 $(D\Delta-\Delta D)(a_0^{\nu_0}\ldots a_d^{\nu_d})=\sum_{i=0}^d(d-2i)\nu_i(a_0^{\nu_0}\ldots a_d^{\nu_d})=(dm-2p)(a_0^{\nu_0}\ldots a_d^{\nu_d}).$
\end{prop0}

Consider now a covariant $F(a,x)$ of bidegree $(m,q)$. When $q=0$ it is an invariant. Its minimal weight $p$ satisfies
$\frac{1}{2}(dm-q)=p$ and its maximal weight satisfies $\frac{1}{2}(dm+q)=p+q$.

A minimal weight for a   covariant $F(a,x)$ of bidegree $(m,q)$ is a isobaric polynomial $F_0(a)$
of degree $m$ and weight $\frac{1}{2}(dm-q)$ which satisfies one of the two equivalent conditions
$DF_0=0$ or $\Delta^{q+1}F_0=0$.

In the same way, a maximal weight for a   covariant $F(a,x)$ of bidegree $(m,q)$ is a isobaric polynomial $F_q(a)$
of degree $m$ and weight $\frac{1}{2}(dm+q)$ which satisfies one of the two equivalent conditions
$D^{q+1}F_q=0$ or $\Delta F_q=0$.

Indeed the weights of the covariant $F$ are segments centered around $(\frac{dm}{2}, \frac{dm}{2})$. 

The basic property we use is Leibniz rule for $\Delta$ (it holds for $D$ as well)

$$\Delta(fg)=f\Delta(g)+\Delta(f)g$$ which implies in turn

$$\Delta(f^k)=kf^{k-1}\Delta(f)$$

\begin{prop0}\label{prop:d+1} Let $f=\sum_{i=0}^da_i{d\choose i}x^{d-i}y^i$.
Then  $\Delta(f)$ is divisible by $x$ and the ratio $(\frac 1x\Delta(f))$ 
coincides with the partial derivative $\partial_y f$.

In particular $\Delta^{d+1}(f)=0$.
\end{prop0}

General formula by Hilbert for the covariant is
\begin{equation}\label{eq:Hilbcov}f^{\nu+1-1/\mu}\Delta^{\nu+1}(f^{1/\mu})\end{equation}

Note that when $f=g^\mu$ we have that $\Delta^{\nu+1}(f^{1/\mu})=\Delta^{\nu+1}g$ and this last expression vanishes by Prop. \ref{prop:d+1}. The rational exponent $1/\mu$ has to be understood formally,and the whole expression can be expanded to an equivalent expression with integer exponents,
 see next Example \ref{exa:hilbsextic}. The expression (\ref{eq:Hilbcov}) is a typical formula of Classical Invariant Theory. A compact and elegant
formulation that magically produces a cascade of nontrivial polynomials and scalars. 

It  has degree $\nu+1$ in $a_i$, which is clear from the definition.
It can be shown it has degree $(d-2)(\nu+1)$ in $x, y$, this is not trivial, we will check this fact in next examples
and explain it when we sketch the proof of Theorem \ref{thm:hilbertmain}.

Hilbert main Theorem in \cite{Hil86} is the following

\begin{thm0}[Hilbert]\label{thm:hilbertmain}
Given a binary form $f$ of degree $n=\mu\nu$, 
  there exists  a homogeneous polynomial $g$ of degree $\nu$ such that 
$$f=g^\mu$$ if and only if $$f^{\nu+1-1/\mu}\Delta^{\nu+1}(f^{1/\mu})=0.$$
\end{thm0}

\begin{exa}\label{exa:hilbsextic}Let see the case $\nu=2$, $\mu=3$, so the condition that a binary sextic is a cube.

The condition is
 
$$f^{8/3}\Delta^3(f^{1/3})$$
We compute

$$\Delta(f^{1/3})=\frac 13 f^{-2/3}\Delta f$$

$$\Delta^2(f^{1/3})=-\frac 29 f^{-5/3}(\Delta f)^2+\frac 13f^{-2/3}\Delta^2f$$

\begin{align*}\Delta^3(f^{1/3})=\frac {10}{27} f^{-8/3}(\Delta f)^3-\frac 49f^{-5/3}\Delta f\Delta^2f
-\frac 29f^{-5/3}\Delta f\Delta^2f+\frac {1}{3} f^{-2/3}\Delta^3 f =\\
\frac 13f^{-8/3}\cdot \left(\frac {10}{9} (\Delta f)^3-2f\Delta f\Delta^2f+f^{2}\Delta^3 f \right)\end{align*}

This computation shows we have to multiply by $f^{8/3}$ to get a polynomial. For binary sextics Hilbert covariant reads (up to scalar multiples) as
 $$\left(\frac {10}{9} (\Delta f)^3-2f\Delta f\Delta^2f+f^{2}\Delta^3 f \right)$$

It is a polynomial of degree $18$ in $(x, y)$ and degree $3$ in $a_i$.
Actually this polynomial is divisible by $x^6$, after performing such a division
it becomes a polynomial of degree $12$ in $(x, y)$, which is proportional to the covariant $T$ found in Prop. \ref{prop:sexticube}. To check this claim it is enough to check that the coefficient of $x^{18}$ is 
proportional to 
$2a_1^3-3a_0a_1a_2+a_0^2a_3$ as in (\ref{eq:minT}).
\end{exa}

In general ($\nu=2$, any $\mu$), the analogous computations are

$$\Delta(f^{1/\mu})=\frac {1}{\mu} f^{\frac{1-\mu}{\mu}}\Delta f$$

$$\Delta^2(f^{1/\mu})=\frac {1-\mu}{\mu^2}  f^{\frac{1-2\mu}{\mu}}(\Delta f)^2+\frac {1}{\mu} f^{\frac{1-\mu}{\mu}}\Delta^2 f$$

{\footnotesize
\begin{align*}\Delta^3(f^{1/\mu}) & =\frac {(1-\mu)(1-2\mu)}{\mu^3}  f^{\frac{1-3\mu}{\mu}}(\Delta f)^3+\frac {2(1-\mu)}{\mu^2} f^{\frac{1-2\mu}{\mu}}(\Delta f)(\Delta^2 f)
+\frac {(1-\mu)}{\mu^2} f^{\frac{1-2\mu}{\mu}}(\Delta f)(\Delta^2f)+\frac {1}{\mu} f^{\frac{1-\mu}{\mu}}(\Delta^3 f) \\
& =\frac{1}{\mu}f^{\frac{1-3\mu}{\mu}}\cdot \left(\frac {(1-\mu)(1-2\mu)}{\mu^2} (\Delta f)^3+\frac{3(1-\mu)}{\mu}f(\Delta f)(\Delta^2f)+f^{2}(\Delta^3 f) \right)\end{align*}}

 Again we have to multiply by $f^{\frac{3\mu-1}{\mu}}$ to get a polynomial. The expression in parenthesis has degree $3$ in the coefficients $a_i$ of the binary form $f(x,y)$.
Moreover it is divisible by $x^6$, after performing this 
it remains of degree $3d-6$ in the variables $x, y$ and it is now a covariant. \\

\textit{Sketch of proof of Theorem \ref{thm:hilbertmain}:}
Given a polynomial $f\in\sym^{\nu\mu}\CC^2$ we look for conditions such that
there exists $g\in\sym^{\nu}\CC^2$ such that $g^\mu=f$.
We are looking for conditions such that $f^{1/\mu}$ is a polynomial of degree $\nu$.
It is well known that a function $h(t)$ in a coordinate $t$ is a polynomial of degree $\nu$
if and only if $\left(\frac{d}{dt}\right)^{\nu+1}(h(t))=0$.

The idea of Hilbert is pretty clear. He considers formally $f^{1/\mu}$  and imposes the condition
$$\Delta^{\nu+1}(f^{1/\mu})=0$$ to guarantee that  $f^{1/\mu}$ is a polynomial, thanks to
Proposition \ref{prop:d+1}.
In order to form a covariant, such that we get a polynomial,  it turns out that we have to multiply this operator by a convenient power of $f$, as
$f^{\nu+1-1/\mu}\Delta^{\nu+1}(f^{1/\mu})$.

We already explained the details in the particular cases $\left(d=4, \nu=2\right)$, $\left(d=6, \nu=2\right)$.

We just sketch some considerations in the general case, referring to \cite{AC} for a detailed proof,
where they consider a refinement of our operator $\Delta$ involving also $(x, y)$.

The minimal weight of our covariant is 
$$a_0^{\nu+1-1/\mu}\Delta^{\nu+1}a_0^{1/\mu}$$
This has degree and weight $\nu+1$ in the $a_i$. Hence the degree $q$ satisfies
$\frac{1}{2}d(\nu+1)-q=\nu+1$ which gives $q=(d-2)(\nu+1)$.

Note $\Delta(a_0)=da_1$, $\Delta^2(a_0)=d(d-1)a_2$ and $\Delta^i(a_0)=\frac{d!}{(d-i)!}a_i$.

Hence for $\nu=1$ we get $\mu=d$ and the minimal weight $a_0^{2-1/d}\Delta^2(a_0^{1/d})$ is proportional to
$$-a_1^2+a_0a_2$$
which is proportional to the minimal weight of the Hessian.

For $\nu=2$  we get $\mu=d/2$ and that $\Delta^3(a_0^{1/\mu})$ is proportional to
$$(1-d/2)(1-d)d^2a_1^3+3(d/2)(1-d/2) d(d-1)a_0a_1a_2+(d/2)^2(d-1)(d-2)a_0^2a_3$$
which is proportional to

\begin{equation}\label{eq:minT}2a_1^3-3a_0a_1a_2+a_0^2a_3\end{equation}
which in turn is proportional to the minimal weight of
$T=f_xH_y-f_yH_x$.

Note in particular the two vanishing

$D(-a_1^2+a_0a_2)=0$

and

$D(2a_1^3-3a_0a_1a_2+a_0^2a_3)=0.$

These are just particular cases of 
$$D(a_0^{\nu+1-1/\mu}\Delta^{\nu+1}a_0^{1/\mu})=0$$ which is a nontrivial exercise proved by Hilbert in
\cite{Hil86} or in \cite[(III) \S 1.6]{Hilb} or in \cite[Lemma 4.1]{AC}.
\qed
\vskip 0.3cm

\begin{rem0}\label{rem:Hermitereci}
We discuss the degree of the generators of the ideal $I_{X_{d,\mu,n}}$.
Hilbert constructs for $n=1$ elements of the ideal in degree {$\ge d/\mu+1$} and he proves,
as we saw above, 
that the elements of the Hilbert covariant, cut the variety set-theoretically.
Abdesselam and Chipalkatti proved in \cite[Theorem 1.2]{AC} that the Hilbert covariant
cut $X_{d,\mu,n}$ scheme-theoretically.
Recently, Raicu, Sam, Weyman and Yang have proved in \cite{RSWY} that the ideal 
is generated  by the elements constructed by Hilbert, so they describe completely the ideal 
$I_{X_{d,\mu,1}}$.

The generators of the ideal have degree $\ge \nu+1$ for every $n\ge 1$.
Let $X=X_{d,\mu,n} \subset\PP(\sym^{\nu}\CC^{n+1})$, we want to understand why 
$I_X(\nu)=0$. We may assume $n=1$ by restriction (for similar restriction arguments see \cite[Theor. 3.8]{FGOV}). Consider that polynomials of degree $\nu$
on ${\PP(\sym^{d}\CC^2)}$ restrict to polynomials of degree $d$ on $\PP(\sym^{\nu}\CC^{n+1})\simeq X$.
Consider the exact sequence
$$0\rig{}{\mathcal I}_X(\nu)\rig{}\O_{\PP(\sym^{d}\CC^2)}(\nu)\rig{}\O_X(\nu)$$
which at $H^0$-level gives

$$0\rig{} I_X(\nu)\rig{}\sym^{\nu}\sym^d\CC^2\rig{\phi}\sym^d\sym^{\nu}\CC^2$$

In \cite{RSWY} it is proved that $\phi$ is an isomorphism, which implies that $I_X(\nu)=0$.
Hermite reciprocity guarantees that the target and the source of $\phi$ are isomorphic
but it is quite subtle how the map $\phi$ works in detail, which is described in \cite{RSWY}. 
See also \cite{RS}, \cite{AC07}, \cite[\S 1.3]{O12}.
\end{rem0}
\vskip 0.4cm

The following Proposition shows that the Hessian of a power has a simple expression. 
We prefer to state it at the end of this section to avoid interruptions around Hilbert paper, and also because it has an independent interest.
A particular case, when $m=2$, was proved in \cite[Prop. 2.4]{COCDR}, with an argument from Invariant Theory.
\begin{prop0}\label{prop:hesspower}
Let $g=g(x_0,\ldots, x_n)$ be a  homogeneous polynomial of degree $m\ge 2$.
Then
$$Hess(g^k)=\frac{k^{n+1}(km-1)}{m-1}g^{(n+1)(k-1)}Hess(g)$$
If $m=1$ then $Hess(g^k)=0$.
\end{prop0}
\begin{proof}
Consider the following partial derivatives of $g$.

$(g^k)_i=kg^{k-1}g_i\quad i=0,\ldots, n$

$(g^k)_{ij}=kg^{k-1}g_{ij}+k(k-1)g^{k-2}g_ig_j\quad i, j=0,\ldots, n$

Then
$$Hess(g^k)=\det\left( kg^{k-2}(gg_{ij}+(k-1)g_ig_j)_{i, j}\right) =k^{n+1}g^{(n+1)(k-2)}\det\left( (gg_{ij}+(k-1)g_ig_j)_{i, j}\right)$$

The matrix  $(gg_{ij}+(k-1)g_ig_j)_{i, j}$ is the sum of two matrices, where the second one has rank $1$.
Expanding its determinant, only two terms survive and we get the following expression

\begin{equation}\label{eq:auxgk}g^{n+1}Hess(g)-(k-1)g^{n}\det\begin{pmatrix}
0&\nabla g\\
(\nabla g)^t&(g_{ij})\end{pmatrix}.\end{equation}

This trick is sometimes called (in a slightly different form) the "matrix determinant lemma" in the literature.

By Euler identity
$(m-1)g_i=\sum_{j=0}^nx_jg_{ij}$

We perform column operations, summing to the zero column
the $j$-th column multiplied by $(-x_j)$,  getting for the last matrix in (\ref{eq:auxgk})

$\begin{pmatrix}-\frac{1}{m-1}\sum_{j=0}^nx_jg_j&g_0&\ldots& g_n\\
\vdots&\vdots&&\vdots\\
g_i-\frac{1}{m-1}\sum_{j=0}^nx_jg_{ij}&g_{i0}&\ldots&g_{in}\\
\vdots&\vdots&&\vdots\end{pmatrix}$

which, using again Euler identity, simplifies to
$\begin{pmatrix}-\frac{m}{m-1}g&\nabla g\\
0&(g_{ij})\end{pmatrix}$

The column operations leave the determinant unchanged, so putting together we get

$$Hess(g^k)=k^{n+1}g^{(n+1)(k-2)}\left(g^{n+1}Hess(g)+\frac{(k-1)m}{m-1}g^{n+1}Hess(g)\right) =$$

$$=\frac{k^{n+1}(km-1)}{m-1}g^{(n+1)(k-1)}Hess(g)$$
\end{proof}

\begin{rem0}\label{rem:hesgk} F. Russo observed that $g$ and $g^k$ define the same polar map. The Hessian is the Jacobian of the polar map, so it is clear that $Hess(g^k)$ and $Hess(g)$ may differ only by a multiple of $g$.
Hence the focus of Prop. \ref{prop:hesspower} is in the computation of the scalar multiple.
\end{rem0}

\section{Hyperbolic Triangle Groups and Modular Forms}\label{sec:hyptriangle}

It is well known that the upper plane $\H=\left\{w\in\CC| (\mathrm{Im}\ w)>0\right\} $ is biholomorphic to the disk $\Delta=\left\{z\in\CC| |z|<1\right\}$ through the map 

\begin{equation}\label{eq:cayleytrans}\begin{array}{ccl}\H&\to&\Delta\\
w&\mapsto&\frac{w-i}{w+i}=z\end{array}\end{equation}
which is called the Cayley transform.

Automorphisms of $\H$ are defined by $PSL(2,\R)=\left\{w\mapsto\frac{aw+b}{cw+d}| a, b, c, d\in\R, ad-bc=1\right\}$.

Through the Cayley transform, automorphisms of $\Delta$ correspond to automorphisms of $\H$ and
we will identify both of them with $PSL(2,\R)$. In the context of automorphisms of $\Delta$, $PSL(2,\R)$ appears in the form
$\left\{z\mapsto\frac{az+b}{\bar{b}z+\bar{a}}| a, b\in\CC, |a|^2-|b|^2=1\right\}$.

The Cayley transform is moreover an isometry if $\H$ is equipped with the Poincaré metric $\frac{dwd\bar{w}}{(\mathrm{Im}\ w)^2}$
and the disk $\Delta$ with the metric $\frac{4dzd\bar{z}}{(1-|z|^2)^2}$. These metrics have both negative constant curvature, equal to $-1$.
All the automorphisms described above are also isometries for these metrics.

Let $\tau$, $\tau'$ two numbers in the half-plane $\H$. 
Consider the two lattices $\Lambda=\langle 1, \tau\rangle_{\Z}$,  $\Lambda'=\langle 1, \tau'\rangle_{\Z}$, 
It is well known that the two elliptic curves $\CC/\Lambda$ and $\CC/\Lambda'$ are isomorphic if and only if
$\tau'=\frac{a\tau+b}{c\tau+d}$ for some $\begin{pmatrix}a&b\\c&d\end{pmatrix}\in SL(2,\Z)$. Hence the quotient $\H /SL(2,\Z)$
can be interpreted as the moduli space of elliptic curves. This quotient is a subtle object, for the presence of automorphisms
at the special points corresponding to harmonic and equianharmonic curves, which prevents the existence of a universal family of elliptic curves.
Just at level of GIT quotient, this quotient is isomorphic to $\PP^1$. We saw this fact with the moduli space of ternary cubics (Theorem \ref{thm:inv_ternary_cubics})
but it is useful to recall the alternative analytic approach on $\H$.

The group $PSL(2,\Z)$ is generated by $T=\begin{pmatrix}1&1\\0&1\end{pmatrix}$, corresponding to $z\mapsto z+1$ and
$S=\begin{pmatrix}0&-1\\1&0\end{pmatrix}$ corresponding to $z\mapsto -\frac 1z$. 

Recall that the function $q(z)=e^{2\pi iz}$ is a $\Z$-covering $\H\to\Delta^*$. 
For every meromorphic function on $\H$, which is invariant by $T$, namely satisfying $f(z+1)=f(z)$, there exists a meromorphic function
$\tilde f$ defined on $\Delta^*$ such that $f(z)=\tilde f(e^{2\pi iz})$. In other words, the invariance of 
$f$ allows us to write $f$ as a function of $q=e^{2\pi iz}$.
We say that such $f$ is holomorphic at $\infty$ if $\tilde f$ can be extended as a homorphic function to the origin.

Consider $g=\begin{pmatrix}a&b\\c&d\end{pmatrix}\in SL(2,\Z)$ and compute $gz=\frac{az+b}{cz+d}$
\begin{equation}\label{eq:form}\frac{d(gz)}{dz}=(cz+d)^{-2}\end{equation}
This computation motivates the following
\begin{defn0} A modular form of weight $k\in\Z$ is a holomorphic function on $\H$
which satisfies
\begin{equation}\label{eq:modular}f\left(\frac{az+b}{cz+d}\right)=(cz+d)^{2k}f(z)\quad\forall\begin{pmatrix}a&b\\c&d\end{pmatrix}\in SL(2,\Z)\end{equation}
and it is moreover holomorphic at infinity in the above sense. This is meaningful since
(\ref{eq:modular}) is equivalent to $f(z+1)=f(z)$ and $f(\frac 1z)=z^{2k}f(z)$.
\end{defn0}

They are called "forms" since, as well explained by Serre in \cite[VII 2.1]{Ser} the equation (\ref{eq:form}) may be written as
$$\frac{f(gz)}{f(z)}=\left(\frac{d(gz)}{dz}\right)^{-k}$$
or
$$f(gz)\left(d(gz)\right)^k=f(z)\left(dz\right)^k.$$

The basic examples of modular form are the Eisenstein series
$$G_k(z)=\sum_{(m. n)\in\Z^2\setminus(0,0)}\frac{1}{(mz+n)^{2k}}\qquad k\in\Z, k\ge 1$$

We already met the Eisenstein series in the proof of Prop. \ref{prop:G2G3}. $G_k$ is a modular form of weight $k$, since the action of $SL(2,\Z)$ corresponds just to a reordering of the series,
which converges absolutely.

Denote by $M_k$ the vector space of  modular forms of weight $k$.
There are obvious maps $M_a\times M_b\to M_{a+b}$, defining a graded ring structure on
$\oplus_{k\ge 0}M_k$.

The main result, which shows that this graded ring is analogous to the invariant ring for binary quartics  (Theorem \ref{thm:inv_binary_quartics}) or
for ternary cubics  (Theorem \ref{thm:inv_ternary_cubics}), is the following

\begin{thm0}\cite[Coroll. 2 in \S VII 3.2]{Ser}\label{thm:G2G3}
The graded ring $\oplus_{k\ge 0}M_k$ is isomorphic to the polynomial ring $\CC[G_2, G_3]$.
\end{thm0}

Notice that the condition $G_2=0$ corresponds once again to the equianharmonic case while the condition $G_3=0$
corresponds to the harmonic case, see Prop. \ref{prop:G2G3}.
Notice moreover that Theorem \ref{thm:G2G3} has nontrivial consequences on the Eisenstein series, for example
$G_4$ is a scalar multiple of $G_2^2$ and $G_5$ is a scalar multiple of $G_2G_3$.

It is useful to recall the 

\begin{thm0}[{\bf Riemann Uniformization Theorem}]\label{thm:rut}\*

\begin{itemize}

\item{} Every Riemann surface of genus zero is isomorphic to the sphere $\PP^1(\CC)=S^2$, where the Fubini-Study metrix coincides, up to a scalar factor, with the Euclidean metric of $S^2$
with positive constant Gaussian curvature.

\item{} Every Riemann surface of genus one is isomorphic to a quotient $\CC/\Lambda$ as in Section \ref{sec:ternary_cubics} and can be equipped with the metric coming from the flat metric of $\CC$.

\item{} Every Riemann surface of genus $g\ge 2$ is isomorphic to a quotient $\Delta/\Gamma$ where $\Gamma\subset PSL(2,\R)$ is a group of isometries of the Poincaré metric. In particular
the Riemann surface is equipped with a  metric of negative constant Gaussian curvature.
\end{itemize}
\end{thm0}
The disk $\Delta$ can be tesselated in a way completely analogous to what we did on the sphere. The structure of the quotient $\Delta/\Gamma$ 
in Riemann Uniformization Theorem is analogous, at least in the examples over $\bar{\QQ}$, at what we did for the polyhedral groups.
This starts the beautiful story of the hyperbolic triangle groups, that we introduce with the following two pictures

\begin{figure}[H]
	\centering
\hskip 1cm \includegraphics[width=3cm]{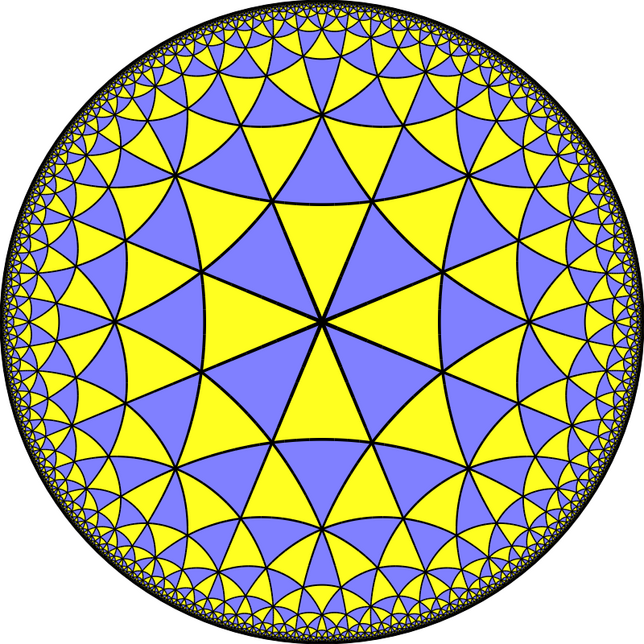}\qquad
\includegraphics[width=3cm]{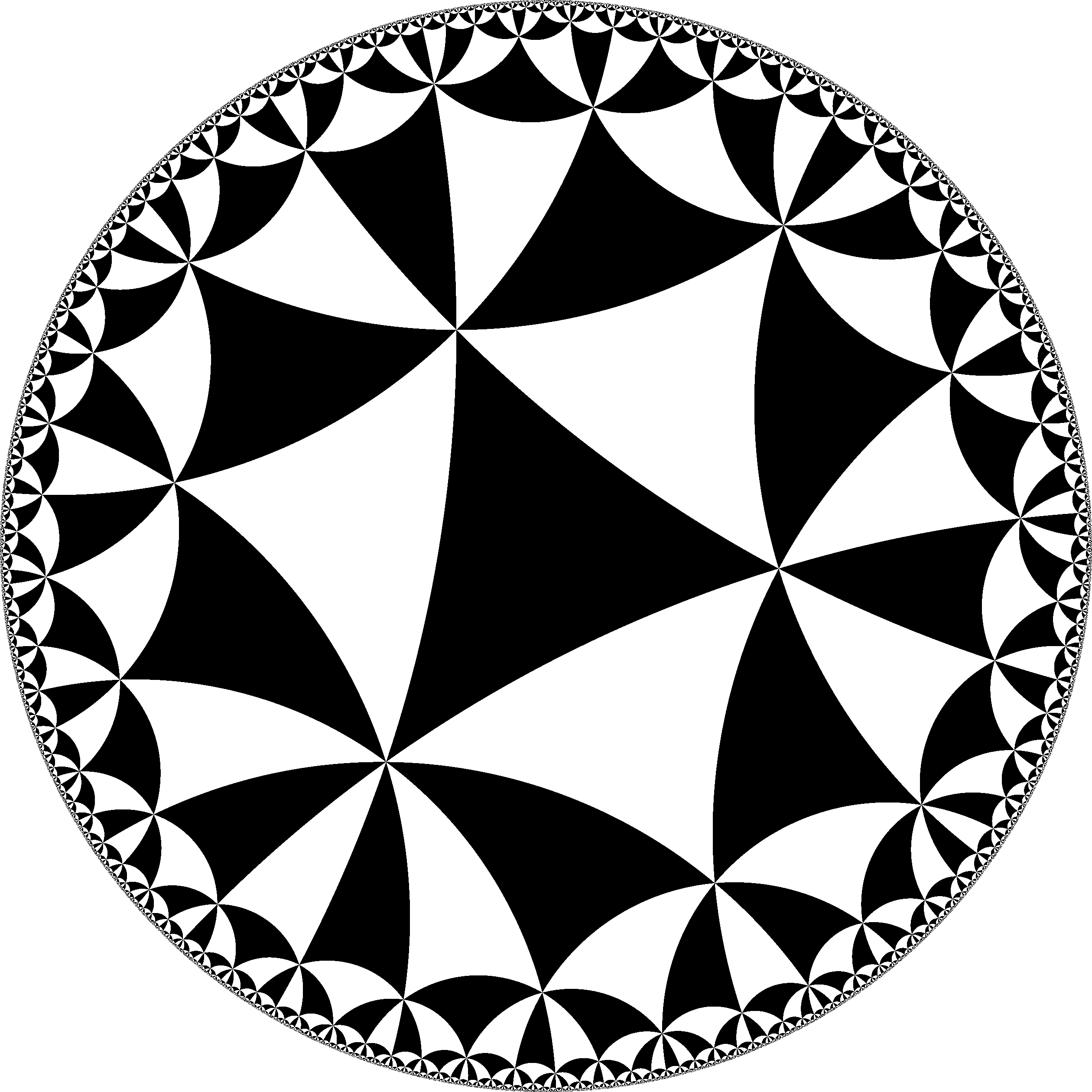}
\caption{Tessellations of the disk with angles respectively $\pi/3$, $\pi/3$, $\pi/4$ (recorded as $(3,3,4)$)
and $\pi/4$, $\pi/4$, $\pi/5$ ($(4,4,5)$).}\label{fig15}
\end{figure}

They can be seen as tessellations of $\Delta$ or, equivalently, as tesselations of $\H$, due to the Cayley transform (\ref{eq:cayleytrans}).

Most of the pictures in this section appeared in a textbook, for the first time, in the influential treatise \cite{FK} by Fricke and Klein, which, along with other publications by the same authors, marked a breakthrough in the portrayal of mathematics. It is impressive that such pictures were drawn without a computer, which is obviously used today. Some original hand-drawn sketches by Grothendieck, found in his archives, featured the same tessellations and are available on the web under "grothendieck artistic drawing".

Consider a geodesic triangle with sides $L_0$, $M_0$, $N_0$, respectively opposite to angles $\pi/l$,  $\pi/m$, $\pi/n$ . This is meaningful on both $\PP^1(\CC)$, $\CC$ and $\H$, where we can consider $L_0$, $M_0$, $N_0$ as whole "lines" (geodesics) where the sides lie.
Let $L$, (resp. $M$, $N$) be the reflections with respect to  $L_0$, (resp. $M_0$, $N_0$). The group  generated by $L$, $M$, $N$, will be called
$T^*(l, m, n)$.

Group theoretically, it is given by \cite[Theor. 2.8]{Magn}.

$$T^*(l, m , n):=\langle L, M, N | L^2=M^2=N^2=1, (LM)^n=(MN)^l=(NL)^m=1\rangle$$

The subgroup of index $2$ in $T^*(l, m , n)$ which consists of even compositions of $L$, $M$, $N$, (they correspond to orientation-preserving maps)
is denoted by $T(l, m, n)$ and it is called a {\it triangle group.}

\begin{prop0}\cite[Theor. 2.10]{Magn} $T(l, m, n)$ is generated by $A=LM$ and $B=MN$ and it is group-theoretically defined as

$$T(l, n, m):=\langle A, B | A^n=B^l=(AB)^m=1\rangle$$

\end{prop0}

Given a tessellation of $\H$ with triangles with angles respectively $\pi/l$, $\pi/m$, $\pi/n$,
such that 
$$\frac 1l +\frac 1m +\frac 1n<1,\quad\textrm{ this is the inequality opposite to (\ref{eq:pqr}),}$$
 we have an action of $T(l, n, m)$ on $\H$.
A basic result in differential geometry guarantees that the area of any triangle has now a sign difference with respect to 
Proposition \ref{prop:areatriang}. Precisely we have
\begin{prop0}\label{prop:areahyp}
The area of a geodesic triangle on $\H$, with angles
$\frac{\pi}{l}$, $\frac{\pi}{m}$, $\frac{\pi}{n}$ is equal to
$ \pi-\frac{\pi}{l} -\frac{\pi}{m}-\frac{\pi}{n}$ . In particular we have the inequality
$$\frac 1l+\frac 1m+\frac 1n<1$$
\end{prop0}

The fundamental domain for $T(l, m, n)$ is given by the union of two adjacent triangles with angles
$\pi/l$, $\pi/m$, $\pi/n$.

We have the following isomorphism

$$\H/T(l, m, n) =\PP^1(\CC)$$

 \begin{figure}[H]
\centering
 \includegraphics[width=6cm]{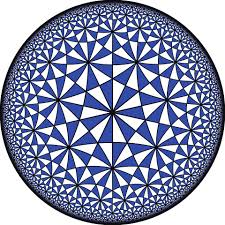} 
\caption{Bisected Hexagonal group, $(2, 3, 7)$. It has a subgroup $\Gamma\subseteq T(2, 3, 7)$ of index $168$ giving 
 the Klein quartic $X(7)$ as the quotient $\H/\Gamma$.}\label{fig16}
 \end{figure}

The following theorem of Belyi, for genus $g\ge 2$, can be seen as an arithmetic refinent of Riemann Uniformization Theorem \ref{thm:rut}.
 \begin{thm0}[Belyi]
 Denote by $\overline{\QQ}$ the field of algebraic numbers. Let $X$ be a Riemann surface of genus $g\ge 2$. The following are equivalent
 \begin{enumerate}
 \item $X$ is defined over $\overline{\QQ}$,
 
 \item $X\simeq \H/\Gamma$ with $\Gamma\subset T(l, m, n)$ of finite index, $l, m, n<\infty$,
 
 \item $X\simeq \left(\H/\Gamma\right)^c$ with $\Gamma\subset PSL(2,\Z)=T(2, 3, \infty)$
 of finite index, where the vertices corresponding to $\infty$ are called cusps and need a compactification (see Remark \ref{rem:yn}).
 \end{enumerate} 
 \end{thm0}
 
The standard formulation of Belyi Theorem  says that 1. is equivalent to get $X$ as ramified cover of $\PP^1$
at the three points  $\{0, 1, \infty\}$.
Actually, the finite map $\H/\Gamma\to \H/T(l, m, n)=\PP^1(\CC)$ is a ramified covering, ramified at most in $3$
 points, corresponding to the vertices of the hyperbolic triangle with angles $\pi/l, \pi/m, \pi/n$. On $\PP^1(\CC)$ the three points  can be assumed to be $\{0, 1, \infty\}$ .
 For details see \cite[Theorem 3 \S 3.3]{Wol}.

This insight allows us to consider many classical examples.
 
The principal congruence group is $$\Gamma(N)=\left\{\begin{pmatrix}a&b\\c&d\end{pmatrix}\in PSL(2,\Z) |
\begin{pmatrix}a&b\\c&s\end{pmatrix}=\begin{pmatrix}1&0\\0&1\end{pmatrix}\mathrm{\ mod\ }N\right\}$$

Two related groups are

 $$\Gamma_1(N)=\left\{\begin{pmatrix}a&b\\c&d\end{pmatrix}\in PSL(2,\Z) |
\begin{pmatrix}a&b\\c&s\end{pmatrix}=\begin{pmatrix}1&*\\0&1\end{pmatrix}\mathrm{\ mod\ }N\right\}$$

 $$\Gamma_0(N)=\left\{\begin{pmatrix}a&b\\c&d\end{pmatrix}\in PSL(2,\Z) |
\begin{pmatrix}a&b\\c&s\end{pmatrix}=\begin{pmatrix}*&*\\0&*\end{pmatrix}\mathrm{\ mod\ }N\right\}$$

Obviously

$$\Gamma(N)\subset\Gamma_1(N)\subset\Gamma_0(N)$$
 
\begin{figure}[H]
  	\centering
 \includegraphics[width=5cm]{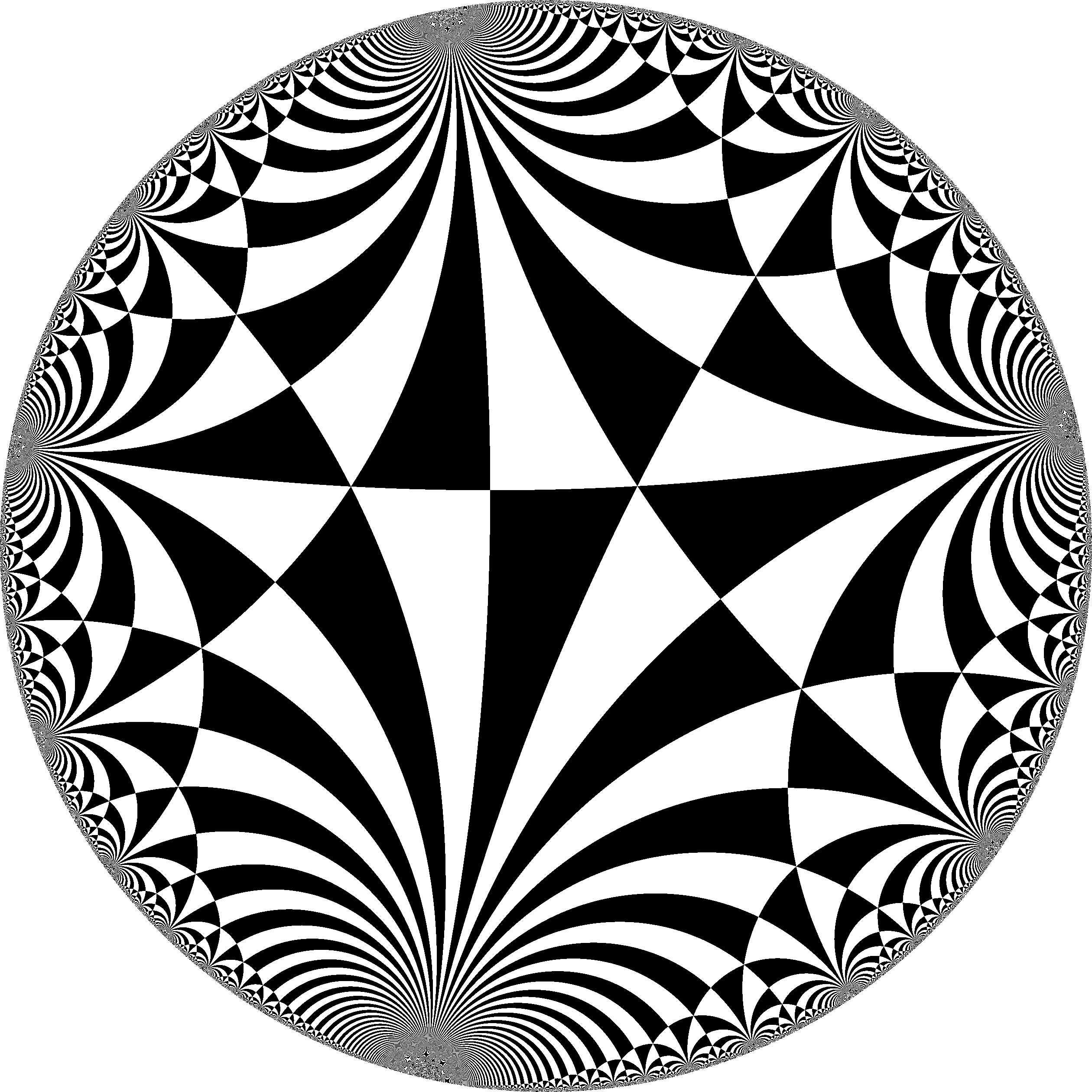} 
 \qquad  \includegraphics[width=5cm]{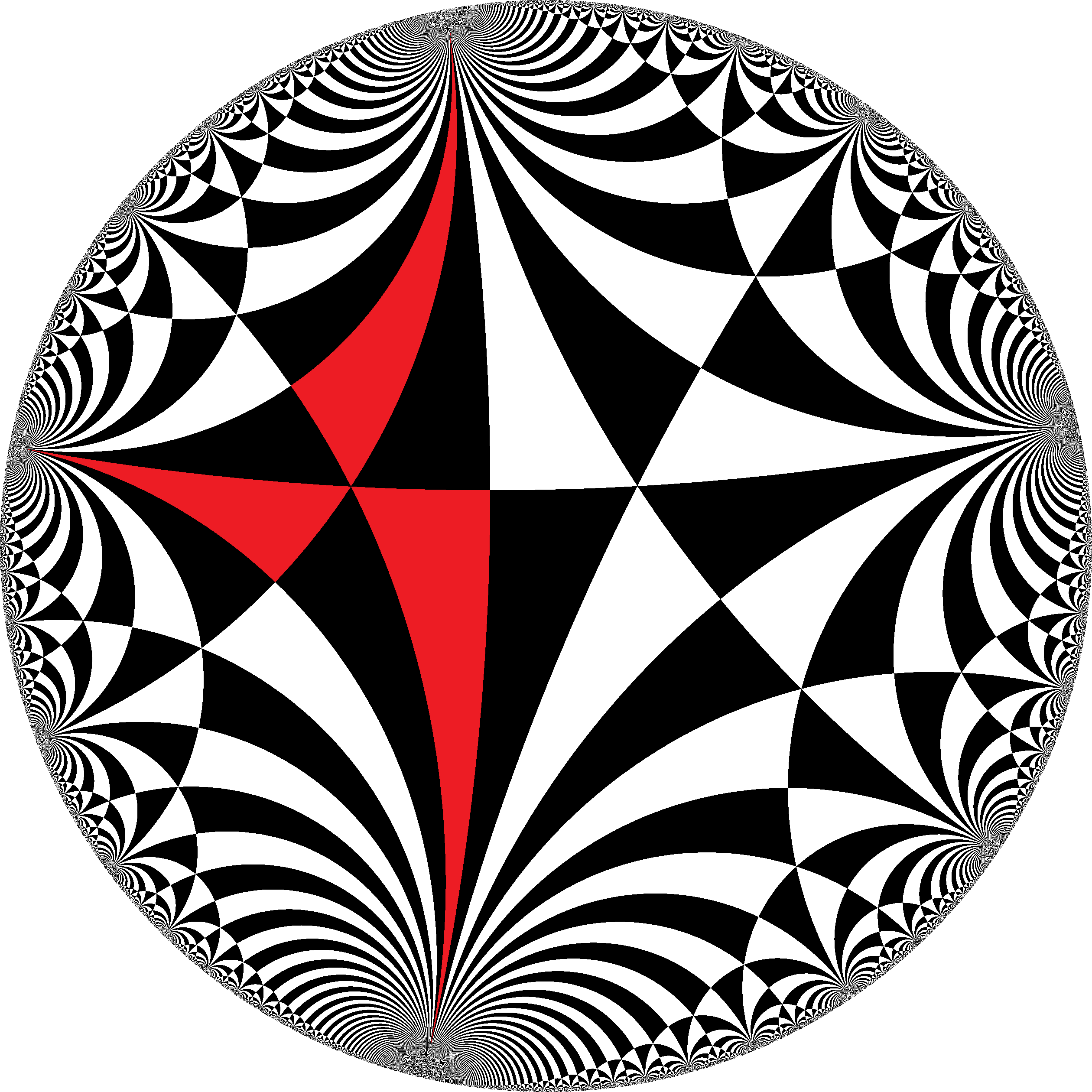} 
  \caption{On the right,
the modular group of level $2$ $\Gamma(2)=T(\infty, \infty, \infty)$,  as a subgroup of index $6$ of $T(2, 3, \infty)=PSL(2,\Z)$, on the left. }\label{fig:Gamma2}
\end{figure}

 \begin{rem0}\label{rem:yn}We denote $Y(N)= PSL(2,\Z) /\Gamma(N)$. This is a noncompact Riemann surface, which can be compactified
 by adding some points corresponding to the boundary points on $\bar{\Delta}$, which can be understood
 by Figure \ref{fig:Gamma2} concerning $\Gamma(2)$. Any geodesic triangle on $\Delta$ with a zero angle (corresponding to $\pi/\infty$) has a vertex on the boundary. The compactified surface is denoted $X(N)$. Note the difference with the previous tesselations, where the triangles had their vertices in the interior of the disk.
 \end{rem0}

 The Figure \ref{fig:Gamma2} shows on the left the fundamental domain (the union of a white and a black adjacent triangles) for $PSL(2,\Z)=\Z_2\ast\Z_3=\{\alpha, \beta | \alpha^2=\beta^3=1\}$, 
 
 By replacing $\Gamma(N)$ with $\Gamma_1(N)$ (resp. $\Gamma_0(N)$) one gets the Riemann surfaces
 $Y_1(N)$ (resp. $Y_0(N)$) and their compactifications  $X_1(N)$ (resp. $X_0(N)$).
 
 Many of these examples appeared already in our story. It can proved that $X(n)$ parametrizes elliptic curves with a fixed basis of their $N$-torsion.
In particular $X(3)$ parametrizes elliptic curves with their $3$-torsion points, that correspond to flexes. So $X(3)$ has genus zero, it corresponds to the parameter $s$ in 
the Hesse pencil (\ref{eq:hessepencil}).
The curve $X(6)$ has genus $1$, it is the Fermat anharmonic cubic,
 while $X(7)$ is isomorphic to the Klein quartic $x^2y+y^2z+z^2x$. This is the quartic with automorphism group of largest size, namely 168. A proof of this bound by the structure of triangle groups can be found in the lovely web page \cite{Baez} by J. Baez.
 
 The genus of $X(N)$ can be computed by an arithmetic formula, and can be found in  A001767
 \cite{AOEIS}, \url{https://oeis.org/A001767}.
 Analogously, the genus of $X_1(N)$ is listed in A029937  \url{https://oeis.org/A029937}. and the genus of $X_0(N)$ is listed in A001617,  \url{https://oeis.org/A001617}.

 \begin{figure}[H]
\centering
 \includegraphics[width=6cm]{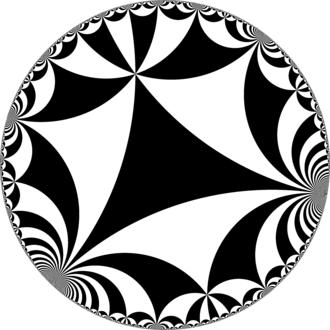} 
\caption{Triangle group, $(5, \infty, \infty)$}\label{fig18}
 \end{figure}

\begin{figure}[H]
\centering
 \includegraphics[width=8cm]{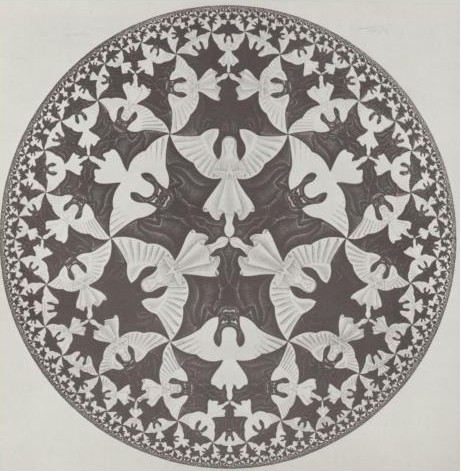} 
\caption{Limit of circle IV(Devils and Angels) by Escher. Here each  half devil joint with an adjacent half angel make a fundamental 
domain for $T(2,6,6)$.}\label{fig:Escher}
 \end{figure}

 We close these lectures with the tesselation of the hyperbolic disk painted by Escher in Figure \ref{fig:Escher}. It reminds us of Weyl's quote  \cite[pag. 500]{Weyl} {\em ``In these days the angel of topology and the devil of abstract algebra
fight for the soul of each individual mathematical domain''}.
 
 \section{Appendix: Pfaffians and Exercises, by V. Galgano}

\subsection{The Pfaffian}\label{appendix:Pfaffians}

We turn our focus on an invariant of skewsymmetric forms: the Pfaffian. Given a complex vector space $V$ of dimension $N$, we consider its skewsymmetric square $\bigwedge^2V$ whose elements are $2$-(skew)forms of the form $\tilde\omega=\sum_{i}v_i\wedge w_i$ as the vectors $v_i,w_i \in V$ vary. Fixing a basis $(e_1,\ldots, e_{N})$ of $V$ allows us to define a linear isomorphism between $\bigwedge^2V$ and the vector space of $N\times N$ skewsymmetric matrices, which to a $2$-form $\tilde\omega=\sum_{i<j}\omega_{ij}e_i\wedge e_j$ associates the skewsymmetric matrix $\omega=(\omega_{ij})$. Due to this, we will often use the exterior square to denote skewsymmetric matrices as well. \\
\indent The natural action of $SL(V)$ on $V$ induces the action on $\bigwedge^2V$ defined by $g \cdot ( v\wedge w ) = (g\cdot v)\wedge (g\cdot w)$, which corresponds to the action by congruence $g \cdot \omega = g\omega g^t$. The Pfaffian is a $SL(V)$-invariant of $\bigwedge^2V$, defined only for even $N=\dim V$. Therefore, we treat separately the cases of $V$ having even or odd dimension.

\paragraph*{Even-dimensional case.} Fix $N=2n$. The $n$-th exterior power $ \tilde\omega^{\wedge n}$ of a $2$-form $\tilde\omega\in \bigwedge^2V$ is a $2n$-form in $\bigwedge^{2n}V \simeq \mathbb C^\vee$. Fixing a basis $(e_1,\ldots, e_{2n})$ of $V$ implies fixing a generator $e_1\wedge \ldots \wedge e_{2n}$ of $\bigwedge^{2n}V$. The {\em Pfaffian} $\Pf(\omega)\in \mathbb C$ of the skewsymmetric matrix $\omega$ is the complex scalar such that
\[ \tilde\omega^{\wedge n} := (n!) \Pf(\omega) e_1\wedge \ldots \wedge e_{2n} \ .\]
\noindent Observe that the coefficient of $e_{i_1}\wedge \ldots \wedge e_{i_{2k}}$ in $\tilde \omega^{\wedge k}$ is the same as in $(\tilde\omega|_{i_1,\ldots, i_{2k}})^{\wedge k}$ since it is obtained by all possible partitions of $\{i_1,\ldots, i_{2k}\}$ into $k$ subsets of two indices each. In particular, the coefficient of $e_{i_1}\wedge \ldots \wedge e_{i_{2k}}$ in $\tilde \omega^{\wedge k}$ is $k!\Pf(\omega|_{i_1,\ldots, i_{2k}})$. This leads to the expression

\begin{equation}\label{eq:lower ext powers and Pfaffians of minors}
	\tilde \omega^{\wedge k} = (k!)\sum_{i_1,\ldots, i_{2k}}\Pf(\omega|_{i_1,\ldots, i_{2k}})e_{i_1}\wedge \ldots \wedge e_{i_{2k}} \ .
\end{equation}

\begin{exa}
	\[ \Pf \begin{pmatrix} 0 & a \\ -a & 0 \end{pmatrix} = a   \ \ \ \ , \ \ \ \ 
	\Pf\begin{pmatrix} 
	0 & a & b & c \\ 
	-a & 0 & d & e \\ 
	-b & -d & 0 & f \\
	-c & -e & -f & 0
	\end{pmatrix} = af - be+ cd \ . \]
\end{exa}

\begin{exe}
	For any $g \in GL(2n)$ it holds $\Pf(g \omega g^t) = \det(g)\Pf(\omega)$. In particular, 
	\[ \Pf \in \big[\sym^n\big(\textstyle{\bigwedge^2V^\vee}\big)\big]^{SL(V)} \ . \]
\end{exe}
\begin{proof}
	\begin{align*}
		& (n!)\Pf(g \omega g^t)e_1\wedge \ldots \wedge e_{2n} = (\widetilde{g\omega g^t})^{\wedge n} = (g \cdot \tilde \omega)^{\wedge n} = g \cdot (\tilde\omega)^{\wedge n} =\\
		& = (n!)\Pf(\omega) (g \cdot e_1\wedge \ldots \wedge e_{2n}) = (n!)\Pf(\omega)\det(g)e_1\wedge \ldots \wedge e_{2n} \ . \qedhere
	\end{align*}
\end{proof}

\begin{rem0}
	The Pfaffian defines a degree-$n$ hypersurface in $\mathbb P(\bigwedge^2V)$: in here, the notion of Pfaffian is preferrable to the notion of determinant as the latter one defines non-reduced schemes.
\end{rem0}

\begin{rem0}
	From a representation-theoretical perspective, the Pfaffian of a generic $2n\times 2n$ skewsymmetric matrix is the unique $SL(2n)$-invariant in the (reducible) $SL(2n)$-module $\sym^n(\bigwedge^2\mathbb C^{2n})^\vee$.
\end{rem0}

Let $J_1:= {\scriptsize \begin{pmatrix} 0 & 1 \\ -1 & 0 \end{pmatrix}}$ and let $J_{n,r}= \diag(J_1, \ldots , J_1, \bold{0},\ldots, \bold{0})$ the $2n\times 2n$ block-diagonal matrix having the first $r$ diagonal blocks equal to $J_1$. 

\begin{exe}
	The following facts holds:
	\begin{enumerate} 
		\item $\Pf(J_{n,n})=1$;
		\item there exists $g \in GL(2n)$ such that $g\omega g^t=J_r$ for a certain $r$;
		\item skewsymmetric matrices have even rank;
		\item $\Pf(\omega)\neq 0$ if and only if $\omega$ has full rank.
	\end{enumerate}
	{\em Hint: For (2) assume $n=2$ and $\omega_{12}\neq 0$, and show that there exists $g$ such that $(g\omega g^t)_{1i}=0$ for any $i\geq 3$.}
\end{exe}
\begin{proof}
	The matrix $J_{n,n}$ corresponds to the $2$-form $\tilde J= \sum_{i=1}^n e_{2i-1}\wedge e_{2i}$ whose $n$-th wedge power is $(n!)e_1\wedge \ldots \wedge e_{2n}$, hence the first statement.\\
	For the other statements, we use an inductive argument. Write
	\[ \omega = \begin{pmatrix}
	\begin{matrix} 0 \ & \omega_{12} \\ -\omega_{12} \ \ \ \ & 0 \end{matrix} & 
	\begin{matrix} \omega_{13} & \ldots & \omega_{1 \, 2n} \\
	\omega_{23} & \ldots & \omega_{2 \, 2n} \end{matrix} \\
	\begin{matrix} -\omega_{13} & -\omega_{23} \\
	\vdots & \vdots \\
	-\omega_{1 \, 2n} & -\omega_{2 \, 2n} \end{matrix} & \omega'
	\end{pmatrix} .\]
	Unless $\omega=0$ (that is $\omega$ has rank at least $1$), we can assume $\omega_{12}\neq 0$. Consider the invertible matrix
	\[ g' := \begin{pmatrix}
	\begin{matrix} 1 \ \ \ \ \ \ \ \ \  & 0 \\ 0 \ \ \ \ \ \ \ \ \ & \omega_{12}^{-1} \end{matrix} & 
	\begin{matrix} 0 & \ldots & 0 \\
	0 & \ldots & 0 \end{matrix} \\
	\begin{matrix} \omega_{23}\omega_{12}^{-1} & -\omega_{13}\omega_{12}^{-1} \\
	\vdots & \vdots \\
	\omega_{2 \, 2n}\omega_{12}^{-1} & -\omega_{1 \, 2n}\omega_{12}^{-1} \end{matrix} & Id_{2n-2}
	\end{pmatrix} \in GL(2n).\] 
	Then 
	\[ g' \omega (g')^t = \begin{pmatrix}
	\begin{matrix} 0 & 1 \\ -1 & 0 \end{matrix} & 
	\begin{matrix} 0 & \ldots & 0 \\
	0 & \ldots & 0 \end{matrix} \\
	\begin{matrix} 0 & \ 0 \\
	\vdots & \ \vdots \\
	0 & \ 0 \end{matrix} & \omega''
	\end{pmatrix} = \diag(J_1, \omega'').\]
	By inductive hypothesis, $\omega''$ has even rank $2(r-1)$ and there exists $g''\in GL(2n-2)$ such that $g''\omega'' (g'')^t= J_{r-1}$. We conclude that $\omega$ has even rank $2r$ and $g:= g'\diag(Id_2,g'')\in GL(2n)$ is such that $g \omega g^t=J_r$. In particular, $\omega$ has full rank if and only if $\Pf(\omega)\neq 0$.
	
\end{proof}

\begin{rem0}
	It is surprising that the canonical form $g\omega g^t=J_r$ holds over any field (even $\mathbb Q$), unlike the symmetric counterpart (the Sylvester canonical form) that holds only over $\mathbb R$. Indeed, the Sylvester form needs square roots, while the canonical form for skewsymmetric matrices only needs unit elements.
\end{rem0}

\begin{exe}
	The following facts are equivalent: 
	\begin{enumerate}
		\item $\omega$ has rank $\geq 2k$;
		\item $\tilde\omega^{\wedge k} \neq 0$;
		\item $\omega$ has a non-singular $2k\times2k$ principal submatrix.
	\end{enumerate}
\end{exe}
\begin{proof}
	The equivalence between $(2)$ and $(3)$ follows from \eqref{eq:lower ext powers and Pfaffians of minors}. Given $\omega$ of rank $2r$, there exists $g \in GL(2n)$ such that $J_{n,r}=g\omega g^t$. In particular, $\widetilde{J_{n,r}}^{\wedge k}= (g \cdot \tilde{\omega})^{\wedge k}= g \cdot \tilde{\omega}^{\wedge k}$, and $\omega^{\wedge k}\neq 0$ if and only if $\widetilde{J_{n,r}}^{\wedge k} \neq 0$, if and only if $r \geq k$. 
\end{proof}

The Pfaffian admits the recursive formula
\begin{equation}\label{eq:recursion Pfaffian}
	\Pf(\omega) = \sum_{j=2}^{2m}(-1)^j\omega_{1j}\Pf(\omega_{\widehat{[1j]}})
\end{equation}
where $\omega_{\widehat{[1j]}}$ denotes the principal sub-matrix of $\omega$ obtained by removing the $1$st and $j$th rows and columns.

\begin{exe}\label{exer:pf-det}
	For every skewsymmetric matrix $\omega$ of even size, it holds 
	\[ \det(\omega)=\Pf(\omega)^2 \ . \]
\end{exe}
\begin{proof}
	If $\omega$ has not full rank, both values are zero. If $\omega$ has full rank, there exists $g \in GL(2n)$ such that $g\omega g^t=J_{n,n}$ and one gets
	\[ \det(g)^2\Pf(\omega)^2 = \Pf(g\omega g^t)^2 = \Pf(J_{n,n})^2 = 1 = \det(g \omega g^t) = \det(g)^2\det(\omega) \ .  \qedhere \] 
\end{proof}

\begin{rem0}
	Historically, Cayley \cite{cayley1849} showed that the determinant of any skewsymmetric matrice (no assumption on the parity of the size) is a square, and defined the Pfaffian precisely as the square root of the determinant. In particular, this relation motivates the fact that the Pfaffian of skewsymmetric matrices of odd size is defined to be zero.
\end{rem0}

\paragraph*{Odd dimensional case.} Fix $N=2n+1$. The Pfaffian of a $N\times N$ skewsymmetric matrix is defined to be zero. Given $\omega \in \bigwedge^2\mathbb C^N$, its {\em adjoint} (or adjugate) is the $N\times N$ matrix $\ad\omega$ such that 
\[ \omega \cdot \ad\omega =\det(\omega)I \ . \]
The matrix $\ad\omega$ is the transpose of the cofactor matrix of $\omega$, that is
\[ \ad\omega = \left( (-1)^{i+j}\det(\Omega_{ji}) \right)_{i,j} \]
where $\Omega_{ij}$ is the submatrix of $\omega$ obtained by erasing the $i$-th row and the $j$-th column.

\begin{exe}
	For $\omega \in \bigwedge^2\mathbb C^{2n+1}$ (odd size), the adjoint $\ad\omega$ is symmetric.
\end{exe}
\begin{proof}
	Since $\omega$ is skewsymmetric, it holds $\Omega_{ji}=-\Omega_{ij}^t$, hence 
	\[ (\ad\omega)_{ji}= (-1)^{i+j}\det(\Omega_{ij})=(-1)^{i+j}\det(\Omega_{ji})=(\ad\omega)_{ij} \ . \qedhere \] 
\end{proof}

\begin{exe}\label{exe:rank of adjoint is 1}
	For every $\omega \in \bigwedge^2\mathbb C^{2n+1}$, it holds $\rk(\ad\omega)\leq 1$, and equality holds if and only if $\rk(\omega)=2n$.
\end{exe}
\begin{proof}
	Note that $\ad\omega =0$ if and only if all $2n\times 2n$ minors of $\omega$ vanish, if and only if $\rk(\omega) \leq 2n-1$. It remains to prove that if $\rk(\omega)=2n$, then $\rk(\ad\omega)=1$. \\
	Assume $\rk(\omega)=2n$, so $\dim \ker(\omega)=1$. Since $\det(\omega)=0$, it holds $\omega \cdot \ad\omega=\bold{0}$, that is the columns of $\ad\omega$ lie in $\ker(\omega)$. Then $\rk(\ad\omega)=1$. 
\end{proof}

\begin{exe}\label{exe:CiCj}
	Let $C_i:= \Pf(\Omega_{ii})$. Then $\det(\Omega_{ij})=C_iC_j$.
\end{exe}
\begin{proof}
	We can assume $\rk(\ad\omega) =1$, otherwise the statement is trivial. We know that $\ad\omega$ is symmetric, hence there exists $\bold{v}=(v_1,\ldots, v_{2n+1})^t \in \mathbb C^{2n+1}$ such that $\ad\omega = \bold{v} \cdot \bold{v}^t$. In particular, we get 
	\[
	v_i^2 =\det(\Omega_{ii})= C_i^2 \ \ \ \ \ , \ \ \ \  
	v_iv_j = (-1)^{i+j}\det(\Omega_{ji}) \ .
	\]
	Then $C_i^2C_j^2 = v_i^2 v_j^2 = \det(\Omega_{ij})^2$, so that $\det(\Omega_{ij})=\pm C_iC_j$. This is a polynomial identity, so it has to be satisfied by any evaluation at a skewsymmetric matrix. In particular, one concludes the thesis by checking for $\tilde\Omega = \diag(J_1,\ldots, J_1,0)$, where the $2\times 2$ block $J_1$ appears $n$ times.
\end{proof}

\begin{exe}
	The vector $(C_1, - C_2, C_3, \ldots, - C_{2n}, C_{2n+1})$ lies in $\Ker(\omega)$.
\end{exe}
\begin{proof}
	Observe that the fact that the principal Pfaffians (i.e. of diagonal submatrices) lower the rank of $\omega$ is a consequence of Exercise \ref{exe:rank of adjoint is 1}.\\
	Now, let $C:=(C_1, - C_2, C_3, \ldots, - C_{2n}, C_{2n+1})$ and $C':= \omega C \in \mathbb C^{2n+1}$. Then:
	\[C'_i = \sum_{j< i}(-1)^{j}\omega_{ji}C_j + \sum_{j>i}(-1)^{j+1}\omega_{ij}C_j \ . \]
	Multiplying by $(-1)^{i+1}C_i$ gives
	\begin{align*} 
		(-1)^{i+1}C_iC_i' & = (-1)^{i+1}\left( \sum_{j< i}(-1)^{j}\omega_{ji}C_iC_j + \sum_{j>i}(-1)^{j+1}\omega_{ij}C_iC_j \right) \\
		& \stackrel{\clubsuit}{=} \sum_{j< i}(-1)^{i+j+1}\omega_{ji}\det(\Omega_{ij}) + \sum_{j>i}(-1)^{i+j+2}\omega_{ij}\det(\Omega_{ij}) = \det(\omega) = 0  
	\end{align*} 
	where the equality $(\clubsuit)$ follows from Exercise \ref{exe:CiCj}.
\end{proof}

\begin{exe}
	Let $D_k:= \{ [\omega] \in \mathbb P(\bigwedge^2\mathbb C^{2n+1}) \ | \ \rk(\omega)\leq 2k \}$ be the determinantal variety defined by the vanishing of all $(2k+2) \times (2k+2)$ subPfaffians of $\omega$. Then
	\begin{enumerate}
		\item the singular locus of $D_k$ is $D_{k-1}$;
		\item $D_k$ has codimension $\binom{2n+1-2k}{2}$.
	\end{enumerate}
\end{exe}
\begin{proof}
	Denote $J_{n,k}^\circ := \diag(J_{n,k}, 0)$ of size $(2n+1)\times (2n+1)$. A generic point $[\omega] \in D_k$ has rank $2k$. Then there exists $g \in GL(2n+1)$ such that $g \omega g^t = J_{n,k}^\circ$. In particular, the orbit $GL(2n+1)\cdot [J_{n,k}^\circ]$ is dense in $D_k$, so they have the same dimension. \\
	\indent The orbit $GL_{2n+1}\cdot [J^\circ_{n,k}]$ is dense, hence its points are smooth in $D_k$. In order to show that $Sing(D_k)$ coincides with $D_{k-1}= D_k \setminus (GL_{2n+1}\cdot [J_{n,k}^\circ])$, it is enough to observe that the critical points of a degree-$m$ Pfaffian (i.e. the points at which the gradient vanishes) are precisely the points at which all sub-Pfaffians of degree $m-1$ vanish, due to the recursive formula for Pfaffians in \eqref{eq:recursion Pfaffian}.\\
	\indent Let's now compute the codimension of $D_k$. Given the map
	\[ \begin{matrix}
	\psi: & GL(2n+1) & \longrightarrow & D_k \\
	& g & \mapsto & gJ_{n,k}^\circ g^t
	\end{matrix} \]
	the fiber at $[J_{n,k}^\circ]$ is 
	$\psi^{-1}([J_{n,k}^\circ]) = \{ h\in GL_{2n+1} \ | \ hJ_{n,k}^\circ h^t = J_{n,k}^\circ\}$. The matrices $h$ in the fiber are determined by the relation
	
	\[ \begin{pmatrix} 
	A & B \\ C & D 
	\end{pmatrix} \begin{pmatrix}
	J_k & 0 \\ 0 & 0
	\end{pmatrix} \begin{pmatrix}
	A^t & C^t \\ B^t & D^t
	\end{pmatrix} = \begin{pmatrix}
	AJ_k & 0 \\ CJ_k & 0
	\end{pmatrix} \begin{pmatrix}
	A^t & C^t \\ B^t & D^t
	\end{pmatrix} = \begin{pmatrix}
	AJ_kA^t & AJ_kC^t \\
	CJ_kA^t & CJ_kC^t 
	\end{pmatrix} \ .
	\]
	This implies that
	\[ \begin{cases}
	AJ_kA^t = J_k\\
	CJ_kC^t = \bold{0}\\
	AJ_kC^t = \bold{0} \\
	\end{cases} \ . \]
	Since $CJ_kC^t=\bold{0}$, then $C$ cannot be invertible. But on the other hand, $AJ_kA^tJ_k^t=I_{2k}$, so $A$ is invertible. Then $AJ_K$ is so, and we get $C=\bold{0}$.\\
	Thus $\psi^{-1}(J^\circ_{n,k})$ is given by matrices of the form {\scriptsize $\begin{pmatrix} A & B \\ 0 & C \end{pmatrix}$} where $A \in Sp(2k)$ is symplectic, $B\in Mat(2k \times (2n+1-2k))$ and $C\in GL(2n+1-2k)$. One deduces the codimension of $D_k$ by applying the fiber dimension theorem. 
\end{proof}

\paragraph*{Pfaffians and $SL(V)$-representations.} We know assume some knowledge on the theory of highest weights for representations of $SL(V)$. Let $V=\mathbb C^N$. We denote by $V_\lambda$ the irreducible representation of $SL(V)$ of highest weight $\lambda \in \mathbb Z^{N-1}$. If $\lambda = (a_1,\ldots , a_{N-1})$, then $V_{\lambda}$ lies in the tensor product (see \cite[\S 15.3]{FH})
\[\textstyle{\sym^{a_1}V \otimes \sym^{a_2}(\bigwedge^2V) \otimes \ldots \otimes \sym^{a_{N-1}}(\bigwedge^{N-1}V)}.\]  Observe that $\sym^{a_1}V=V_{(a_1,0,\ldots, 0)}$ and $\bigwedge^iV=V_{(0,\ldots,1_i,\ldots 0)}$ are irreducible for all $i$, while $\sym^{a_i}(\bigwedge^iV)$ is not irreducible for every $i\geq 2$ and $a_i \geq 2$.

We focus on covariants of degree $2$ of $\bigwedge^2V$. They lie in the $SL(V)$-module $\sym^2(\bigwedge^2V^\vee)$, which decomposes into the irreducible modules
\[ \sym^2(\textstyle{\bigwedge^2V^\vee}) = V_{(0,2,0, \ldots)} \oplus V_{(0,0,0,1,0,\ldots)} = V_{(0,2,0, \ldots)} \oplus \textstyle{\bigwedge^4V^\vee}. \]
The summand $V_{(0,2,0,\ldots)}$ can be understood as kernel of the $SL(V)$-equivariant map
\[\begin{matrix}
\psi: & \sym^2(\textstyle{\bigwedge^2V^\vee}) & \longrightarrow & \bigwedge^4V^\vee \\
& \alpha \cdot \beta & \mapsto & \alpha\wedge \beta
\end{matrix}\]
where $\alpha,\beta$ are $2$-forms in $\bigwedge^2V^\vee$. The kernel is the direct sum of the following two subspaces:
\begin{enumerate}
	\item[(S1)] $\langle (e_i\wedge e_j)(e_h \wedge e_k) \mid \{i,j\}\cap \{h,k\} = \emptyset  \rangle_{\mathbb C}$;
	\item[(S2)] $\big\langle (e_{i_1}\wedge e_{i_2})(e_{i_3}\wedge e_{i_4}) + (e_{i_1}\wedge e_{i_3})(e_{i_2}\wedge e_{i_4}) , (e_{i_1}\wedge e_{i_4})(e_{i_2}\wedge e_{i_3}) + (e_{i_1}\wedge e_{i_3})(e_{i_2}\wedge e_{i_4}) \mid \{i_1,i_2,i_3,i_4\}\in \binom{[N]}{4} \big\rangle_{\mathbb C}$.
\end{enumerate}
The above map $\psi$ is a nice way to see the Pfaffians (of size $4\times 4$) as coordinates of the equivariant projection onto $\bigwedge^4V^\vee$, that is covariants of degree $2$ of $\bigwedge^2V$. Indeed, given a $2$-form $\tilde{\omega}\in \bigwedge^2V^\vee$ its second exterior power is
\[ \tilde{\omega}^{\wedge 2} = \sum_{I \in \binom{[N]}{4}} \Pf(M_I)e_{i_1}\wedge e_{i_2}\wedge e_{i_3}\wedge e_{i_4} \in \textstyle{\bigwedge^4V^\vee}.\]

Note that for $N=4$ the summand $V_{(0,0,0)}=\bigwedge^4V^\vee \simeq \mathbb C$ is one-dimensional, spanned by the unique $SL(4)$-invariant of degree $2$, which is the Pfaffian of a generic $4\times 4$ skewsymmetric matrix. For $N\geq 5$, there are no $SL(N)$-invariants in degree $2$, since $[\sym^2(\bigwedge^2V^\vee)]^{SL(V)}$ does not contain trivial representations. However, the $4\times 4$ Pfaffians span all degree-$2$ covariants of $\bigwedge^2V$.

\subsection{Binary and ternary forms}\label{appendix:binary and ternary}

\paragraph*{Harmonic quadruples.} Let $P,X,Y$ be three collinear points in $\mathbb R^2$ such that $X$ lies on the segment $PY$. For any affine conic $\mathcal{C}$ passing through $X$ and $Y$, let $\ell_1,\ell_2$ be the lines passing through $P$ and tangent to $\mathcal{C}$ at two points $A,B$. Let $Q$ be the intersection point of $AB$ and $XY$. The point $Q$ is said the {\em harmonic mean} of $X$ and $Y$ with respect to $P$, while $AB$ is the {\em polar} with respect to the pole $P$.

\begin{figure}[h!]
	\centering
	\includegraphics[scale=0.8]{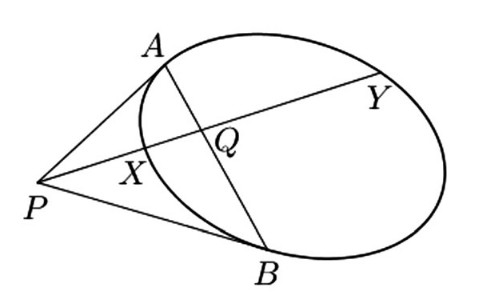}
	\caption{This figure is taken from \cite{DCG}.}\label{fig20}
\end{figure}

\begin{exe}
	The point $Q$ is independent on the conic $\mathcal{C}$. In particular, it only depends on the cross ratio $PX:PY = QX:QY$.
\end{exe}
\begin{proof}
	We can assume $P=[0:0:1]$, $X=[1:0:1]$ and $Y=[\lambda:0:1]$ in $\mathbb P^2_{\mathbb R}$, for $\lambda \gneq 1$. Let $\mathcal{C} \in \sym^2\mathbb R^3$ be a conic defined by a $3\times 3$ symmetric matrix $S$. Note that $s_{33}\neq 0$ since $P \notin \mathcal{C}$. The $\mathcal{C}$ passes through $X$ and $Y$ if and only if $X^tSX = Y^tSY=0$, that is 
	\begin{align*}
		s_{11} + 2s_{13} + s_{33} & = 0 \\
		\lambda^2 s_{11} + 2\lambda s_{13} + s_{33} & = 0 \ . 
	\end{align*}
	Multiplying the first equation by $\lambda^2$ and subtracting the second equation gives 
	\[ 2\lambda(\lambda - 1)s_{13} + (\lambda^2-1)s_{33}=0 \ \implies \ \frac{s_{13}}{s_{33}} = - \frac{\lambda+1}{2\lambda} \ . \]
	The point $Q=[q:0:1]$ is polar to $P$ if $0 = Q^tSP = qs_{13} + s_{33}$, so that $\frac{1}{q} = -\frac{s_{13}}{s_{33}}$. In particular, 
	\[ \frac{1}{q} = \frac{\lambda+1}{2\lambda} \ , \]
	so $q$ (hence $Q$) is independent on the conic. 
\end{proof}

\paragraph*{Binary quartics.} Consider the space of binary quartics $\sym^4\mathbb C^2$ in the variables $x,y$, with the action of $SL(2)$ via linear change of coordinates. A general $f\in \sym^4\mathbb C^2$ is of the form 
\[ f(x,y)=a_4x^4 + 4a_3 x^3y + 6a_2x^2y^2 + 4a_1xy^3 + a_0y^4 \ . \] 
The ring $\mathbb C[\sym^4\mathbb C^2]^{SL(2)}$ of $SL(2)$-invariants of $\sym^4\mathbb C^2$ is generated by the degree-$2$ and degree-$3$ forms
\begin{align*} 
	I(f) & = a_0a_4 - 4a_1a_3 + 3a_2^2 \\
	J(f) & = a_0a_2a_4 + 2a_1a_2a_3 - a_2^3 -a_0a_3^2 - a_4a_1^2 \ .
\end{align*}

\noindent A binary quartic $f$ is {\em harmonic} if $J(f)=0$, and {\em equianharmonic} if $I(f)=0$.\\
The determinant of the hessian matrix of a binary quartic $f \in \sym^4\mathbb C^2$ is the covariant of order $4$ (with respect to $x,y$) and degree $2$ (with respect to $a_0,\ldots, a_4$)
\begin{align}\label{eq:hessian binary quartic}
	H(f) & = f_{xx}f_{yy}-f_{xy}^2\\
	& = (a_4a_2-a_3^2)x^4 + 2(a_4a_1-a_3a_2)x^3y + (a_4a_0 + 2a_3a_1-3a_2^2)x^2y^2 + \nonumber\\
	& \ \ \ + 2(a_3a_0-a_2a_1)xy^3 + (a_2a_0-a_1^2)y^4 \ . \nonumber
\end{align}

\begin{rem0} The origins of the invariants $I$ and $J$ must be sought in Cayely's work in mid XIX century. In his work on symbolic methods for finding invariants, Cayley introduced certain differential operators called {\em transvectants}. The invariant $I(f)$ is (up to scalar) the {\em 4th transvectant} $( f,f)_4$. The invariant $J(f)$ is obtained similarly. It is the iterated transvectant $(f, (f, f )_2 )_4$. Note that $( f,f)_2=H(f)$. Explicit computations can be found in \cite[\S 6]{olver1999classical}. \\
	Cayley's transvectants are a sort of predecessor of Hilbert's operators $D$ and $\Delta$ \cite[\S I.5]{Hilb}, which provide a tool for generating invariants via derivations.
\end{rem0}

\begin{exe}\label{exercise: generic binary quartic}
	Show that a generic $f\in \sym^4\mathbb C^2$ is projectively equivalent to one of the form 
	\[ F_t=x^4 + y^4 + 6tx^2y^2 \ .  \]
	Moreover, $F_t$ is {\em harmonic} if $t\in \{0, \pm 1\}$, and {\em anharmonic} if $t =\pm \frac{i}{\sqrt{3}}$.
\end{exe}
\begin{proof}
	Dehomogeneising $f$ with repsect to $y$ gives $f(x,1)=a_4x^4 + 4a_3x^3 + 6a_2x^2 + 4a_1x +a_0$. We can assume $a_4=1$ up to rescaling. Write $f(x,1)=\prod_{i\in [4]}(x-\alpha_{i})$. Consider the action of $SL(2)$ on $\mathbb P^1$. Then there exists $g \in SL(2)$ sending the four roots of $f(x,1)$ in the four roots of $F_t(x,1)$. These will depend on a parameter $t$ since the moduli space of four points on $\mathbb P^1$ (also, the moduli space of genus-$0$ curves with four marked points) is isomorphic to $\mathbb C$.\\
	The second part of the statement follows from the fact that invariants of $F_t$ are
	\[ I(F_t) = 1 + 3t^2 \ \ \ \ \ , \ \ \ \ \  J(F_t) = t(1-t^2) \ . \] 
\end{proof}

\begin{exe}\label{exe:sec2 and J invariant}
	The variety $V(J)$ of harmonic forms in $\mathbb P(\sym^4\mathbb C^2)$ coincides with the secant variety $\sigma_2(\nu_4(\mathbb P^1))$ of the degree-$4$ rational normal curve.
\end{exe}
\begin{proof}
	The forms obtained from the quartic $F_t$ for $t \in \{0,\pm 1\}$ can be written as $x^4+y^4$, $\frac{1}{2}(x+y)^4 + \frac{1}{2}(x-y)^4$ and $\frac{1}{2}(x+iy)^4 + \frac{1}{2}(x-iy)^4$. These lie in the same orbit $SL(2)\cdot [x^4+y^4]$, which is dense in $\sigma_2(\nu_4(\mathbb P^1))$. On the other hand, any generic harmonic form can be obtained from $x^4+y^4$ via $SL(2)$-action, and any other harmonic form is a degeneration of this, so $SL(2)\cdot [x^4+y^4]$ is dense in $V(J)$ as well.
\end{proof}

\begin{exe}\label{exercise: I,J,Hess}[also see \cite{CO}, \S 4]
	The following facts hold:
	\begin{enumerate}
		\item $f=g^2$ for $g\in \sym^2\mathbb C^2$ nondegenerate $\implies$ $J(f)=0$, $I(f)\neq 0$ and $H([f]) = [f]$;
		\item $f=\ell^4$ for $\ell\in \mathbb C^2$ $\implies$ $J(f)=I(f)=0$ and $H([f])$ is not defined;
		\item The Hessian map restricted to $V(J)$ is dominant on the variety of squares. In particular, $J(f)=0$ $\iff$ $H(f)$ is a square;
		\item The Hessian map restricted to $V(I)$ is a birational involution. In particular, $I(f)=0$ $\iff$ $H(H([f])) =[f]$. 
	\end{enumerate}
\end{exe}

\begin{proof}\*
	\begin{enumerate}
		\item[\textit{1,2.}] 
		If $f=g^2$ for $g \in \sym^2(\mathbb C^2)$, a direct computation shows that
		\begin{align*} 
			H(g^2) 
			= 12 H(g) g^2 \ .
		\end{align*}
		Since $H(g)=0$ if and only if $g=\ell^2 \in \nu_2(\mathbb P^1)$ for some $\ell \in \mathbb C^2$, at the projective level the Hessian map $H: \mathbb P(\sym^4\mathbb C^2) \dashrightarrow \mathbb P(\sym^4\mathbb C^2)$ is not defined over $\nu_4(\mathbb P^1)$, i.e. on $[f]=[\ell^4]$, while it is the identity on all $[g^2]$ for nondegenerate $g\in \sym^2\mathbb C^2$.
		
		Moreover, any $f=g^2$ is in the $SL(2)$-orbit of either $x^4$ or $(xy)^2$. Both such orbits lie in the secant variety $\sigma_2(\nu_4(\mathbb P^1))$, so in particular in the zero locus of $J$.
		
		Finally, we know that $J(f)=I(f)=0$ if and only if $f$ has a root of multiplicity at least $3$. In particular, $V(I,J)$ is the tangential variety of $\nu_4(\mathbb P^1)$, given by the union of the curve itself and the orbit $SL(2)\cdot x^3y$.
		
		\item[\textit{3.}] If $f$ is harmonic, then we can assume $f=x^4+y^4$ that has Hessian $12x^2y^2$, and being square is a $SL(2)$-invariant property. On the other hand, up to projective equivalence we can assume $f$ to be of the form $F_t$ in Exercise \ref{exercise: generic binary quartic}, whose hessian polynomial is $H(F_t)=12^2 \big( tx^4 + ty^4 + (1 - 3t^2)x^2y^2 \big)$. For $t=0$ such polynomial is $(12xy)^2$ and $F_t$ is the Fermat quartic, which is harmonic. Assuming $t\neq 0$ and imposing $H(F_t)$ to be a square, one gets $t = 0 , \pm 1 , \pm \frac{1}{3}$. The first three cases are the harmonic ones, while the last two cases are such that $F_t$ is a square itself.
		
		\item[\textit{4.}] A generic equianharmonic is projectively equivalent to $F_t$ in Exercise \ref{exercise: generic binary quartic} for $t=\pm \frac{i}{\sqrt{3}}$. A direct computation shows that the Hessian of $F_t$ is $H(F_t)=tx^4 + ty^4 + (1-3t^2)x^2y^2$, and in particular $H(F_{\frac{i}{\sqrt{3}}}) = F_{-\frac{i}{\sqrt{3}}}$. This shows that the Hessian map is an involution on $V(I)$. Conversely, writing the generic form of the pencil $F_t$ and imposing that the Hessian is an involution on it leads to the vanishing of the invariant $I$. \qedhere
		
	\end{enumerate}	
	
\end{proof}

\paragraph*{Ternary cubics.} Consider the space of ternary cubics $\sym^3\mathbb C^3$ in the variables $x,y,z$, with the action of $SL(3)$ via linear change of coordinates. The invariant ring 
\[ \mathbb C[\sym^3\mathbb C^3]^{SL(3)}=\langle \mathcal{A}, \mathcal{T} \rangle\]
is generated by two forms of degree $4$ and $6$, called {\em Aronhold} and {\em Clebsch} invariants respectively.\\
Any generic $g \in \sym^3\mathbb C^3$ can be transformed via $SL(3)$ into the Weierstrass form (cf. \cite[Chap.\ 9]{Dolg})
\[ G= y^2z - x^3 -pxz^2 -qz^3 \ .\]
The cubic $G$ is {\em harmonic} if $q=0$ (equivalently, $\mathcal{T}(G)=0$), and {\em equianharmonic} if $p=0$ (equivalently, $\mathcal{A}(G)=0$). For instance, $y^2=x^3\pm x$ are harmonic, while Aronhold proved that a generic equianharmonic ternary cubic is projectively equivalent to the Fermat cubic $x^3+y^3+z^3$.

\begin{exe}
	The $3$rd secant variety $\sigma_3(\nu_3(\mathbb P^2))$ of the $3$rd Veronese surface is the hypersurface $V(\mathcal{A})$ defined by the Aronhold invariant.
\end{exe}
\begin{proof}
	Similar to the proof of Exercise \ref{exe:sec2 and J invariant}, it is enough to consider the Fermat cubic $x^3+y^3+z^3$. Indeed, its $SL(3)$-orbit is dense both in $\sigma_3(\nu_3(\mathbb P^2))$ and in $V(\mathcal{A})$, since any equianharmonic ternary cubic either is in such orbit or is a degeneration of it.
\end{proof}

\begin{exe}[also cf.\ \cite{CO}, \S 5] The following facts hold:
	\begin{enumerate}
		\item $\mathcal{T}(G)=0$ if and only if $H(H(G))$ is proportional to $G$;
		\item $\mathcal{A}(G)=0$ if and only if $H(G)$ is given by three lines.
	\end{enumerate}
\end{exe}
\begin{proof} 
	The first statement can be checked by using the Weierstrass form of harmonic ternary cubics. The second statement can be deduced by looking at the Fermat cubic. 
\end{proof}

\begin{rem0}
	THe above two exercise highlight another interesting dichotomy between binary quartics and ternary cubics. Indeed, the secant variety of $\nu_4(\mathbb P^1)$ corresponds to harmonic forms in $\sym^4\mathbb C^2$, while the third secant variety of $\nu_3(\mathbb P^2)$ to equianharmonic forms in $\sym^3\mathbb C^3$. On the other hand, the Hessian map for binary quartics $\mathbb P^4 \dashrightarrow \mathbb P^4$ is an involution along the equianharmonic forms, while the Hessian map for ternary cubics $\mathbb P^9 \dashrightarrow \mathbb P^9$ is an involution along the harmonic forms.
\end{rem0}

\paragraph*{The Aronhold invariant as Pfaffian.}  Let $\End_0\mathbb C^3$ denotes the space of traceless $3\times3$ matrices. For any $v^3=v \otimes v \otimes v \in \sym^3\mathbb C^3$, consider the linear map
\[ A_{v^3}: \End_0\mathbb C^3 \longrightarrow \End_0\mathbb C^3 \]
such that 
\[ \big(A_{v^3}(M)\big)(w) = \big( M(v)\wedge v \wedge w \big) v \ \ \ \ \ \forall M \in \End_0\mathbb C^3, \ \forall w \in \mathbb C^3 \ . \]
By linearity, the map $A_{v^3}$ can be extended to $A_f$ for any $f\in \sym^3\mathbb C^3$. Then $A_f$ is described by the $9\times 9$ matrix 
\[ \begin{pmatrix}
0 & hess(f_z) & -hess(f_y)\\
-hess(f_z) & 0 & hess(f_x)\\
hess(f_y) & -hess(f_x) & 0
\end{pmatrix} \]
where $hess(f_x)$ is the hessian matrix of the partial derivative $f_x$ (similarly for $f_y,f_z$).

\begin{exe}[cf.\ \cite{ottaviani2009invariant}, Theorem 2.1]
	Let $B_i$ be the $8\times 8$ principal submatrix obtained by removing the $i$th row and column of $A_f$. Then $\Pf(B_i)= \Pf(B_j)$ for all $i,j$. Moreover, the Aronhold invariant $\mathcal{A}$ coincides with such Pfaffians.
\end{exe}

\section*{Acknowledgements}
We thank the Lie-St\o{}rmer Center, K. Kohn, A. Oneto, P. Santarsiero and E. Turatti, along with C. Reiner, for the organization of the School in the stimulating Artic atmosphere, and all the participants for the interest.
GO thanks P. Aluffi and D. Faenzi for historical and mathematical help on group actions and Fano varieties. Special thanks to  M. Andreatta and R. Pignatelli for explanations about Fano's last Fano, see Remark \ref{rem:lastfano}.
C. Raicu's insight kept us away from a slippery argument in Remark \ref{rem:Hermitereci}.
We owe Remark \ref{rem:hesgk} to F. Russo.
GO is member of Italian GNSAGA of INdAM and he is supported by European Union’s HORIZON–MSCA-2023-DN-JD program under the Horizon Europe (HORIZON) Marie Sk{\l}odowska-Curie Actions, grant agreement 101120296 (TENORS). VG is member of the Italian national group GNSAGA-INdAM.

\begin{flushleft}

{\bf AMS Subject Classification: 13A50, 14N05, 15A72}\\[2ex]

%
Giorgio~OTTAVIANI,\\
Dipartimento di Matematica e Informatica ``U. Dini'', Università di Firenze\\
viale Morgagni 67/A, 50134 Firenze, ITALY\\
e-mail: \texttt{giorgio.ottaviani@unifi.it}\\[2ex]

Vincenzo~GALGANO,\\
Max-Planck Institute of Molecular Cell Biology and Genetics, \\
\& Faculty of Mathematics, Technische Universität Dresden\\
Pfotenhauerstrasse 108, 01307 Dresden, GERMANY \\
e-mail: \texttt{galgano@mpi-cbg.de}\\[2ex]

%
\end{flushleft}

\end{document}